\documentclass[11pt,reqno]{amsart}

\usepackage{amsmath}
\usepackage{amsthm}
\usepackage{mathrsfs}
\usepackage{graphicx}
\usepackage{mathtools}
\usepackage{amsfonts}
\usepackage{amssymb}
\usepackage{hyperref}
\usepackage[margin=1 in]{geometry}
\usepackage{enumerate}
\usepackage[normalem]{ulem}
\usepackage{tikz}
\usepackage{tikz-cd}
\usepackage{xcolor}
\usetikzlibrary{matrix,arrows,positioning,automata}
  \usepackage{cancel}
\usepackage{dsfont}

  \usepackage{scalerel,stackengine}
\stackMath
\newcommand\reallywidehat[1]{%
\savestack{\tmpbox}{\stretchto{%
  \scaleto{%
    \scalerel*[\widthof{\ensuremath{#1}}]{\kern-.6pt\bigwedge\kern-.6pt}%
    {\rule[-\textheight/2]{1ex}{\textheight}}
  }{\textheight}%
}{0.5ex}}%
\stackon[1pt]{#1}{\tmpbox}%
}

\numberwithin{equation}{section}
\newcommand{\N}{\mathbb N}
\newcommand{\Z}{\mathbb Z}

\newcommand{\R}{\mathbb R}
\def\E{\mathbb E}

\newcommand{\linf}{L^{\infty}}
\newcommand{\sF}{\mathcal{F}}
\newcommand{\bP}{\mathbb{P}}
\newcommand{\sL}{\mathcal{L}}

\newcommand{\sG}{\mathcal{G}}
\newcommand{\sA}{\mathcal{A}}
\newcommand{\sP}{\mathcal{P}}

\newcommand{\bx}{\bm{x}}
\newcommand{\bX}{\boldsymbol X}

\def\XXint#1#2#3{{\setbox0=\hbox{$#1{#2#3}{\int}$}
\vcenter{\hbox{$#2#3$}}\kern-.5\wd0}}

\newcommand{\T}{\mathbb{T}}

\numberwithin{equation}{section}
\newtheorem{thm}{Theorem}[section]
\newtheorem{lem}[thm]{Lemma}

\newtheorem{prop}[thm]{Proposition}
\newtheorem{assumption}[thm]{Assumption}
\theoremstyle{definition}

\newtheorem{rmk}[thm]{Remark}

\def\smallnegint{\mathop{\int\mkern-13mu
        \raise.5ex\hbox{${\scriptscriptstyle\diagup}$}}\nolimits}
\def\ds{\displaystyle}

\def\div{\operatorname{div}}
\def\tr{\operatorname{tr}}

\def\bx{{\boldsymbol x}}
\def\by{{\boldsymbol y}}
\def\ssetminus{\,\raise.4ex\hbox{$\scriptstyle\setminus$}\,}

\newcommand{\be}{\begin{equation}}
\newcommand{\ee}{\end{equation}}

\newcommand{\bc}{\begin{case}}
\newcommand{\ec}{\end{cases}}
\newcommand{\bs}{\begin{split}}
\newcommand{\es}{\end{split}}

\newcommand{\norm}[1]{\left\Vert#1\right\Vert}

\newcommand{\bm}[1]{\boldsymbol #1}

\renewcommand{\bar}{\overline}
\renewcommand{\tilde}{\widetilde}
\renewcommand{\hat}{\widehat}
\renewcommand{\ln}{\log}
\renewcommand{\Re}{\mathcal{R}}

\def \be {\begin{equation}}
\def \ee {\end{equation}}

\def \E {\mathbb{E}}

\def \R {\mathbb{R}}

 \renewcommand{\i}{{\mathrm{i}}}
\renewcommand{\tilde}{\widetilde}

\newcommand{\pr}{\mathcal{P}}
\newcommand{\cC}{\mathcal{C}}

\newcommand{\cN}{\mathcal{N}}

\newcommand{\leb}{\text{Leb}}
\newcommand{\bbF}{\mathbb{F}}
\newcommand{\bk}{\bm k}
\newcommand{\trip}[1]{{\left\vert\kern-0.25ex\left\vert\kern-0.25ex\left\vert #1 
    \right\vert\kern-0.25ex\right\vert\kern-0.25ex\right\vert}}
\newcommand{\bY}{\bm Y}

\title{On the optimal rate for the convergence problem in mean field control}

\begin{document}
 \author[S. Daudin]{Samuel Daudin 
\address{(S. Daudin) Universit\'e C\^ote d'Azur, CNRS, Laboratoire J.A. Dieudonn\'e, 06108 Nice, France
}\email{samuel.daudin@univ-cotedazur.fr
}}
 
\author[F. Delarue]{ Fran\c{c}ois Delarue 
\address{(F. Delarue) Universit\'e C\^ote d'Azur, CNRS, Laboratoire J.A. Dieudonn\'e, 06108 Nice, France
}\email{francois.delarue@univ-cotedazur.fr
}}

 \author[J. Jackson]{Joe Jackson \address{(J. Jacskon) Department of Mathematics, The University of Texas at Austin, United States}\email{jjackso1@utexas.edu}}

\thanks{
 S. Daudin and F. Delarue acknowledge the financial support of the European Research Council (ERC) under the European Union’s Horizon 2020 research and innovation programme (ELISA project, Grant agreement No. 101054746). Part of the work was achieved during J. Jackson's visit to 
 Université Côte d'Azur between January and March 2023. His visit was also supported by the ERC Grant No. 101054746. J. Jackson is supported by the NSF under Grant No. DGE1610403. Any opinions, findings and conclusions or recommendations expressed in this material are those of the authors and do not necessarily reflect the views of the NSF}
\maketitle

\begin{abstract}
The goal of this work is to obtain optimal rates for the convergence 
problem in mean field control. Our analysis covers cases where the solutions to the limiting problem may not be unique nor stable. Equivalently the value function of the limiting problem might not be differentiable on the entire space. Our main result is then to derive sharp rates of convergence in two distinct regimes. When the data is sufficiently regular, we obtain rates proportional to $N^{-1/2}$, with $N$ being the number of particles. When the data is merely Lipschitz and semi-concave with respect to the first Wasserstein distance, we obtain rates proportional to $N^{-2/(3d+6)}$. Noticeably,  
the exponent $2/(3d+6)$ is close to $1/d$, which is the optimal rate of convergence for uncontrolled particle systems driven by data with a similar regularity. The key argument in our approach consists in mollifying the value function of the limiting problem in order to 
produce functions that
 are almost classical sub-solutions to the limiting Hamilton-Jacobi equation (which is a PDE set on the space of probability measures). These sub-solutions can be projected onto finite dimensional spaces and then compared with the value functions 
associated with the particle systems.  In the end, this comparison 
is used to prove the most demanding bound in the estimates. The key challenge therein is thus to exhibit an appropriate form of mollification.
We do so by employing sup-convolution within a convenient functional Hilbert space. 
To make the whole easier, we limit ourselves to the periodic setting. 
We also provide some examples to show that our results are sharp up to some extent.
\vspace{4pt}

\vspace{4pt}

\noindent \textsc{AMS Classification (2020)}: {Primary: 49N80, 65C35; Secondary: 49L35.}
\end{abstract}

\setcounter{tocdepth}{1}
\tableofcontents

\section{Introduction}

\subsection{A short review of mean field control and games.} Mean field control theory and its twin, mean field game theory, aim at the asymptotic study of equilibria within large populations of weakly interacting agents. Typically, each agent controls a $d$-dimensional state process which is impacted by a Brownian noise. In mean field control, equilibria are understood in a cooperative sense, while in mean field games, they are understood in a competitive sense. The limiting formulations, which arise as the number of players increases to infinity, are distinct: in the cooperative case, we arrive at an optimal control problem set on the Wasserstein space, while the competitive case leads to a well-known fixed point problem.
We refer to 
 \cite{HuangCainesMalhame1,Huang2006,Lasry2006,LasryLions2,LasryLions,Lionscollege2} for earlier  
contributions and \cite{cardaliaguetporretta-cetraro,CarmonaDelarue_book_I,CarmonaDelarue_book_II,gomes_survey}
for
surveys or monographs. 

For almost twenty years, both theories have made parallel and profound advances. In particular, great progress has been made in understanding the infinite-dimensional partial differential equations which describe the relevant value functions - the value of the optimization problem in the case of control, and the value of the equilibrium in the case of games. We refer to 
 \cite{ben-fre-yam2015,CardaliaguetDelarueLasryLions,gan-swi,Lionsvideo} for some key contributions in this direction. In the case of mean field control (with periodic data), the value function is a map 
 \begin{align*}
     U = U(t,m) : [0,T] \times \pr(\T^d) \to \R, 
 \end{align*}
 where $\pr(\T^d)$ denotes the set of probability measures on the $d$-dimensional torus $\T^d = \R^d / \Z^d$. Roughly speaking, $U(t,m)$ denotes the value of the limiting optimization problem when the continuum of agents is distributed according to $m$ at the initial time $t$. We postpone a definition of $U$ as a value function to Subsection \ref{subsec:formulation}, but we mention already that $U$ is expected to solve a first-order Hamilton-Jacobi equation on the space of probability measures, of the form
 \begin{equation} 
 \label{hjbinf} \tag{HJB($\infty$)} 
 \begin{cases}
  \ds - \partial_t U(t,m)   
     - \int_{\T^d} \Delta_{x} \frac{\delta U}{\delta m}(t,m,x) dm(x)  \vspace{.25cm} \\ \ds
     \hspace{1.5cm} + \int_{\T^d} H\Bigl(x, D_x \frac{\delta U}{\delta m}(t,m,x) \Bigr) dm(x) = \sF(m), \quad (t,m) \in (0,T) \times \pr(\T^d), 
     \vspace{.25cm} \\ \ds
     U(T,m) = \sG(m), \quad m \in \pr(\T^d),
     \end{cases}
 \end{equation}
 for a Hamiltonian $H = H(x,p) : \T^d \times \R^d \to \R$ which is typically regular and convex in the second variable, and costs $\sF,\sG : \pr(\T^d) \to \R$.
We refer to Section \ref{sec:mainresults} for more details including the definition of the linear derivative
  $\frac{\delta U}{\delta m}$. For games, the value of the equilibrium problem is instead expected to solve the master equation, which resembles an infinite-dimensional system of first-order hyperbolic equations on the space of probability measures. In both cases, the study of the solution is subtle, but several regimes are known under which the relevant infinite-dimensional PDE admits a classical solution. For control, this is the case if the coefficients are convex in the measure argument and regular, see
 \cite{CardaliaguetDelarueLasryLions,cha-cri-del_AMS,gan-swi}. For games, the convexity condition has to be replaced by a monotonicity condition, see 
 \cite{Bertucci,CardaliaguetDelarueLasryLions,CarmonaDelarue_book_II,cha-cri-del_AMS,Lionsvideo}. 
 
 The importance of the regularity of the value function was explained in the book  \cite{CardaliaguetDelarueLasryLions} by Cardaliaguet, Delarue, Lasry and Lions: the existence of a regular value makes it possible to obtain an optimal rate for the convergence of the values of finite-player games towards the value of their mean field counterparts. In this approach, the bounds on the regularity of the value play an essential role, and the resulting convergence rate is linear in the number $N$ of agents in the finite system. The approach used for games has subsequently been extended to mean field control when the solution of the Hamilton-Jacobi equation is regular, with a rate of the same order, see \cite{germain_pham_warin_2022}. In particular, in the case of control, bounds on the second-order `Lions derivative' $D^2_{mm} U$ allow one to conclude that 
 \begin{align} \label{smoothcase}
     |U(t,m_{\bx}^N) - V^N(t,\bx) | \leq C/N, 
 \end{align}
 where $V^N : [0,T] \times (\T^d)^N \to \R$ denotes the value function for the corresponding $N$-particle control problem, and is, under mild assumptions on the data, the unique classical solution of the Hamilton-Jacobi-Bellman equation  
 \begin{equation} \ds
\label{hjbn} \tag{HJB($N$)}
 \begin{cases} \ds
  - \partial_t V^N(t,\bx) - \sum_{i = 1}^N \Delta_{x^i} V^N(t,\bx)
  + \frac{1}{N} \sum_{i = 1}^N H\bigl(x^i,N D_{x^i} V^N(t,\bx)\bigr) = \sF(m_{\bx}^N), 
  \\
\hspace{240pt}  \quad (t,\bx) \in (0,T) \times (\T^d)^N, 
\vspace{.25cm}
\\
\ds V^N(T,x) = \sG(m_{\bx}^N), \quad \bx \in (\T^d)^N.
 \end{cases}
 \end{equation}
 We note that here and throughout the paper we use the notation $m_{\bx}^N = \frac{1}{N} \sum_{i = 1}^N \delta_{x^i}$ when $\bx = (x^1,...,x^N) \in (\T^d)^N$. The argument leading to \eqref{smoothcase} is relatively simple - if $U$ is smooth, then explicit computation shows that $[0,T] \times (\T^d) \ni (t,\bx) \mapsto U(t,m_{\bx}^N)$ is a solution of \eqref{hjbn} up to an error term which is of order $1/N$ provided that $D_{mm} U$ is bounded, which, by the comparison principle, gives \eqref{smoothcase}. 
 We note that here and throughout the paper we use the notation $m_{\bx}^N = \frac{1}{N} \sum_{i = 1}^N \delta_{x^i}$ when $\bx = (x^1,...,x^N) \in (\T^d)^N$. The argument leading to \eqref{smoothcase} is relatively simple - if $U$ is smooth, then explicit computation shows that $[0,T] \times (\T^d) \ni (t,\bx) \mapsto U(t,m_{\bx}^N)$ is a solution of \eqref{hjbn} up to an error term which is of order $1/N$ provided that $D^2_{mm} U$ is bounded, which, by the comparison principle, gives \eqref{smoothcase}. 

 Thanks to the contributions discussed above, the convergence problem is now well-understood when the relevant infinite-dimensional PDE has a smooth solution, and the existence of smooth solutions in turn is well-understood under certain (fairly restrictive) convexity or monotonicity assumptions. Answering similar questions in the absence of structural conditions like convexity and monotonicity is now one of the main objectives of the theory of mean field games and mean field control. This issue is somewhat easier to understand for control problems than games, simply because it is easier to identify and characterize an optimizer than a fixed point. In fact, several recent works have been published on the Hamilton-Jacobi equation \eqref{hjbinf} in the absence of classical solutions. Most of them aim to understand viscosity solutions and, in particular, to obtain a comparison principle allowing the identification of the value of the mean field control problem as the unique viscosity solution in a class and sense as broad as possible, see 
\cite{ Burzoni,conforti,CossoGozziKharroubiPham,soner2022viscosity,WuZhang}. 
 Typically, these results cover the setting where $U$ is Lipschitz but may not be differentiable, which is expected to be the case when $\sF$ and $\sG$ are regular but not convex and thus optimizers may not be unique. In comparison, there are much fewer general results on the master equation of mean field games in the absence of uniqueness of the equilibria: the work of \cite{cecchin2022weak} gives a possible approach in the case of potential games, which are, by definition, derived from a control problem. At this stage, there are no general results on the convergence of the value functions of finite games to a possible value of the mean field game outside the analysis of  \cite{CardaliaguetDelarueLasryLions}. The best that is known are compactness results, see for instance \cite{djete2022extended,Fischer2017,Lacker2}. They suffice to establish convergence of the value functions at measures where the equilibrium is unique, but
questions of selection remain very challenging when uniqueness does not hold. 

\subsection{Our motivation.} 
The goal of the present paper is to understand the rate of convergence of the value functions $V^N$ to $U$ in the non-convex setting. Of course, as already discussed above, a quantitative answer to this convergence problem is already known when $\sF$ and $\sG$ are convex and sufficiently regular, with the (optimal) rate $1/N$. The convergence of $V^N$ to $U$ in the non-convex setting has received significant attention in the literature in recent years, and qualitative results
have been obtained in \cite{Lacker2017,DjetePossamaiTan} (see also 
\cite{CavagnariLisiniOrrieriSavare,FornasierLisiniOrrieriSavare,GangboMayorgaSwiech}
for deterministic dynamics or dynamics with a sole common noise). More recently, 
a first quantitative result outside the convex setting has been obtained, under `natural' assumptions, in a work \cite{cardaliaguet2023algebraic} by the first and last author with Cardaliaguet and Souganidis. The main result of \cite{cardaliaguet2023algebraic} (when specialized to the periodic setting) is the estimate 
\begin{align} \label{cdjsresult}
    |U(t,m_{\bx}^N) - V^N(t,\bx)| \leq C N^{- \gamma(d)},  
\end{align}
with $C$ independent of $N$ and $\gamma(d)$ depending only on the dimension $d$.
The value of the exponent $\gamma(d)$ is not given explicitly, but the calculations can be followed step by step, and it is clear that $\gamma$ decreases faster than (any multiple) of $1/d$. 

The result of \cite{cardaliaguet2023algebraic} shows that one can indeed have an algebraic convergence rate of convergence even when optimizers of the limiting problem are not unique. But it leaves open a very natural question, which we aim to investigate in this work:
\begin{align*}
    \text{What is the optimal rate of convergence of $V^N$ to $U$ in the non-convex setting?}
\end{align*} 
As far as we can tell, before the present paper there was not a clear conjecture about what the optimal rate should be, and in particular whether (or under what circumstances) it should be possible to obtain a dimension-free rate. Let us first emphasize that without convexity we cannot hope to obtain the rate $1/N$, as the calculations leading to \eqref{smoothcase} in the smooth case clearly indicate that the rate $1/N$ is tied to second-order regularity of $U$ in the measure variable, which we cannot expect without convexity of $\sF$ and $\sG$.
To gain some intuition, notice that because $U$ is only Lipschitz, the rate should be compared 
 to the one for the uncontrolled case with data which is only Lipschitz continuous, and this should in turn be related to the rate observed in the convergence of uncontrolled weakly
interacting particle systems. Due to the underlying statistical averaging phenomena, the latter
convergence rate is actually related to fundamental results in probability theory on the convergence
of empirical measures of an
$N$-sample. In this context, there are two rates of convergence which play an especially important role:
\begin{enumerate}[(a)]
    \item the rate $N^{-1/d}$ governs convergence of empirical measures for the Kantorovich-Rubinstein distance (also called the $1$-Wasserstein distance and denoted $d_1$) when\footnote{\label{foo:1} The exponent becomes $1/2$ when $d=1,2$, with an additional logarithmic correction in the rate of convergence when $d=2$, but we feel better to stick to the reference value $1/d$ throughout the introduction as it makes the presentation easier.}
 $d \geq 3$,
 see \cite{AtajiKomlosTusnady,DereichScheutzowSchottstedt,FournierGuillin}, 
 and also describes the typical minimal distance, in dimension $d$, between two particles within a cloud of uniformly drawn particles.
 \item The rate $N^{-1/2}$ corresponds to the central limit theorem, whose transposition to the convergence of empirical measures nevertheless requires some precautions, see \cite{fernandez-meleard,jourdain-meleard,meleard,sznitman1985fluctuation,tanaka1981central}.
\end{enumerate}
The difference between $(a)$ and $(b)$ lies in the class of test functions used to measure the convergence rate of the empirical measure: in $(a)$, the test functions are Lipschitz continuous, while in $(b)$, the test functions are much smoother. To illustrate point $(b)$, it is worth observing that the main fluctuation results in the literature on particle systems are stated in 
``sufficiently negative" Sobolev (Hilbertian) spaces. For instance, in \cite{meleard,fernandez-meleard}, 
fluctuations are 
estimated in the dual of a space of functions admitting $s:=1+ 2 \lfloor d/2 \rfloor$ generalized derivatives that are square integrable with respect to some heavy (polynomially) tailed measure $\nu$. 
Convergence of the fluctuations is obtained in a similar but larger space, obtained 
by replacing 
$s=1+ 2 \lfloor d/2 \rfloor$  by $s=4 + 2 \lfloor d/2 \rfloor$ 
and by changing accordingly the polynomial decay of the underling reference measure $\nu$. 
As we are working on the 
torus, the description of $\nu$ does not really matter here. Still, it is worth stressing that the mollification procedures 
implemented in the present paper 
also rely on the properties of the Hilbert space $H^{-s}$ for $s > d/2 + 1$, i.e., the dual of the Hilbert space of functions with $s$ generalized derivatives in $L^2$,
and could be used to recover the fact (hence already proven in 
\cite{meleard,fernandez-meleard})
that fluctuations are on average of order $N^{-1/2}$ when measured in the Hilbert space $H^{-s}$ for $s > d/2 + 1$. 

In other words, if $(\xi^i)_{i = 1,...,N}$ are i.i.d random variables with common law $m$, then
\begin{itemize}
 \item $
d_1(m_{\boldsymbol \xi}^N, m) := \sup_{g \, \text{1-Lip}} \int g \, d(m_{\boldsymbol \xi}^N - m)$
is typically of order $N^{-1/d}$, while
\item $
\|m_{\boldsymbol \xi}^N - m \|_{-s} := \sup_{ \|g\|_{s} = 1} \int g  \, d(m_{\boldsymbol \xi}^N - m)$
is typically of order $N^{-1/2}$ when $s > d/2 + 1$. 
\end{itemize}

The heuristic discussion above suggests the following conjecture: if we work under conditions on the data ($\sF$, $\sG$, and $H$) which guarantee only that $U$ is Lipschitz with respect to $d_1$, then the optimal rate should be $N^{-1/d}$, the size of typical fluctuations of empirical measures as measured with respect to $d_1$. If, on the other hand, the value function $U$ is Lipschitz with respect to a much weaker metric, like the one generated by $\| \cdot \|_{-s}$ for $s$ large enough, then it should be possible to obtain the rate $N^{-1/2}$, the size of typical fluctuations of emprical measures as measured with respect to this weaker metric.
Our objective is to verify this conjecture as far as possible. In order to do so, we limit our analysis to the periodic setting: this avoids any technicalities about the decay at infinity of the various functions that we manipulate. 
We also emphasize that this question has already been solved 
 for mean field control problems on a finite state space, see \cite{CecchinFinite,Kolokoltsov}: the convergence rate is shown to be  
 $1/\sqrt{N}$ when cost coefficients are non-convex in the measure argument, which is consistent  with case $(b)$ right above, keeping in mind that the regularity of the test functions does not matter in this case since the state space is finite. 

\subsection{Our results}
We work with two sets of conditions of the data $H$, $\sF$, and $\sG$. Assumption \ref{assump:d1} gives minimal conditions under which we can establish that $U$ is Lipschitz and semi-concave (defined below) with respect to $d_1$, and when Assumption \ref{assump:d1} is in force we say that we are in the ``$d_1$-regular case". Assumption \ref{assump:Hs} gives minimal conditions under which we can establish that $U$ is Lipschitz and semi-concave (defined below) with respect to $\|\cdot\|_{-s}$ for some $s > d/2 + 2$, and when Assumption \ref{assump:Hs} is in force we say that we are in the ``$H^{-s}$-regular case". 
The role of semi-concavity in both cases is outlined in Subsection 
\ref{subse:method}. We also refer to Remark 
\ref{rmk:semi-concave}.

Our contributions in the two ``$H^{-s}$ and $d_1$-regular cases''  can be summarized as follows.
\newline \newline 
 \textbf{Rates of convergence:} 
 In the $H^{-s}$-regular case, we obtain in Theorem \ref{thm:mainhs} the estimate 
 \begin{align} \label{hsrateintro}
    V^N(t,\bx) - C/N \leq U(t,m_{\bx}^N) \leq V^N(t,\bx) + C/\sqrt{N}, 
 \end{align}
 and in particular $|V^N(t,\bx) - U(t,m_{\bx}^N)| \leq C/\sqrt{N}$. This exactly matches the conjectured rate in the $H^{-s}$-regular case.
 
 In the $d_1$-regular case, we find (again taking $d \geq 3$ for simplicity and recalling footnote \ref{foo:1} for the peculiar cases $d=1,2$) that, for each $\eta > 0$, 
  there is a constant $C > 0$ such that
 \begin{align} \label{d1rateintro}
  V^N(t,\bx) - C N^{-1/d} \leq U(t,m_{\bx}^N) \leq V^N(t,\bx) + CN^{ -
\beta(d) + \eta}, \quad \beta(d) :=
  \frac{2}{3d + 6}.
 \end{align}
 In particular, we see that for $d$ large, the estimate on $|V^N(t,\bx) - U(t,m_{\bx}^N)|$ is roughly of order $N^{-2/(3d)}$. This is obviously slightly worse than the conjectured rate of $N^{-1/d}$, but still represents a significant improvement on existing results.
 
Finally, in the case where the coefficients are convex and Lipschitz with respect $d_1$, we establish the estimate (again taking $d \geq 3$ for simplicity) 
\begin{align} \label{convcase}
    0 \leq V^N(t,\bx) - U(t,m_{\bx}^N) \leq CN^{-1/d}.
\end{align}
Even though the convex case is much simpler to analyze, this result seems to be new. Indeed, since $\sF$ and $\sG$ are only assumed to be $d_1$-Lipschitz, $U$ is not expected to be smooth, and so the well-known argument based on ``projecting" $U$ fails. We instead use purely control-theoretic arguments to obtain \eqref{convcase}. We note in particular that the observation that $V^N \leq U$ in the convex regime appears to be new. 
In view of Example 1 presented in Subsection \ref{subsec:examples}, the rate in Proposition \ref{convcase} is sharp.
\newline \newline 
\textbf{Examples:} Of course, to provide a complete picture we need also to verify through examples that the conjectured rates $N^{-1/d}$ and $N^{-1/2}$ cannot be improved. In Subsection \ref{subsec:examples}, we first demonstrate in Example 1 that $N^{-1/d}$ is indeed the best possible rate when $\sF$ and $\sG$ are just $d_1$-Lipschitz. We note that it is easy to construct such an example if $H = 0$ (there is no control), but this does not rule out the possibility that strict convexity or coercivity of the Hamiltonian $H$ somehow benefits the convergence rate. Our example uses the ``model" Hamiltonian $\frac{1}{2} |p|^2$, and the idea is to use the Cole-Hopf transform to analyze the functions $V^N$ - this ultimately leads to an interesting probabilistic analysis related to the ``coupon-collector problem", which is presented in Section \ref{sec:convexandexproofs}.
Noticeably, this counter-example works due to the infinite dimensional nature of the optimal control problem that is treated here; to the best of our understanding, similar constructions, but in the Euclidean setting, would not provide interesting examples
in the study of vanishing viscosity for finite-dimensional Hamilton-Jacobi equations.

Indeed, for the $H^{-s}$-regular case, we show (in Example 2) that the convergence problem is in fact related to a question of vanishing viscosity for finite-dimensional Hamilton-Jacobi equations, with the viscosity being of order $1/N$.
And, remarkably, there are two distinct situations for this latter problem: $(1)$ when the costs are convex and smooth enough, the convergence rate of the value functions is linear in the viscosity, see 
\cite{FlemingSmall,FlemingSouganidis}; $(2)$ in general (but under reasonable regularity assumptions), the convergence rate is linear in the square root of the viscosity, 
see \cite{CrandallLionsvanishing,Evansvanishing,Lions82}, and this rate is 
claimed to be optimal. As for the latter point, it is however fair to say that the precise conditions under which the root of the viscosity is the optimal rate 
are rather unclear to us and, in particular, we do not know whether these conditions cover or not the type of Hamiltonians 
we use below. As we announced in the previous paragraph, the construction based 
on the same Cole-Hopf transformation as the one used below for proving that $N^{-1/d}$ is indeed the optimal rate in the ``$d_1$-regular case'' does not provide a relevant example in the vanishing viscosity problem for finite-dimensional Hamilton-Jacobi equations. Anyway,  
we believe that the rate obtained in \eqref{hsrateintro}  is the optimal rate in the $H^{-s}$-regular case.
\newline \newline 
\textbf{Regularity:} While the convergence problem is our main focus, we also provide some new regularity results for the value function $U$. In the $d_1$-regular case, we establish in Proposition \ref{prop:uregularityd1} that the value function $U$ is Lipschitz and semi-concave with respect to the metric $d_1$. This result is expected, except for the fact that we obtain the result without assuming that $\sF$ and $\sG$ are differentiable. In fact we only assume that $\sF$ and $\sG$ are Lipschitz and semi-concave with respect to $d_1$, so our assumptions on $\sF$ and $\sG$ in this result are in some sense optimal (certainly we cannot have a global Lipschitz and semi-concavity estimate for $U$ without assuming one for $\sG$). Without differentiability of $\sF$ and $\sG$ we do not have access to the usual description of optimizers in terms of a forward-backward PDE system, and we must instead proceed by a mollification procedure. We also prove in Proposition \ref{prop:upropshs} that in the $H^{-s}$-regular case, $U$ is Lipschitz and semi-concave with respect to $\| \cdot \|_{-s}$. This result appears to be new, and relies on stability estimates for a certain Fokker-Planck equation in negative Sobolev spaces which are presented in the Appendix.
\newline \newline
\textbf{Regularization procedures on the Wasserstein space:} In order to achieve our main convergence and regularity results, we implement several regularization techniques for functions on the space $\pr(\T^d)$. All three of these methods have appeared in some way in the literature before, but here we study for the first time their interplay with the Hamilton-Jacobi equation \eqref{hjbinf}. In particular, much of our analysis in Section \ref{sec:hardinequalities} is focused on analyzing the degree to which the regularization procedures preserve subsolutions of \eqref{hjbinf}. We believe the techniques we employ in this section could be useful, for example in studying the comparison principle for equations like \eqref{hjbinf}, where a procedure for approximating subsolutions by more regular subsolutions would clearly be useful.

\subsection{Our method} 
\label{subse:method}
As discussed above, one of the examples presented in Subsection \ref{subsec:examples} shows that the convergence problem we are considering can be viewed as an infinite-dimensional analogue of an evanescent viscosity problem in which the intensity of the evanescent noise would be $1/\sqrt{N}$. Accordingly our proof, both in the $d_1$-regular and $H^{-s}$-regular cases, takes up essential ideas from the analysis of finite-dimensional Hamilton-Jacobi equations with an evanescent viscosity, as written, for example, in the notes \cite{Calder}. 

The upper bounds in both \eqref{hsrateintro} and \eqref{d1rateintro} are more challenging, and we refer to them as the “hard inequalities". The basic idea for proving the hard inequality in the $H^{-s}$-regular case (i.e. the upper bound in \eqref{hsrateintro}) is as follows. Let us first suppose that $\Phi = \Phi(t,m) : [0,T] \times \pr(\T^d) \to \R$ is a function such that 
\begin{enumerate}[(i)]
    \item $\Phi(T,m) \leq \sG(m) + C_0$,
    \item $\Phi$ is a subsolution of \eqref{hjbinf}, up to an error of order $C_1$,
    \item $\|\tr(D^2_{mm} \Phi) \|_{\linf} \leq C_2$.
\end{enumerate}
Then, at least formally, the function $[0,T] \times (\T^d)^N \ni (t,\bx) \mapsto \Phi(t,m_{\bx}^N)$ is a subsolution of the equation describing $V^N$, up to an error of order $C_1 + \frac{C_2}{N}$, so that the comparison principle gives 
\begin{align} \label{phiestintro}
    \Phi(t,m_{\bx}^N) \leq V^N(t,\bx) + C_0 + T(C_1 + C_2/N).
\end{align}
This argument is made precise under appropriate regularity conditions on $\Phi$ in Proposition \ref{prop:projections}.

At first, the estimate \eqref{phiestintro} appears to be useless in the non-convex setting, since as discussed already we do not expect 
$D^2_{mm} U$ to be bounded. But suppose that we manage to produce regularization $(U^{\epsilon})_{\epsilon > 0}$ of $U$ such that for some $C > 0$, 
\begin{enumerate}
    \item $\| U - U^{\epsilon} \|_{\linf} \leq C \epsilon$,
    \item $U^{\epsilon}$ is a subsolution of \eqref{hjbinf} up to an error of order $C \epsilon$, and 
     \item $\|D^2_{mm} U^{\epsilon} \|_{\linf} \leq \frac{C}{\epsilon}$ 
\end{enumerate}
Then, the estimate \eqref{phiestintro}, together with the triangle inequality, gives 
\begin{align*}
    U(t, m_{\bx}^N) \leq  V^N(t,\bx) + C \Big( \epsilon + \frac{1}{N \epsilon} \Big),
\end{align*}
for some (new) constant $C$ independent of $N$.
Choosing $\epsilon = \frac{1}{\sqrt{N}}$ gives the upper bound in \eqref{hsrateintro}. 

Of course, the question is how to produce functions $U^{\epsilon}$ satisfying properties (1)-(3). In finite dimensions (and even in Hilbert spaces, see \cite{Lasry1986}), sup-convolution is known to be a convenient way to create $\cC^{1,1}$ regularity while preserving sub-solution properties 
(see for instance \cite[Chapter II]{Bardi}). This motivates our choice to define $U^{\epsilon}$ by 
\begin{align} \label{uepsdefintro}
    U^{\epsilon}(t,m) = \sup_{\nu} \Big\{U(t,\nu) - \frac{1}{2 \epsilon} \|m - \nu \|_{-s}^2 \Big\},
\end{align}
with the supremum taken over the space of probability measures on 
$\T^d$, and with $s > 0$ chosen appropriately.
Proposition \ref{prop:dmlipschitz} explains why $H^{-s}$ is a good Hilbert space to work with - roughly speaking, when $s > d/2 + 1$, bounds on $U^{\epsilon}(t,\cdot)$ in $\cC^{1,1}(H^{-s})$ (which arise naturally from the sup-convolution procedure) imply bounds on the $\linf$-norm of $D^2_{mm} U^{\epsilon}$ as a consequence of Sobolev embeddings. Together with arguments taken mostly from \cite{Lasry1986} and presented in Proposition \ref{prop:supconvproperties}, this makes it possible to verify that $U^{\epsilon}$ in fact satisfies the estimates appearing in (1) and (3) above, provided that the original value function $U$ is Lipschitz and semi-concave with respect to $\| \cdot \|_{-s}$. This explains why, in the $H^{-s}$-regular case, we must verify in Proposition \ref{prop:upropshs} that $U$ is Lipschitz and semi-concave with respect to $\| \cdot \|_{-s}$ when Assumption \ref{assump:Hs} is in force. Verifying that the functions $U^{\epsilon}$ satisfy point (2), the subsolution property, turns out to be much more subtle, and this is handled in Subsection \ref{subsec:uepsanalysis}. This is in fact the main challenge in the $H^{-s}$ regular case.

Of course, there are several technical issues to overcome when implementing this argument which we have ignored in the above outline, for example even after regularizing we do not in fact have access to $D^2_{mm} U^{\epsilon}$, so bounds on $\|
D^2_{mm} U^{\epsilon}\|_{\linf}$ have to be understood as bounds on the Lipschitz constant of $D_m U(t,m,x)$ in $m$ (with respect to $d_1$), and we must verify that bounds on the Lipschitz constant of $D_m U^{\epsilon}$ are in fact enough to execute the 
sketch of proof outlined above.

In the $d_1$-regular case, the general strategy for the upper bound in \eqref{d1rateintro} is the same - we want to approximate the value function $U$ by functions which are smoother, but are still close to being sub-solutions to \eqref{hjbinf}. But this time, we cannot directly apply sup-convolution in $H^{-s}$, because $U$ is only Lipschitz with respect to $d_1$. Instead, we start with a linear mollification procedure, which transforms $U$ (which is only regular with respect to $d_1$) into some function $U^{\delta}$ which is regular with respect to $\| \cdot \|_{-s}$ for any $s$. Then we choose an appropriate $s$ and apply sup-convolution in $H^{-s}$ to produce $U^{\delta,\epsilon}$, which has the required regularity. The idea, like in points (1)-(3) above, is then to estimate the distance $\|U - U^{\delta, \epsilon}\|_{\linf}$, the amount by which $U^{\delta, \epsilon}$ fails to be a sub-solution to \eqref{hjbinf}, and the size of $\|D^2_{mm}U^{\delta, \epsilon}\|_{\infty}$, all as functions of the parameters $\delta$ and $\epsilon$, and then choose $\delta$ and $\epsilon$ as appropriate functions of $N$ to conclude.

Unfortunately, there is a last difficulty, which is that, roughly speaking, we could only find an efficient estimate for the subsolution property of $U^{\delta,\epsilon}$ at measures $m$ which were bounded from below by a constant depending (explicitly) on $\delta$ and $\epsilon$, i.e. for measures $m$ such that $m \geq c(\delta,\epsilon) \leb$ for some appropriate $c$ (with $\leb$ denoting the Lebesgue measure on the torus). We refer to Lemma \ref{lem:udelteepsproperties} for a precise statement of this result. The constraint on $m$ arises when proving that 
$U^{\delta,\epsilon}$ inherits the regularity properties of $U^\delta$, which is in fact just possible where the supremum defining $U^{\delta, \epsilon}$ is achieved in the interior of $\pr(\T^d)$, with the interior being here defined with respect to the $L^\infty$ norm. 
This step is the heart of the analysis in the $d_1$-regular case, and requires a somewhat demanding technical result, see Proposition \ref{prop:threshold}. The fact that we can only obtain a good estimate at such measures necessitates a final transformation, in which we replace $U^{\delta, \epsilon}$ by $U^{\delta, \epsilon, \lambda}(t,m) = U^{\delta, \epsilon}(t, (1- \lambda) m + \lambda \leb)$, and it is the function $U^{\delta, \epsilon, \lambda}$ which, in the $d_1$-regular case, ultimately plays the role that $U^{\epsilon}$ played in the $H^{-s}$-regular case. The analysis of $U^{\delta}$, $U^{\delta, \epsilon}$, and $U^{\delta, \epsilon, \lambda}$ is carried out in subsection \ref{subsec:uepsdelta}.

For the lower bounds in \eqref{hsrateintro} and \eqref{d1rateintro}, called the `easy inequalities', we use a different, more control-theoretic argument. Very schematically, the idea is that optimal strategies identified in the asymptotic mean field regime are admissible in the finite setting, i.e., they can be played by a finite number of players. The converse is false  and explains why the bounds in \eqref{hsrateintro} and \eqref{d1rateintro} are not symmetric (see also \cite{cardaliaguet2023algebraic}, where the two inequalities are also treated separately).

\subsection{Further prospects} 
In the end, our regularization argument provides a different proof from \cite{cardaliaguet2023algebraic}: in the latter, more effort is spent on the particle system itself; here, we mostly work with the Hamilton-Jacobi equation. Compared to \cite{cardaliaguet2023algebraic}, we work in a more restrictive setting (in particular we work with periodic data and we do not address the important issue of common noise), but we obtain much sharper estimates. Despite this improvement, the obtained bounds are not all optimal. The most interesting possible improvement, in our view, would be to obtain the bound 
\begin{align*}
|U(t,m_{\bx}^N) - V(t,\bx) | \leq CN^{-1/d}
\end{align*} in the $d_1$-regular case (with the same corrections as in footnote \ref{foo:1} when 
$d=1,2$), which amounts to improving the lower bound in \eqref{d1rateintro} to match the upper bound. Such a result would confirm our conjecture about the optimal rate of convergence in the $d_1$-regular case. The reason why we obtain a slightly worse rate in \eqref{d1rateintro} is related to the difficulty of estimating the subsolution property of $U^{\delta, \epsilon}$. Roughly speaking, we could have a more efficient estimate of the subsolution property of $U^{\delta,\epsilon}$, which would circumvent the need for $U^{\delta, \epsilon ,\lambda}$ and lead to a rate closer to the optimal $N^{-1/d}$, if we could prove that the sup-convolution $U^{\delta, \epsilon}$ preserves the initial regularity of $U^{\delta}$ with respect to metrics other than $H^{-s}$ - in particular, it is not clear to us whether the sup-convolutions $U^{\delta, \epsilon}$ preserves the $d_1$-Lipschitz constant of $U^{\delta}$ on the whole space (and not only on the set of measures satisfying $m \geq c(\delta, \epsilon) \leb$, as we just explained). This would be immediate if $U^{\delta}$ were defined over the entire Sobolev space $H^{-s}$, but things are much more difficult here, since the original problem is set on the  smaller space of probability measures, which is of empty interior for any negative Sobolev norm. As  discussed above, we partially circumvent this issue through Proposition \ref{prop:threshold}, which requires a somewhat involved technical analysis, and it is possible that the latter could be refined. In any case, it seems clear to us that some new ideas are required in order to obtain the optimal rate $N^{-1/d}$ using our techniques.

A related possibility to obtain the optimal $N^{-1/d}$ rate in the $d_1$-regular case would be to find another way to regularize the value function with similar features: the convergence rate should be explicit, the regularized value function should be regular enough (with explicit bounds) and should be a subsolution of the Hamilton-Jacobi equation up to an explicit residual term. In this regard, the sup-convolution is very convenient, but passing through a Hilbertian structure may seem somewhat unintuitive. At this stage, we have no intuition about the possible existence of a `better' regularization technique. What is certain is that the same Sobolev spaces have been used in different contexts related to ours: not only in the analysis of fluctuations for non-controlled particle systems, as already mentioned above, but also in the study of a comparison principle for viscosity solutions of the Hamilton-Jacobi equation in the very recent work \cite{soner2022viscosity}. We are also convinced that the regularization technique we use could allow us to reprove a similar comparison principle. This leads us to believe that obtaining regularization methods (by sup-convolution or otherwise), specifically adapted to the Hamilton-Jacobi equation \eqref{hjbinf}, is of broader interest than that of this work.

Another related question is to understand what our results say about the convergence of the optimal trajectories, i.e., the convergence of the optimal trajectories for the $N$-particle problem towards the optimal trajectories for the limiting. This problem is quite subtle, since optimal trajectories may not be unique when $\sF$ and $\sG$ are not convex. However, a recent result by Cardaliaguet and Souganidis \cite{cardaliaguet-souganidis:2} identifies an open and dense set $\mathcal O \subset [0,T] \times \pr(\T^d)$ where $U$ is locally of class $\cC^{1,1}$, and shows that optimal trajectories which start in $\mathcal O$ remain there, and in particular, $U$ is $\cC^{1,1}$ in a tube around an optimal trajectory which starts in $\mathcal O$. As shown in \cite{cardaliaguet-souganidis:2}, this leads to a propagation of chaos result for initial conditions in $\mathcal O$, provided that one already has a rate of convergence of $V^N$ to $U$. Lemma 3.3 of \cite{cardaliaguet-souganidis:2} explains clearly how the convergence rate of $V^N$ to $U$ impacts the rate of propagation of chaos, since $R_N$ in that statement is just $\sup_{(t,\bx)} |V^N(t,\bx) - U(t,m_{\bx}^N)|$. In particular, our estimates can be immediately ``plugged into" Lemma 3.1 of \cite{cardaliaguet-souganidis:2} to improve the main results of that paper (at least in the case of periodic data). The only unsatisfying point is that in the $H^{-s}$-regular case, the resulting rate of propagation of chaos will be dimension-free (in fact, the rate is $N^{-1/4}$), but only when measured up to stopping times $\tilde{\tau}^N$ which satisfies 
\begin{align} \label{stoppingtimesest}
    \bP[\tilde{\tau}^N < T] \leq CN^{-1/(d+8)}.
\end{align}
So, dimension still enters the propagation of chaos result through the asymptotic behavior of the stopping times. In order to get a truly dimension-independent propagation of chaos statement, one would need to replace the stopping times appearing in Lemma 3.3 of \cite{cardaliaguet-souganidis:2} with stopping times satisfying a dimension-free analogue of \eqref{stoppingtimesest}. To do this, we believe it would suffice to estimate the radius of the aforementioned tube as measured with respect to $\|\cdot\|_{-s}$ rather than $d_1$, but we do not pursue this analysis here. 

Of course, it would also be interesting to extend the analysis to the Euclidean (non-periodic) case.
We believe that the weighted Sobolev spaces used in \cite{fernandez-meleard,meleard} in the analysis of the CLT for uncontrolled particle systems may also be useful here. The case with a common noise (which, as mentioned above, was treated in \cite{cardaliaguet2023algebraic}) would also deserve some attention. Finally, it is certainly worth noting that, to our knowledge, the case of (non-cooperative) games remains beyond the reach of the methods developed here.
\raggedbottom

\subsection{Organization of the paper}
The article is organized as follows. Section 2 begins with a discussion of relevent notation and function spaces in Subsection \ref{subsec:notation}, and then Subsections\ref{subsec:formulation} and \ref{subsec:mainresults} contain the problem statement and main results. We also discuss the aforementioned examples in Subsection \ref{subsec:examples}, and discuss the connection between analysis in $\pr(\T^d)$ and analysis in $H^{-s}$ in Subsection \ref{subsec:connectionprhs}. Section 3 is devoted to establishing several fundamental properties of the value function in both the $d_1$-regular and $H^{-s}$-regular cases. In Section \ref{sec:regularization}, we propose and study several regularization methods for functions defined on the space of probability measures, including the sup-convolution method as well as the convolution method introduced in \cite{cecchin2022weak}. Section \ref{sec:hardinequalities} is the true heart of the paper, with the establishment of the “hard inequalities", i.e. the upper bounds in \eqref{hsrateintro} and \eqref{d1rateintro}. Section \ref{sec:easyinequalities} contains the proofs of the corresponding “easy inequalities", while Section \ref{sec:convexandexproofs} contains some arguments related to the examples and the estimate \eqref{convcase} in the convex case. Finally, the Appendix contains a number of auxiliary results for finite dimensional linear PDEs and Hamilton-Jacobi-Bellman equations as well as mollification argument for functions defined on the space of probability measures which is borrowed from \cite{cecchin2022weak}.

\section{Preliminaries and main results}

\label{sec:mainresults}

\subsection{Function spaces and notation} \label{subsec:notation}

Firstly, as mentioned above, $\T^d = \R^d/ \Z^d$ is the $d$-dimensional flat torus, whose general element is written $x = (x_1,...,x_d)$. We write $\bx = (x^1,...,x^N)$ for the general element of $(\T^d)^N$, with each $x^i = (x^i_1,...,x^i_d) \in \T^d$. We denote by $\pr(\T^d)$ the set of probability measures on $\T^d$. We endow this space with the Wasserstein metric $d_1$, defined by 
\begin{align*}
d_1(m_1,m_2) = \sup_{g \,\, 1-\text{Lip}} \int g(x) \bigl(m_1 - m_2\bigr)(dx), 
\end{align*}
the supremum being taken over $1$-Lipschitz functions $g : \T^d \to \R$. Of course, by Kantorovich duality we have an equivalent definition in terms of couplings which we also use when convenient.

We are next going to describe several spaces of functions. First, we mention that we will use the notation {\rm Lip}($A$) for the set of Lipschitz functions on a space $A$ (the metric on $A$ always being understood from context). Next, we fix some notation for multi-indices. We define multi-indices to be tuples $\bm j = (j_1,...,j_d)$ with each $j_i \in \{0,1,2,...\}$, and for $f : \T^d \to \R$ we interpret $D^{\bm j} f$ as 
\begin{align*}
    D^{\bm j} f = D_{x_1}^{j_1} f ... D_{x_d}^{j_d} f. 
\end{align*}
We write $|\bm j| = \sum_{i = 1}^d j_i$ for the order of the multi-index $\bm j$. We include the possibility that $\bm j = (0,...,0)$, in which case $D^{\bm j} f = f$. 

For $k \in \N$ (with $\N := \{1,2,...\}$), we define $\cC^k = \cC^k(\T^d)$ to be the space of functions $\T^d \to \R$ with continuous derivatives up to order $k$. We endow $\cC^k$ with the norm 
\begin{align*}
    \| f \|_{\cC^k} = \sum_{0 \leq | \bm j | \leq k} \|D^{ \bm j} f\|_{\linf}, 
\end{align*}
where $\| \cdot \|_{\linf}$ denotes the usual supremum norm. 

For $f \in L^2(\T^d)$, we define the Fourier coefficients $(\hat{f}^{\bk})_{\bk \in \Z^d}$ by 
\begin{align*}
    \hat{f}^{\bk} = \int_{\T^d} e^{\i 2 \pi \bk \cdot x} f(x) dx, \quad \i^2=-1.
\end{align*}
Here we are writing $\bk = (k_1,...,k_d)$ for an element of $\Z^d$.
For $s \geq 0$, we define the Sobolev space $H^s$ as the set of $f \in L^2(\T^d)$ such that 
\begin{align*}
    \norm{f}^2_s \coloneqq \sum_{\bk \in \Z^d} |\hat{f}^{\bk}|^2 \big(1 + \sum_{i = 1}^d |k^i|^{2s}) < \infty. 
\end{align*}
The space $H^s$ is a Hilbert space, with the norm $\| \cdot \|_s$ arising from the inner product 
\begin{align*}
\langle f, g \rangle_s = \sum_{\bk \in \Z^d} \hat{f}^{\bk} \overline{\hat{g}^{\bk}} \big(1 + \sum_{i = 1}^d |k_i|^{2s}),
\end{align*}
with the bar above denoting the usual conjugate of a complex number. 
We emphasize that this is not the usual definition of $\| \cdot \|_s$ -- it would be more typical to replace $(1 + \sum_{i = 1}^d |k_i|^{2s})$ with $(1 
 +|\bk|^2)^s$. However, the norms generated by these two choices are equivalent, and in particular, this means that we have access to the usual Sobolev embeddings, e.g. for $s > d/2$ we have $ \| f \|_{\linf} \leq C \|f\|_{s}$
for a constant $C$ independent of $f$. The reason for using this particular choice of norm will become clear in the proof of Lemma \ref{lem:faadibruno}.

Next, for $s \geq 0$ we define the space $H^{-s}$ to be the dual of $H^s$, i.e. the set of bounded linear functionals $q : H^s \to \R$. We define the Fourier coefficients $(\hat{q}^{\bm k})_{\bm k \in \Z^d}$ of $q \in H^{-s}$ by $\hat{q}^{\bm k} = q ( x \mapsto e^{\i 2 \pi {\bm k} \cdot x}).$ It is easy to check that the norm $\| \cdot \|_{-s}$ on $H^{-s}$ inherited by duality is induced by the inner product
\begin{align} \label{hminus}
    \langle p, q \rangle_{-s} \coloneqq \sum_{\bm k \in \Z^d} \hat{p}^{\bk} \overline{\hat{q}^{\bk}}\Big(1 + \sum_{i = 1}^d |k^i|^{2s} \Big)^{-1}.
\end{align}
For $f \in H^s$, we denote $f^*$ for the dual element of $H^{-s}$, i.e. $\langle f^* , q \rangle_{-s} = \langle f, q \rangle_{s,-s} \coloneqq q(f)$ for all $q \in H^{-s}$.
Likewise, given $q \in H^{-s}$, we write $q^*$ for the element of $H^s$ such that 
    $\langle q^*, f \rangle_s = q(f)$ for all $f \in H^s$.
Using \eqref{hminus}, it is easy to check that
\begin{align} \label{dualformula}
    \widehat{f^*}^{ \bk} = \Big(1 + \sum_{i = 1}^d |k^i|^{2s} \Big) \hat{f}^{\bk}, \quad f \in H^s, \text{ and }
     \widehat{q^*}^{ \bk} = \Big(1 + \sum_{i = 1}^d |k^i|^{2s} \Big)^{-1} \hat{q}^{\bk}, \quad q \in H^{-s}.
\end{align}
We follow Chapter 5 of \cite{CarmonaDelarue_book_I} for the definition of the linear derivative of a function $\pr(\T^d) \to \R$. In particular, given a continuous function $\phi = \phi(m) : \pr(\T^d) \to \R$, we say that $\phi \in \cC^1(\pr(\T^d))$ if there is a continuous map $\frac{\delta \phi}{\delta m} = \frac{\delta \phi}{\delta m} (m,y) : \pr(\T^d) \times \T^d \to \R$ such that 
\begin{align} \label{def:linderiv}
    \phi(m_2) = \phi(m_1) + \int_0^1 \int_{\T^d}\frac{\delta \phi}{\delta m} (t m_2 + (1-t) m_1, y) (m_2 - m_1)(dy)dt.
\end{align}
The equation \eqref{def:linderiv} determines $\frac{\delta \phi}{\delta m}$ only up to an $m$-dependent constant. We will make the additional normalization convention
\begin{align} \label{normalization}
    \int_{\T^d} \frac{\delta \phi}{\delta m}(m,x)dx=0, 
\end{align}
under which $\frac{\delta \phi}{\delta m}$ is indeed unique. Thus if $\phi \in \cC^1(\pr(\T^d))$ then we can refer to the (unique) continuous function $\frac{\delta \phi}{\delta m}$ satisfying both \eqref{def:linderiv} and \eqref{normalization} as the linear derivative of $\phi$. Notice that our normalization convention differs to the usual one (which requires instead that the derivative at $m$ has zero mean with respect to $m$). The reason is that, under our convention, the zeroeth Fourier coefficient of $\frac{\delta \phi}{\delta m}(m, \cdot)$ is zero, which simplifies an argument in the proof of Proposition \ref{prop:threshold}.
If $\phi$ has a linear derivative $\frac{\delta \phi}{\delta m}$ which is $C^1$ in $y$ for each $m$, then we can define the L-derivative $D_m \phi : \pr(\T^d) \times \T^d \to \R^d$ by
\begin{equation}
    D_m \phi(m,y) = D_y \frac{\delta \phi}{\delta m}(m,y). 
\label{eq:newlabel:Dm}
\end{equation}
Unlike $\frac{\delta \phi}{\delta m}$, $D_m \phi$ is uniquely defined if it exists, see again \cite[Chapter 5]{CarmonaDelarue_book_I}.

\subsection{Problem formulation} \label{subsec:formulation}

We fix a time horizon $T > 0$ and a filtered probability space $
\big(\Omega, \sF, \bbF = (\sF_t)_{0 \leq t \leq T}, \bP \big) 
$
satisfying the usual conditions, hosting independent $d$-dimensional Brownian motions $W$ and $(W^i)_{i \in \N}$, and such that $\sF_0$ is atomless. Our data consists of three functions
\begin{align*}
\sF, \sG : \pr(\T^d) \to \R, \quad L = L(x,a) : \T^d \times \R^d \to \R, 
\end{align*}
and the Lagrangian $L$ determines a Hamiltonian $H$ via the usual formula
\begin{align} \label{hamdef}
H(x,p) = \sup_{a \in \R^d} \big\{-L(x,a) - a \cdot p \big\}.
\end{align}
For $N \in \N$, $V^N$ is defined by
\begin{align} \label{vndef}
V^N(t_0,\bx_0) = \inf_{\bm \alpha \in \sA^N} \E\bigg[\int_{t_0}^T \Big(\frac{1}{N} \sum_{i = 1}^N L(X_t^i, \alpha_t^i) + \sF(m_{\bX_t}^N) \Big)dt + \sG(m_{\bX_T}^N) \bigg], 
\end{align}
for $(t_0,{\bm x}_0) \in [0,T] \times {\mathbb T}^d$,
subject to 
\begin{align} \label{xdef}
    dX_t^i = \alpha_t^i dt + \sqrt{2} dW_t^i, \quad t_0 \leq t \leq T, \quad X_{t_0}^i = x_0^i,
\end{align}
where the infimum is taken over the set $\sA^N$ of all square-integrable and progressively measure $(\R^d)^N$-valued processes $\bm \alpha = (\alpha^1,...,\alpha^N)$. We note that the Brownian motions $(W^i)_{i=1,\cdots,N}$ and the processes $(X^i)_{i=1,\cdots,N}$ in \eqref{xdef} are understood as taking values in $\R^d$, but they determine $\T^d$-valued processes in a canonical way (by composition with quotient map $\R^d \to \R^d / \Z^d$), which allows us to regard $m^N_{\bX_t}$, for each $t \in [0,T]$, as a probability measure on $\T^d$.
The function $U$, meanwhile, is given by 
\begin{align} \label{def:U}
U(t_0,m_0) = \inf_{(m,\alpha)} \bigg\{\int_{t_0}^T \Big(\int_{\R^d} L \big(x,\alpha(t,x) \big) m_t(dx) + \sF(m_t) \Big)dt + \sG(m_T) \bigg\}, 
\end{align}
where the infimum is taken over all pairs $(m,\alpha)$ with $m = (m_t)_{t_0 \leq t \leq T} \in C([t_0,T], \pr(\T^d))$ ($\mathcal{P}(\T^d)$ being equipped with the $d_1$ distance) and 
$\alpha  : (t,x) \in [t_0,T] \times \T^d \mapsto \alpha(t,x) \in \R^d$ satisfying (in a weak distributional sense) the Fokker-Planck equation
\begin{align*}
\partial_t m_t = \Delta_{x} m_t - \div_{x}\bigl(m_t \alpha(t,\cdot) \bigr), \quad \text{ in } [t_0,T) \times \T^d, \quad
m_{t_0} = m_0,
\end{align*}
and the integrability condition
\begin{align*}
\int_{t_0}^T \int_{\T^d}|\alpha(t,x)|^2 m_t(dx)dt < \infty.
\end{align*}
We note that because there is no common noise, there is no subtlety in checking (under mild conditions on the data, and in particular under Assumption \ref{assump:d1} below) that 
\begin{align} \label{udefalt}
U(t_0,m_0) = \inf_{\alpha \in \sA} \E\bigg[\int_{t_0}^T \Big(L(X_t, \alpha_t) + \sF(\sL(X_t)) \Big) dt + \sG(\sL(X_T)) \bigg], 
\end{align}
subject to 
\begin{align} \label{xdefinf}
dX_t = \alpha_t dt + \sqrt{2} dW_t, \quad t_0 \leq t \leq T, \quad X_{t_0} = \xi \sim m_0,
\end{align}
where the infimum is taken over the set $\sA$ of square-integrable progressive processes $\R^d$-valued $\alpha =(\alpha_t)_{t_0 \le t \le T}$ and $\sL(X_t) \in \pr(\T^d)$ denotes, for any $t \in [t_0,T]$, the law of $X_t$ (when viewed as a $\T^d$-valued random variable). Indeed, this is a consequence of a so-called mimicking argument, see e.g. the discussion in \cite[Section 8]{Lacker_superposition}. We note that this alternative formulation is not crucial for any of our arguments, but at times is convenient, e.g. in the presentation of the examples in Subsection \ref{subsec:examples}.

The function $V^N$ is, under mild assumptions on the data (which are satisfied below under each of the two sets of assumptions spelled out in Subsection \ref{subsec:mainresults} below), the unique smooth solution of \eqref{hjbn}.
The function $U$ is expected to satisfy, in an appropriate viscosity sense, the equation \eqref{hjbinf}.

\subsection{Assumptions and main results} \label{subsec:mainresults}

We have two main sets of assumptions, consistent with the $d_1$-regular and the $H^{-s}$-regular cases already discussed at length in the introduction.

\begin{assumption}[Assumptions for the $d_1$-regular case] \label{assump:d1}
We assume that there exists a constant $C \geq 1$ such that
\begin{enumerate}
    \item $H \in C^{2}(\T^d \times \R^d)$, and 
    for each $(x,p) \in \T^d \times \R^d$,
    \begin{align*}
        \frac{1}{C} I_{d \times d} \leq D^2_{pp} H(x,p) \leq C I_{d \times d}.
    \end{align*}
    \item For all $(x,p) \in \T^d \times \R^d$, 
       \begin{equation}
       |D_x H(x,p)| \leq C(1 + |p|).
       \label{AssumptiononD_xH}
       \end{equation}
    \item $\sF$ and $\sG$ are Lipschitz and semi-concave with respect to $d_1$ and to the constant $C$. Semi-concavity with respect to $d_1$ and to $C$ means (explained only for $\sF$) that
\begin{equation} \label{d1semiconcave}
\mathcal{F}((1-\lambda) m_0 + \lambda m_1) \geq (1-\lambda) \mathcal{F}(m_0) + \lambda \mathcal{F}(m_1) - \lambda (1-\lambda) \frac{C}{2}d_1^2(m_0,m_1), 
\end{equation}
for all $m_0, m_1 \in \pr( \T^d)$ and all $\lambda \in (0,1)$.
\end{enumerate}
\end{assumption}

\begin{rmk} \label{rmk:lprops}
Let us remark that the conditions on $H$ in Assumption \ref{assump:d1} are fairly minimal. The strict convexity of $H$ is standard, and ensures, among other things, that $L$ inherits some regularity from $H$. Moreover, the bound $D_{pp}^2 H \leq CI_{d\times d}$ implies that $L$ is coercive, in the sense that there is a constant $C$ such that 
$L(x,a) \geq - C + \frac{1}{C} |a|^2$. The condition \eqref{AssumptiononD_xH} is standard to find $L^{\infty}$ estimates for the optimal controls in the problems defining $U$ and $V^N$. The assumptions on $H$ are satisified for instance if $H(x,p) = \frac{1}{2}|p|^2 + H_1(x) \cdot p$ for some smooth vector field $H_1 : \T^d \rightarrow \R^d$.
\end{rmk}
\begin{rmk}
\label{rmk:semi-concave}
{The $d_1$-Lipschitz continuity of $\sF$ and $\sG$ required in Assumption \ref{assump:d1} is quite natural in this context, but the semi-concavity condition \eqref{d1semiconcave} may be more surprising. But a simple computation shows that \eqref{d1semiconcave} holds (for some constant $C$) if $\sF$ admits two linear derivatives $ \frac{\delta \sF}{\delta m}(m,y)$ and $\frac{\delta^2 \sF}{\delta m^2}(m,y,z)$, the latter being defined analogously to the former, and 
    \begin{align*}
    \forall m \in \pr(\T^d), \ \forall y,z \in \T^d, \quad 
     D_y D_z \frac{\delta^2 \sF}{\delta m^2}(m,y,z) \leq C \textrm{\rm I}_d, 
    \end{align*}
    with $\textrm{\rm I}_d$ denoting the $d$-dimensional identity matrix and 
    the symbol $\leq$ standing for the usual comparison between symmetric matrices. 
    The reader will notice that the left-hand side of the inequality is nothing but the second-order Lions derivative $D^2_{mm} \sF(m,y,z)$}. Thus the conditions on $\sF$ and $\sG$ in Assumption \ref{assump:d1} are met, for instance, under the standing assumptions of \cite{cardaliaguet2023algebraic}. A typical example of a function $\sG$ satisfying the conditions of Assumption \ref{assump:d1} would be
    \begin{align*}
        \sG(m) = G\bigl(m(\phi_1),...,m(\phi_k)\bigr),
    \end{align*}
    where $k \in \N$, $G : \R^k \to \R$ is $\cC^2$ with bounded first and second derivatives, and $\phi_1,...,\phi_k : \T^d \to \R$ are $\cC^1$ with bounded derivatives. 
\end{rmk}

\begin{assumption}[Assumptions for the $H^{-s}$-regular case] \label{assump:Hs}
We assume that there is a number $s \in \N$ with $s > d/2 +2$, such that 
\begin{enumerate}
    \item $H$ belongs to $C^{s}(\T^d \times \R^d)$ and satisfies conditions (1) and (2) from Assumption \eqref{assump:d1}.
    \item $\sF$ and $\sG$ belong $\mathcal{C}^1(\mathcal{P}(\T^d))$ with 
\begin{equation}
\sup_{m \in \mathcal{P}} \norm{ \frac{\delta \mathcal{F}}{\delta m}(m,\cdot)}_{C^s} +\sup_{m \in \mathcal{P}} \norm{ \frac{\delta \mathcal{G}}{\delta m}(m,\cdot) }_{C^s} < +\infty{.}
\label{Assumptiononlinearderivative24Mars}
\end{equation}
    Moreover $\mathcal{F}$ and $\mathcal{G}$ are semi-concave with respect to the norm $\| \cdot \|_{-s}$.
\end{enumerate}
\end{assumption}

\begin{rmk} \label{rmk:hslip}
In Assumption \ref{assump:Hs}, we keep the same basic structural conditions on the Hamiltonian, but require in addition ``order $s$ smoothness" of all the data.
It is easy to check, using the definition of the linear derivatives that the condition \eqref{Assumptiononlinearderivative24Mars} implies that $\sF$ and $\sG$ are Lipschitz with respect to $\| \cdot \|_{-s}$. A typical example of a function $\sG$ satisfying condition (2) in Assumption \ref{assump:Hs} would be 
\begin{align*}
        \sG(m) = G(m(\phi_1),...,m(\phi_k))
    \end{align*}
    where $k \in \N$, $G : \R^k \to \R$ is $\cC^2$ with bounded first and second derivatives, and $\phi_1,...,\phi_k \in \mathcal{C}^s(\T^d)$.
\end{rmk}

Our main result in the $d_1$-Lipschitz case is the following. 

\begin{thm} \label{thm:maind1}
    Let Assumption \ref{assump:d1} hold. Then for each $\eta > 0$, there is a constant $C$ such that for each $N \in \N$, we have, 
    for all $(t,\bx) \in [0,T] \times (\T^d)^N$, 
    \begin{align*}
        V^N(t,\bx) - C R_{d,N} \leq U(t,m_{\bx}^N) \leq V^N(t,\bx) + CN^{ -\beta(d)+\eta}, \quad \textrm{\rm with} \ \beta(d):= \frac{2}{3d+ 6},
    \end{align*}
    where 
    \begin{align} \label{fournierguillinrate}
        R_{d,N} \coloneqq \begin{cases}
            N^{-1/2} & d = 1, \\
            N^{-1/2} \log(1 + N) & d = 2, \\
            N^{-1/d} & d \geq 3.
        \end{cases}
    \end{align}
\end{thm}

And now we present our main result in $H^{-s}$-Lipschitz case.

\begin{thm} \label{thm:mainhs}
    Let Assumption \ref{assump:Hs} hold. Then there is a constant $C$ such that for each $N \in \N$, 
    \begin{align*}
        V^N(t,\bx) - \frac{C}{N} \leq U(t,m_{\bx}^N) \leq V^N(t,\bx) + \frac{C}{\sqrt{N}},
    \end{align*}
    for all $(t,\bx) \in [0,T] \times (\T^d)^N$. 
\end{thm}

Finally, we state a (much simpler) result in the convex setting.

\begin{prop} \label{prop:convexlip}
Suppose that Assumption \ref{assump:d1} is in force, and that in addition $\sF$ and $\sG$ are convex. Then there is a constant $C$ such that
\begin{align*}
    0 \leq V^N(t,\bx) - U(t,m_{\bx}^N) \leq CR_{N,d},
\end{align*}
for all $N \in \N$ and all $(t,\bx) \in [0,T] \times (\T^d)^N$, where $R_{N,d}$ is as defined in \eqref{fournierguillinrate}. 
\end{prop}

\subsection{Examples} \label{subsec:examples}

In this section we give two examples which pinpoint the best convergence we can expect. In both cases, it is more convenient to state the examples in the Euclidean space $\R^d$ rather than $\T^d$.
\newline \newline 
\textbf{Example 1:}  In this example we define $V^N : [0,T] \times (\R^d)^N \to \R$ by
\begin{align} \label{vndefex2}
    V^N(t_0,\bx_0) = \inf_{\bm \alpha \in \sA^N} \E\bigg[\int_{t_0}^T \frac{1}{2N} \sum_{i = 1}^N |\alpha^i_s|^2ds + d_1(m_{\bX_T}^N,\cN_{T}) \bigg]{,}
\end{align}
subject to dynamics
\begin{align} \label{xdynamicsex1}
    dX_t^i = \alpha_t^i dt + dW_t^i, \quad X_{t_0}^i = x_0^i, 
\end{align}
with $\cN_{T}$ denoting the standard Gaussian measure on $\R^d$ with mean $0$ and variance $T I_{d \times d}$, and the infimum taken over all square-integrable and progressive processes $\bm \alpha = (\alpha^1,...,\alpha^N)$. We note that the omission of the $\sqrt{2}$ the dynamics \eqref{xdynamicsex1} was intentional, and based purely on notational convenience. Similarly, we define $U  : (t,m) \in [0,T] \times \sP_2(\R^n) \mapsto U(t,m) \in \R$ by
\begin{align} \label{udefex2}
    U(t_0,m_0) = \inf_{\alpha \in \sA} \E\bigg[\int_{t_0}^T \frac{1}{2} |\alpha_t|^2 dt + d_1(\sL(X_T), \cN_{T}) \bigg], 
\end{align}
subject to the dynamics \begin{align*}
    dX_t = \alpha_t dt + dW_t, \quad X_{t_0} = \xi \sim m_0.
\end{align*}
Then we have the following: 
\begin{prop} \label{prop:ex2}
    Let $V^N$ and $U$ be defined by \eqref{vndefex2} and \eqref{udefex2}. Then we have $U(0,\delta_0) = 0$, and if $T$ is large enough then there is a constant $c > 0$ such that for infinitely many values of $N$,
    \begin{align*}
        V^N(t,0) \geq c N^{-1/d}.
    \end{align*}
\end{prop}
Of course, since $m \mapsto d_1(m, \cN_{T})$ is clearly $d_1$-Lipschitz, Proposition \ref{prop:ex2} shows that we cannot have a rate better than $N^{-1/d}$ when the data is only $d_1$-Lipschitz. The proof of Proposition \ref{prop:ex2} is given in Section \ref{sec:convexandexproofs}. We note that the assumption that $T$ is large enough is just to make the proof more straightforward, and it is also clear that the same argument should work for \textit{all} $N$ large enough, rather than for infinitely many $N$ -- we state things this way just to make the proof as transparent as possible.
\newline \newline 
\textbf{Example 2:}
Given two costs functionals $F,G : \R^d \rightarrow \R$, consider the functions $V^N : [0,T] \times (\R^d)^N \to \R$ given by 
\begin{align} \label{vndefex1}
    V^N(t_0,\bx_0) = \inf_{\bm \alpha  \in \sA^N} \E\bigg[\int_{t_0}^T \Big(\frac{1}{N} \sum_{i = 1}^N L(\alpha^i(t,X_t^i)) + F\big(\frac{1}{N} \sum_{i = 1}^N X_t^i\big)  \Big)dt + G\big(\frac{1}{N} \sum_{i = 1}^N X_T^i\big) \bigg],
\end{align}
subject to the dynamics \eqref{xdef}
as well as the function $U : [0,T] \times \sP_2(\R^n) \mapsto U(t,m) \in \R$ given by 
\begin{align} \label{udefex1}
    U(t,m) = \inf_{\alpha} \E\bigg[\int_{t}^T \Big(L(\alpha(s,X_s)) + F(\E[X_s])  \Big)ds + G(\E[X_T]) \bigg],
\end{align}
subject to the dynamics \eqref{xdefinf}. In other words, we have defined $V^N$ and $U$ as in Subsection \ref{subsec:formulation} but on the whole space and with 
\begin{align*}
    \sF(m) = F(\bar{m}), \quad \sG(m) = G(\bar{m}), \quad \bar{m} \coloneqq \int_{\mathbb R} x dm(x).
\end{align*}
Then we have the following result, whose elementary proof is presented in Section \ref{sec:convexandexproofs}.
\begin{prop} \label{prop:ex1}
Suppose that $F$ and $G$ are Lipschitz, and $L = L(a)$ is $C^2$ with $\frac{1}{C}I_{d \times d} \leq  D^2 L \leq C I_{d \times d}$ for some constant $C \geq 1$. Let $V^N$ and $U$ be defined as in \eqref{vndefex1} and \eqref{udefex1}.
Then we have 
\begin{align} \label{ex1claim}
    U(t,m) = u(t,\bar{m}), \quad V^N(t,\bx) = v^N(t,\bar{m_{\bx}^N}) = v^N\Bigl(t,\frac{1}{N} \sum_{i = 1}^N x^i\Bigr), 
\end{align}
where $u$ and $v^N$ are the unique viscosity solutions of the finite-dimensional PDEs 
\begin{align*}
    -\partial_t u(t,x) + H\bigl(D_x u(t,x)\bigr) = F(x), \text{ for $(t,x) \in [0,T) \times \R^d$,} \quad u(T,x) = G(x), \text{ for $x \in \R^d$}.
\end{align*}
and 
\begin{equation*}
\begin{split}
&    -\partial_t v^N(t,x) - \frac{1}{N} \Delta_x v^N(t,x)+ H\bigl(D_{x}v^N(t,x)\bigr) = F(x), \text{ for $(t,x) \in [0,T) \times \R^d$,} 
\\
& v^N(T,x) = G(x), \text{ for $x \in \R^d$},
\end{split}
\end{equation*}
 with $H$ denoting the Hamiltonian defined by \eqref{hamdef}.
 \end{prop}
 
 From Proposition \ref{prop:ex1} it is clear that 
\begin{align*}
\sup_{(t,\bx) \in [0,T] \times (\R^d)^N} |V^N(t,\bx) - U(t,m_{\bx}^N)| = \sup_{(t,x) \in [0,T] \times \R^d} |v^N(t,x) - u(t,x)|, 
\end{align*}
i.e. the convergence rate of $V^N$ to $U$ is completely governed by that of $v^N$ to $u$. As discussed in the introduction, this is a well-studied problem, and classical results show that the
rate is
\begin{align*}
\sup_{(t,x) \in [0,T] \times \R^d} |v^N(t,x) - u(t,x)|  \leq \frac{C}{\sqrt{N}}.
\end{align*}

\subsection{Relationship between analysis on $\pr(\T^d)$ and analysis on $H^{-s}$}
\label{subsec:connectionprhs}

In this paper we work both with the calculus for functions defined on $\pr(\T^d)$ which is commonly used in the setting of mean field control and mean field games, and with the calculus on the space $H^{-s}$ provided by the Hilbertian structure. In this Subsection we discuss both of these notions and their relationship with each other. Much of the analysis in this Subsection will be based on the fact that, as discussed in Subsection \ref{subsec:notation} we have access to the usual Sobolev embedding, and in particular we will use crucially the fact that
\begin{align} \label{embeddings}
H^s \hookrightarrow \linf \text{ if } s > d/2, \quad H^s \hookrightarrow \cC^1 \text{ if } s > d/2 + 1.
\end{align}

Let us mention that we identify a probability measure $m$ with an element of $H^{-s}$ in the usual distributional way: for $f \in H^s$, $s > d/2$, we define $m(f) = \int f dm$, which makes sense thanks to \eqref{embeddings}. In this way we regard $\pr(\T^d)$ as a subspace of $H^{-s}$, and it is easy to check that if $s > d/2$, then $\pr(\T^d)$ is in fact a compact subset of $H^{-s}$. Given $s > 0$ and a function $\phi = \phi(q) : H^{-s} \to \R$ which has a Fr\'echet derivative at a point $q \in H^{-s}$, we use $D_{-s} \phi(q)$ and $\nabla_{-s} \phi(q) : H^{-s} \to H^{-s}$ to denote the derivative and the gradient of $\phi$ at $q$, respectively. That is, $D_{-s} \phi(q)$ is an element of $H^s$ satisfying
\begin{align}
    \phi(p) = \phi(q) + \langle D_{-s} \phi(q) , p - q \rangle_{s,-s} + o(\norm{p-q}_{-s}), \quad 
 p \in H^{-s} , 
 \label{eq:D_{-s}:vs:nabla_{-s}}
\end{align}
while $\nabla_{-s} \phi(q) = \big(D_{-s} \phi(q)\big)^{\ast}$. We say that $\phi \in C^1(H^{-s})$ if $\phi$ is continuously differentiable. 

We would like to understand the relationship between the derivatives $D_{-s} \phi$, $\nabla_{-s} \phi$ of $\phi$ on $H^{-s}$ and the derivatives $\frac{\delta \phi}{\delta m}$, $D_m \phi$ of the restriction of $\phi$ to $\pr(\T^d)$, which we again denote by $\phi$. This is explained in the following Lemma.

\begin{lem} \label{lem:derivativeconnection}
Suppose that $s > d/2 +1$, and that $\phi \in C^1(H^{-s})$. Then $\phi \in \cC^1(\pr(\T^d))$ and
\begin{align*}
    \frac{\delta \phi}{\delta m}(m,y) = D_{-s} \phi(m)(y) - \int_{\T^d} D_{-s} \phi(m)(x) dx, \quad (m,y) \in \pr(\T^d) \times \T^d.
\end{align*}
Moreover, $\frac{\delta \phi}{\delta m}(m,\cdot) \in \cC^1(\T^d)$ for each $m$, and
\begin{align*}
    D_m \phi(m,y) = D \big(D_{-s} \phi(m)\big)(y).
\end{align*}
In addition, $D_m \phi$ is jointly continuous on $\pr(\T^d) \times T^d$.
\end{lem}

\begin{proof}
    If $\phi \in C^1(H^{-s})$ and $m_1,m_2 \in \mathcal{P}(\T^d)$, then 
    \begin{align*}
        \phi(m_2) - \phi(m_1) &= \int_0^1 \langle D_{-s} \phi\big(tm_2 + (1-t) m_1\big), m_2 - m_1 \rangle_{s,-s} dt \\
        &= \int_0^1 \int_{\T^d} D_{-s} \phi\big(tm_2 + (1-t)m_1 \big)(y) (m_2 - m_1)(dy) dt.
    \end{align*}
    To prove the first claim of the lemma, it remains only to check that $\pr(\T^d) \times \T^d \ni (m,y) \mapsto D_{-s} \phi(m) (y)$ is continuous. Uniform Lipschitz continuity in $y$ (uniformly in $m$) follows from the Sobolev embedding $H^{s} \hookrightarrow \cC^{1}(\T^d)$. For continuity in $m$, we use 
    Sobolev embedding again, which says that if $m_k \xrightarrow{k \to \infty} m$ in $d_1$, then 
     $m_k \to m$ in $H^{-s}$ (since $s > d/2 + 1$), so 
    \begin{align*}
        |D_{-s}\phi(m_k)(y) - D_{-s}\phi(m)(y)| \leq \norm{D_{-s}\phi(m_k) - D_{-s} \phi(m)}_{\linf} 
        \leq \norm{D_{-s}\phi(m_k) - D_{-s} \phi(m)}_{s}{,}
    \end{align*}
    which tends to zero by continuity of $D_{-s} \phi$. The argument for the second claim is similar and so is omitted.
\end{proof}

Next, we state a key result, which states how $\cC^{1,1}$ regularity of a function $\phi : q \in H^{-s} \mapsto \phi(q) \in \R$ relates to Lipschitz bounds on the map $m \mapsto D_m \phi$. To state the lemma clearly, it is helpful to introduce the $\cC^{1,1}$ seminorm of 
$\phi$, which is defined by the formula
\begin{equation}
    \big[\phi \big]_{\cC^{1,1}(H^{-s})} = \sup_{q_1 \neq q_2} \frac{\| \nabla_{-s} \phi(q_1) - \nabla_{-s} \phi(q_2)\|_{-s} }{\|q_1 - q_2\|_{-s}}, 
\label{def:seminormc11}
\end{equation}
i.e. $\big[\phi \big]_{\cC^{1,1}(H^{-s})}$ is the Lipschitz constant of the gradient $\nabla_{-s} \phi$.

\begin{prop} \label{prop:dmlipschitz}
    Let $s > d/2 + 1$. Then there is a constant $C$ depending only on $d$ and $s$ such that for any $\cC^1$ map 
   $\phi : q \in H^{-s} \mapsto \phi(q) \in \R$
     and any $m_1,m_2 \in \pr(\T^d)$, $y \in \T^d$, we have 
    \begin{align*}
        |D_m \phi(m_1,y) - D_m \phi(m_2,y)| \leq C \big[ \phi \big]_{C^{1,1}(H^{-s})} d_1(m_1,m_2). 
    \end{align*}
\end{prop}
\begin{proof}
On the one hand, using the Sobolev embedding $H^{s} \hookrightarrow \cC^{1}(\T^d)$
together with Lemma
\ref{lem:derivativeconnection}, 
 we find
\begin{align*}
    |D_m \phi(m_1,y) - D_m \phi(m_2,y)| &= |D \big(D_{-s} \phi(m^1) \big)(y) - D \big(D_{-s} \phi(m_2) \big)(y) | \\
        &\leq \norm{D_{-s} \phi(m_1) - D_{-s} \phi(m_2) }_{\cC^1} \\
        &\leq C\norm{D_{-s} \phi(m_1)  - D_{-s} \phi(m_2) }_{s} \\
        &\leq C\big[ \phi \big]_{\cC^{1,1}(H^{-s})}  \|m_1 - m_2\|_{-s}{.}
\end{align*}
On the other hand, again using the Sobolev embedding $H^s \hookrightarrow \cC^1$, we have $\|m_1 - m_2\|_{-s} \leq C d_1(m_1,m_2)$, and the result follows.
\end{proof}


\color{black}

\section{Properties of the value function}

\label{sec: PropertiesoftheValuefunction}

The goal of this section is to prove regularity properties of the value function $U$ defined in \eqref{def:U} in the $H^{-s}$-regular case (this is Proposition \ref{prop:upropshs}) and in the $d_1$-regular case (this is Proposition {\ref{prop:uregularityd1}}). In this latter case, we use an approximation argument to bypass the existence of linear functional derivatives for the mean field costs, usually required when one wants to use the system of optimality conditions. We also recall some viscosity solutions properties for the value function in Subsection \ref{sec:DPPforvlauefunction}.

\subsection{Optimality conditions}

We start with a standard lemma which characterizes the optimizers of the mean field control problem in terms of a forward-backward system of PDEs.

\begin{lem} \label{lem:optconditions}
    Let Assumption \ref{assump:d1} hold, and assume in addition that $\sF, \sG \in \cC^1(\pr(\T^d))$ with 
    \begin{align*}
        \sup_{m \in \pr(\T^d)} \norm{\frac{\delta \sF}{\delta m}(m,\cdot)}_{\cC^2(\T^d)}  +\sup_{m \in \pr(\T^d)} \norm{\frac{\delta \sG}{\delta m}(m,\cdot)}_{\cC^2(\T^d)} < \infty.  
    \end{align*}
    Then for any $(t_0,m_0) \in [0,T] \times \mathcal{P}(\T^d)$ there is at least one optimizer for the problem \eqref{def:U}, and for any optimizer $(m,\alpha)$ we must have
\begin{enumerate}
\item  $\alpha(t,x) = -D_pH(x,Du(t,x))$ for some solution $(u,m) \in \mathcal{C}^{1,2}([t_0,T) \times \T^d) \times \mathcal{C}([t_0,T], \mathcal{P}(\T^d))$ of
    \begin{align}
        \begin{cases}
            - \partial_t u(t,x)- \Delta_x u(t,x) + H\bigl(x,D_xu(t,x)\bigr) = \frac{\delta \sF}{\delta m}(m_t,x), \quad (t,x) \in [t_0,T) \times \T^d, \\
            \partial_t m(t,x) = \Delta_x m (t,x) + \div_x \Bigl(m D_p H\bigl(x, D_xu(t,x)\bigr)\Bigr), \quad (t,x) \in [t_0,T) \times \T^d, \\
            m_{t_0} = m_0, \quad u(T,x) = \frac{ \delta \sG}{\delta m}(m_T,x), 
        \end{cases}
\label{thesystemofoc29mars2023}
    \end{align}
with the second equation being understood in the weak sense.
\item $\alpha$ satisfies 
$$\|\alpha \|_{\linf} + \sup_{t \in [t_0,T)} \sqrt{T-t}\|D\alpha(t, \cdot) \|_{\linf} \leq C{,}$$
for some $C>0$ depending on $\mathcal{F}$ and $\mathcal{G}$ only through their Lipschitz constants and independent of $t_0$.
\item If in addition Assumption \ref{assump:Hs} is in force, then $u(t,\cdot) :x \in \T^d \mapsto u(t,x) \in \R^d$ belongs to $\mathcal{C}^{s}(\T^d)$ for all $t \in [t_0,T]$ and to $\mathcal{C}^{s+1}(\T^d)$ for all $t \in [t_0,T)$ and we have
\begin{align*}
    \sup_{t \in [t_0,T]} \|\alpha(t,\cdot)\|_{\cC^{s-1}} + \sup_{t \in [t_0,T)} \sqrt{T-t} \|\alpha(t,\cdot) \|_{\cC^{s}} \leq C. 
\end{align*}
\end{enumerate}
\end{lem}

\begin{proof}
    Claim (1) follows easily from standard arguments (see e.g. \cite{BrianiCardaliaguet} Lemma 3.1). Claims (2) and (3) follows easily from Claim (1), thanks to Lemma \ref{lem:linearHJBappendix3avril2023} from the Appendix.
\end{proof}

\subsection{Regularity of the value function}

\subsubsection{The $H^{-s}$-regular case}

The properties of $U$ is the $H^{-s}$-regular case are summarized in the following proposition. 

\begin{prop} \label{prop:upropshs}
    Under Assumption \ref{assump:Hs}, there is a constant $C$, depending only on the various data appearing in the assumption, such that the following holds:
    \begin{enumerate}
    \item The function $U$ is globally $C$-Lipschitz with respect to $H^{-s}$. More precisely, for all $t_1,t_2 \in [0,T]$ and $m_1,m_2 \in \pr(\T^d)$, we have 
    \begin{align*}
    |U(t_1,m_1) - U(t_2,m_2)| \leq C \big( |t_1 - t_2| + \norm{m_1 - m_2}_{-s} \big).
    \end{align*}
    \item For each $t \in [0,T]$, the map $\pr(\T^d) \ni m \mapsto U(t,m)$ is $C$-semi-concave with respect to the norm $\norm{\cdot}_{-s}$. More precisely, for each $t \in [0,T]$, $m_1,m_2 \in \pr(\T^d)$, and $\lambda \in (0,1)$, we have 
\begin{align*}
    U\big(t,(1-\lambda) m_1 + \lambda m_2 \big) \geq (1- \lambda) U(t,m_1) + \lambda U(t,m_2) - \frac{C}{2}\lambda(1 - \lambda) \norm{m_1 - m_2}_{-s}^2. 
\end{align*}
    \end{enumerate}
\end{prop}

To prove Proposition \ref{prop:upropshs} we will need the following stability result.

\begin{lem}
    Assume that $s>d/2$ and $\alpha : (t,x) \in [t_0,T] \times \T^d \mapsto \alpha(t,x) \in \R^d$ satisfies 
    $$\sup_{t_0 \leq t\leq T} \|\alpha(t,\cdot) \|_{s-1} \leq C.$$
Then there is a constant $C'$ depending only on $C$ such that if $(m^i)_{i = 1,2}
\in  \mathcal{C}([t_0,T], \mathcal{P}(\T^d))
$ satisfy in the weak sense
    \begin{align*}
           \partial_t m^i(t,x) = \Delta_{x} m^i(t,x) - \div_{x}\bigl(m^i(t,x) \alpha(t,x)\bigr),
    \end{align*}
    for $(t,x) \in [t_0, T) \times \T^d$, then we have
    \begin{align} \label{fokkerplanck}
        \sup_{t_0 \leq t \leq T} \norm{m^1_t - m^2_t}_{-s} \leq C' \norm{m^1_{t_0} - m^2_{t_0}}_{-s}. 
    \end{align}
\label{stabilityLemmaHs27Mars}
\end{lem}
The proof of Lemma 
\ref{stabilityLemmaHs27Mars}
 is given in Appendix \ref{sec:LinearPDEsappendix}.

\begin{proof}[Proof of Proposition \ref{prop:upropshs}]
\underline{Claim 1.} For regularity in $m$, we fix $t_0 \in [0,T]$ and $m_1,m_2 \in \pr(\T^d)$. Let $\alpha$ be an an optimal control for $(t_0,m_1)$, and let $(m^i_t)_{t_0 \leq t \leq T}$, for $i=1,2$, denote the solutions on $[t_0,T]$ to 
\begin{align*}
  \partial_t m^i(t,x) = \Delta_{x} m^i(t,x) - \div_{x}\bigl(m^i(t,x) \alpha (t,x)\bigr), \quad m^i_{t_0} = m_i. 
\end{align*}
Then we have 
\begin{align*}
    U(t_0,m_1) &= \int_{t_0}^T \Big( \int_{\R^d} L\big(x, \alpha(t,x) \big) m_t^1(dx) + \sF(m_t^1) \Big) dt + \sG(m_T^1) \\
    &=  \int_{t_0}^T \Big( \int_{\R^d} L\big(x, \alpha(t,x) \big) \big(m_t^2 + m_t^1 - m_t^2)(dx) + \sF(m_t^2) + \sF(m_t^1) - \sF(m_t^2) \Big) dt \\ &\qquad+ \sG(m_T^2) + \sG(m_T^1) - \sG(m_T^2)  \\
    &\geq U(t_0,m_2) - C \Big(1 + \int_{t_0}^T \norm{L(\cdot, \alpha(t,\cdot))}_{\mathcal{C}^s} dt \Big)\sup_{t_0 \leq t \leq T} \norm{m^1_t - m_t^2}_{-s} \\
    &\geq U(t_0,m_2) - C \Big(1 + \int_{t_0}^T \frac{1}{\sqrt{T-t}} dt \Big) \sup_{t_0 \leq t \leq T} \norm{m^1_t - m_t^2}_{-s} \\
    &\geq U(t_0,m_2) - C \norm{m_1 - m_2}_{-s}.
\end{align*}
The first inequality uses the fact that, by Lemma \ref{lem:optconditions}, we have 
\begin{align*}
L(x, \alpha(t,x)) = L\Bigl(x, - D_p H\bigl(x,D_{x}u (t,x)\bigr)\Bigr) = - H\bigl(x, D_{x} u(t,x)\bigr) - \alpha(t,x) \cdot D_{x}u(t,x),
\end{align*}
where $(m^1, u)$ is a solution of \eqref{thesystemofoc29mars2023}, and so by Lemma \ref{lem:optconditions} and the regularity of $H$,
$$   \| L(\cdot, \alpha(t,\cdot)) \|_{\cC^{s}} \leq \frac{C}{\sqrt{T-t}}. $$
For the final inequality, we used the bound $\sup_{t_0 \leq t \leq T} \norm{m_t^1 - m_t^2}_{-s} \leq C \norm{m_1 - m_2}_{-s}$ from Lemma \ref{stabilityLemmaHs27Mars}. Up to inverting the roles of $m_0$ and $m_1$, this completes the proof of the Lipschitz regularity in $m$.

For the time regularity we proceed as follows. Given $m_0 \in \mathcal{P}(\T^d)$, $t_0 \in [0,T)$ and $h>0$ such that $t_0+h \in [0,T]$ we can expand
\begin{align} \label{triangleineq}
U(t_0+h,m_0)-U(t_0,m_0) = U(t_0+h,m_0) - U(t_0+h, m^0_{t_0+h}) + U(t_0+h,m^0_{t_0+h}) - U(t_0,m_0),
\end{align}
where $(m^0_t)_{t_0 \leq t  \leq T}$ is an optimal trajectory for $U(t_0,m_0)$ with control $\alpha$. On the one hand, by dynamic programming, using that $\alpha$ is bounded independently from $(t_0,m_0)$ and the boundness of $\mathcal{F}$, we can argue that 
\begin{align*}
    \big|U(t_0+h,m^0_{t_0+h}) - U(t_0,m_0)    \big|&\leq \left|\int_{t_0}^{t_0+h} \int_{\T^d} L(x,\alpha(t,x))m^0_t(dx)dt + \int_{t_0}^{t_0+h} \mathcal{F}(m^0_t)dt \right| 
    \leq Ch
\end{align*}
for some $C$ depending only the growth properties of $L$ and ${\mathcal F}$.
On the other hand, by Lischitz regularity of $U(t,\cdot)$ over $\mathcal{P}(\T^d)$ with respect to $\|.\|_{-s}$ we can argue that 
    $$|U(t_0,m^0_{t_0+h}) - U(t_0,m_0) | \leq C \norm{m^0_{t_0+h} - m_0 }_{-s}.$$
Now, using the Sobolev embedding $H^s \hookrightarrow \cC^2$ (recall that $s > d/2+2$ here), we have 
\begin{align*}
    \int_{\T^d} \phi(x) (m^0_{t_0+h} - m_0)(dx) &= \int_{t_0}^{t_0+h} \int_{\T^d} \left[ D\phi(x) \cdot \alpha(t,x) + \Delta \phi(x) \right] m^0_t(dx) dt \leq Ch \|\phi \|_{s}, 
\end{align*}
so that $\|m^0_{t_0 + h} - m_0\|_{s} \leq Ch$,
which, in light of \eqref{triangleineq}, completes the proof that $U$ is Lipschitz in time.

\underline{Claim 2.} We now turn to the semi-concavity of $U$.
We fix $m_1,m_2 \in \mathcal{P}(\T^d)$, $t_0 \in [0,T]$ and $\lambda \in (0,1)$ and we consider $\alpha$ to be an optimal control for $U(t_0,m_{\lambda})$ with $m_{\lambda} =(1-\lambda)m_1 + \lambda m_2$. By Lemma \ref{lem:optconditions}, we know that $\alpha(t,\cdot)$ is bounded in $H^{s-1}$, uniformly  in $t \in [t_0,T]$. This time we define $(m^1_t)_{t_0 \leq t \leq T}$, $(m^2_t)_{t_0 \leq t \leq T}$ and $(m^{\lambda}_t)_{t_0 \leq t \leq T}$ to be the solutions to 
\begin{equation}
\partial_t m(t,x) = \Delta_{x}m(t,x) - \div_{x} \bigl(m(t,x)\alpha(t,x) \bigr), \hspace{15pt} \mbox{in } [t_0,T) \times \T^d,
\label{Fpe12avril2023}
\end{equation}
starting respectively from $m_1$, $m_2$ and $m_{\lambda}$. By linearity of the Fokker-Planck equation \eqref{Fpe12avril2023} it holds for all $t \in [t_0,T]$,
\begin{equation}
    m^{\lambda}_t = (1-\lambda) m^1_t + \lambda m^2_t.
\label{LinearityFPE10/02/2023}
\end{equation}    
By optimality of $\alpha$ for $U(t_0,m^{\lambda}_{t_0})$ we have

\begin{align*}
    (1-\lambda) U(t_0,m_1) &+ \lambda U(t_0,m_2) - U(t_0, m_{\lambda}) \\
    &\leq (1-\lambda) \int_{t_0}^T \int_{\R^d} L(x,\alpha(t,x))m^1_t(dx)dt + \lambda  \int_{t_0}^T \int_{\R^d} L(x,\alpha(t,x))m^2_t(dx)dt \\
    & \quad -  \int_{t_0}^T \int_{\R^d} L(x,\alpha(t,x))m^\lambda_t(dx)dt \\
    & \quad + \int_{t_0}^T \left[(1-\lambda)\mathcal{F}(m^1_t) + \lambda \mathcal{F}(m^2_t) - \mathcal{F}(m^{\lambda}_t) \right]dt \\
    &\quad + (1-\lambda)\mathcal{G}(m^1_T) + \lambda \mathcal{G}(m^2_T) - \mathcal{G}(m^{\lambda}_T).
\end{align*}
The terms involving the running cost $L$ cancel out since $m^{\lambda}_t$ satisfies \eqref{LinearityFPE10/02/2023}. For the remaining terms, we use the semi-concavity of $\mathcal{F}$ and $\mathcal{G}$ as well as Lemma 
\ref{stabilityLemmaHs27Mars} to deduce that, for some constant $C>0$ depending on the semi-concavity constants of $\mathcal{F}$ and $\mathcal{G}$ as well as $\sup_{t_0 \leq t \leq T} \|\alpha(t,\cdot) \|_{s-1}$ (and independent of $t_0$),
\begin{align*}
    (1-\lambda) &U(t_0,m_1) + \lambda U(t_0,m_2) - U(t_0, m_{\lambda}) \\
    &\leq \frac{C}{2} \lambda(1-\lambda) \sup_{t \in [t_0,T]} \|m_t^1 - m_t^2 \|^2_{-s} \leq \frac{C}{2} \lambda(1-\lambda) \|m_1 - m_2 \|^2_{-s},
\end{align*}
which concludes the proof of the proposition.
    
\end{proof}

\subsubsection{The $d_1$-regular case}

We now investigate the regularity of $U$ under Assumption \ref{assump:d1}. We will use the notation
\begin{equation}
    \norm{m^2 -m^1}_{-2,\infty} := \sup_{\|\phi\|_{\cC^2} \leq 1} \int_{\R^d} \phi(x) (m^2-m^1)(dx)
\label{normW2inftystar}
\end{equation}
for any two probability measures $m^1,m^2 \in \sP(\T^d).$ 
The regularity properties of $U$ can be summarized as follows. 

\begin{prop} \label{prop:uregularityd1}
      Under Assumption \ref{assump:d1}, there is a constant $C$, depending only on the various data appearing in the assumption, such that the following holds:
       \begin{enumerate}
        \item $U$ is $d_1$-Lipschitz with respect to $m$ with constant $C$, i.e. for each $t \in [0,T]$ and each $m_1,m_2 \in \pr(\T^d)$, we have 
        \begin{align*}
        |U(t, m_1) - U(t,m_2) | \leq C d_1(m_1,m_2).
        \end{align*}
        \item $U$ is locally $\|\cdot\|_{-2,\infty}$-Lipschitz: it satisfies, for each $t \in [0,T)$ and each $m_1,m_2 \in \pr(\T^d)$, the estimate 
        \begin{align}
        |U(t, m_1) - U(t,m_2) | \leq \frac{C}{\sqrt{T-t}} \norm{m_2-m_1}_{-2,\infty}.
        \label{W2inftyregu}
        \end{align} 
        \item $U$ is $d_1$-semiconcave with respect to $m$ with constant $C$, i.e. for each $t \in [0,T]$, $\lambda \in (0,1)$, and each $m_1,m_2 \in \pr(\T^d)$, we have
        \begin{align*}
            U(t,(1-\lambda)m_1 + \lambda m_2) \geq (1 - \lambda ) U(t, m_1) + \lambda U(t,m_2) - \frac{C}{2} \lambda (1 - \lambda) d_1^2(m_1,m_2).
        \end{align*}
        \item $U$ is H\"older continuous and locally Lipschitz continuous in time: 
         it satisfies 
        \begin{align*}
            |U(t + h, m) - U(t,m)| \leq C \Bigl( \sqrt{h} \wedge \frac{h}{\sqrt{T-h-t}}\Bigr),
        \end{align*}
        for each $t \in [0,T]$, $h \in (0,T-t)$, and $m \in \pr(\T^d)$. 
    \end{enumerate}
\end{prop}

The rest of this subsection is devoted to proving Proposition \eqref{prop:uregularityd1}. Assumption \ref{assump:d1} is in force, and we also  
assume that $(\mathcal{F}^{n})_{n \geq 0}$ and $(\mathcal{G}^{n})_{n \geq 0}$ are sequences of smooth approximations of $\mathcal{F}$ and $\mathcal{G}$ such that, 

\begin{enumerate}
    \item[($n$.1)] for all $n \geq 0$, $\mathcal{F}^n$ and $\mathcal{G}^n$ are Lipschitz and semi-concave with Lipschitz and semi-concavity constants bounded independently from $n.$
    \item[($n$.2)] $\mathcal{F}^n \rightarrow \mathcal{F}$ and $\mathcal{G}^n \rightarrow \mathcal{G}$ uniformly over $\mathcal{P}(\T^d)$ as $n \rightarrow +\infty$.
    \item[($n$.3)] For all $n \geq 0$, $\mathcal{F}^n, \mathcal{G}^n \in \mathcal{C}^1(\mathcal{P}(\T^d))$ with $\frac{\delta \mathcal{F}^n}{\delta m}(m,\cdot), \frac{\delta \mathcal{G}^n}{\delta m}(m,\cdot) \in \mathcal{C}^2(\T^d)$ for all $m \in \T^d$ with jointly continuous derivatives.
\end{enumerate}
The existence of such approximations is proved in Lemma \ref{lem:regularisation} below. For all $n \geq 0$ we define $U^n$ to be the value function in \eqref{def:U} with costs $\mathcal{F}^n$ and $\mathcal{G}^n$ instead of $\mathcal{F}$ and $\mathcal{G}$. It is plain to check from the control formulation that $U^n$ converges uniformly to $U$ over $[0,T] \times \mathcal{P} (\T^d)$ as $n \rightarrow +\infty$ (thanks to property ($n$.2) right above). 

The proof of Proposition \ref{prop:uregularityd1} will make use of the following stability result,
\begin{lem}
   Suppose that $\alpha : (t,x) \in [t_0,T] \times \T^d \mapsto \alpha(t,x) \in \R^d$ 
   is differentiable in $x$
   on $[t_0,T) \times {\mathbb R}^d$ and
   satisfies
\begin{equation}
    \sup_{t_0 \leq t  \leq T} \|\alpha(t,\cdot) \|_{\linf} +\sup_{t_0 \leq t < T} \sqrt{T-t} \|D \alpha(t,\cdot) \|_{\linf} \leq C,
\label{aprioriestimateoptimalcontrol27mars}
\end{equation}
for some $C>0$, then, there is $C'$ depending only on $C$ such that if $(m^{i})_{i=1,2}
\in  \mathcal{C}([t_0,T], \mathcal{P}(\T^d))
$ satisfy in the weak sense
    \begin{align*}
           \partial_t m^i(t,x) = \Delta_{x} m^i(t,x) - \div_{x}\bigl(m^i(t,x) \alpha(t,x)\bigr),
           \end{align*}
    for all $(t,x) \in [t_0,T) \times \T^d$, then we have
\begin{equation}
\sup_{t\in [t_0,T]} d_1(m_t^1,m_t^2) \leq C' d_1(m_{t_0}^1,m_{t_0}^2),
\label{d1stability8mai}
\end{equation}
and, using notation \eqref{normW2inftystar},
\begin{equation}\sup_{t \in [t_0,T]} \sqrt{t-t_0} d_1(m_t^1,m_t^2) \leq C' \norm{m^2_{t_0}-m^{1}_{t_0}}_{-2,\infty}.
\label{W2inftystartstability8mai}
\end{equation}
\label{lem:stabilityEstimate20Mars2023}
\end{lem}
The proof is given in Appendix \ref{sec:LinearPDEsappendix}. We now prove the main result of the subsection, Proposition \ref{prop:uregularityd1}.


\begin{proof}[Proof of Proposition \ref{prop:uregularityd1}]
\underline{Claim 1:} We first show the analog result for $U^n$. Fix $t_0 \in [0,T]$ and let $(m_t^{i})_{i=1,2}$ be solutions to 
$$\partial_t m^{i}(t,x) -\Delta_{x} m^{i}(t,x) + \div_{x} (\alpha^n(t,x) m^{i}(t,x))=0, \quad m_{t_0}^i = m_i,$$
as in the statement of Lemma \ref{lem:optconditions}, where $\alpha^n$ is an optimal control for $U^n(t_0,m_1)$. Thanks to Lemma \ref{lem:optconditions} and to Property  ($n$.1),  we know that $\alpha^n$ satisfies \eqref{aprioriestimateoptimalcontrol27mars} for some $C$ independent from $n$. By optimality of $\alpha$ for $U^n(t_0,m_1)$ and Lipschitz regularity of $\mathcal{F}^n$ and $\mathcal{G}^n$, inequality 
\begin{align}
\notag   &U^n(t_0,m_1) - U^n(t_0,m_2) \\
\notag    &\quad \geq \int_{t_0}^T \int_{\T^d} L\bigl(x,\alpha^n(t,x)\bigr)(m^1_{t} - m^2_{t})(dx)dt + \int_{t_0}^T \big(\mathcal{F}^n(m^1_t) - \mathcal{F}^n(m^2_t) \big) dt +\mathcal{G}^n(m^1_T) - \mathcal{G}^n(m^2_T) \\
    &\quad \geq -C \int_{t_0}^T 
    \bigl(1 + \|D\alpha^n(t,\cdot) \|_{\linf}) \bigr)d_1(m^1_t,m_t^2) dt - C \int_{t_0}^T d_1(m^1_t,m^2_t) - Cd_1(m^1_T,m^2_T) \label{computation8mai} \\
\notag    &\quad \geq -C \int_{t_0}^T 
    \bigl(1 + \|D\alpha^n(t,\cdot) \|_{\linf}) \bigr)d_1(m^1_t,m_t^2) dt - C \sup_{t_0\leq t \leq T} d_1(m^1_t,m^2_t) \\
\notag    &\quad \geq -C \Bigl(1 + \int_{t_0}^T \|D\alpha^n (t,\cdot) \|_{\linf} dt \Bigr) \sup_{t \in [t_0,T]} d_1(m_t^1,m^2_t) 
\end{align}
holds for some $C$ independent of $n$ (and of $t_0$, $m_1$, $m_2$). Using \eqref{d1stability8mai} in Lemma 
 {\ref{lem:stabilityEstimate20Mars2023}} and then letting $n \rightarrow +\infty$, we conclude that
$$ U(t_0,m_2) - U(t_0,m_1) \leq  Cd_1(m_{t_0}^1,m_{t_0}^2)= Cd_1(m_1,m_2).$$
Reversing the roles of $m_1$ and $m_2$, we conclude that $U$ is Lipschitz continuous in the measure variable. 

\underline{Claim 2:} Using \eqref{W2inftystartstability8mai} from Lemma \ref{lem:stabilityEstimate20Mars2023} instead of \eqref{d1stability8mai} in \eqref{computation8mai} leads to Claim 2. 

\underline{Claim 3:} To prove the semi-concavity, we argue as in the $H^{-s}$-regular case, using the stability of the Fokker-Planck equation given by Lemma \ref{lem:stabilityEstimate20Mars2023}. This allows one to show that there is some $C>0$ independent of $n \geq 0$ such that, for all $m_1,m_2 \in \mathcal{P}(\T^d)$ and all $\lambda \in (0,1)$, it holds 
$$U^n \bigl(
t,(1-\lambda)m_1+\lambda m_2
 \bigr) \geq (1-\lambda) U^n  (
t_0,m_1  ) + \lambda U^n   ( t_0,m_2 ) - \frac{C}{2}\lambda(1-\lambda) d_1^2(m_1,m_2).$$
Using the convergence of $U^n$ toward $U$, we conclude that $U$ is semi-concave in the measure variable.

\underline{Claim 4:} We start by the analog statement for $U^n$ for $n \geq 0$. Let $m \in \mathcal{P}({\mathbb T}^d)$, $t \in [0,T)$ and $h \in (0,T-t]$. Following the same argument as in the proof of the time Lipschitz regularity in Proposition \ref{prop:upropshs} and using \eqref{W2inftyregu}, we infer that
\begin{equation}
     |U^n(t+h,m) - U^n(t,m)| \leq \frac{C h}{\sqrt{T-t-h}}.
\label{TimeLipshitz8mai}
\end{equation}
The global Hölder regularity in time is, on the other hand, standard. Once we have \eqref{TimeLipshitz8mai} it is straightforward, simply by noticing that the latter implies that, for all $m \in \pr(\T^d)$, $t \mapsto U(t,m)$ is absolutely continuous over $[0,T)$ with
$$\left| \frac{d}{dt}U^n(t,m) \right| \leq \frac{C}{\sqrt{T-t}}, \quad \quad \mbox{ for almost all } t\in [0,T),$$
and therefore, integrating in time this leads, for all $m\in \pr(\T^d)$, all $t \in [0,T)$ and all $h \in [0,T-t]$ to
$$|U^n(t+h,m) - U^n(t,h)| \leq C \sqrt{h}.$$
Passing to the limit when $n \rightarrow +\infty$ gives the result. 
\end{proof}

\subsection{Viscosity (sub)-solution property}

\label{sec:DPPforvlauefunction}

Here we record the fact that $U$ is a viscosity solution of the equation \eqref{hjbinf}. In fact, we only need the sub-solution property, and we only need to consider test functions $\phi \in \mathcal{C}^1((0,T) \times H^{-s})$.

\begin{lem} \label{lem:usubsol}
Suppose that $\phi \in \mathcal{C}^1((0,T) \times H^{-s})$ touches $U$ from above at $(t_0,m_0) \in (0,T) \times \mathcal{P}(\T^d)$ that is 
\begin{equation}
    U(t_0,m_0) -\phi(t_0,m_0) = \sup_{(t,m) \in (0,T) \times \mathcal{P}(\T^d)} \Bigl[U(t,m) - \phi(t,m)
    \Bigr],
\end{equation}
then, it holds that 
$$-\partial_t \phi(t_0,m_0) - \int_{\R^d} \Delta_x \frac{\delta \phi}{\delta m}(t_0,m_0,x)m_0(dx)+\int_{\R^d} H\Bigl(x,D_{x} \frac{\delta \phi}{\delta m}(t_0,m_0,x)\Bigr)m_0(dx)  \leq \mathcal{F}(m_0).$$

\end{lem}

The proof is standard, up to using Lemma \ref{lem:derivativeconnection} to check that tests functions $\phi \in \mathcal{C}^1((0,T) \times H^{-s})$ are regular enough for the usual argument based on the chain rule for flows of probability measures to apply. We refer for instance to \cite{CarmonaDelarue_book_II} section 4.4.3 for a proof and the corresponding supersolution property.

\section{Three regularization procedures}
\label{sec:regularization}

In this section, we present and analyze three methods for regularizing a function 
$\Phi : m \in \pr(\T^d) \mapsto \Phi(m) \in \R$.

\subsection{Regularization by mollification of the Fourier coefficients} \label{subsec:mollificationfourier}

\begin{lem}
\label{lem:regularisation}
For two constants $c_1$ and $c_2$, let $\Phi : {\mathcal P}({\mathbb T}^d) \rightarrow {\mathbb R}$ be 
$c_1$-Lipschitz continuous and $c_2$-semi-concave with respect to the $d_1$-Wasserstein distance. 
Then, there exists a sequence of functions $(\Phi^n : {\mathcal P}({\mathbb T}^d) \rightarrow {\mathbb R})_{n \geq 1}$ such that 
\begin{enumerate}
\item $(\Phi^n)_{n \geq 1}$ converges to $\Phi$, uniformly on ${\mathcal P}({\mathbb T}^d)$;
\item for any $n \geq 1$, $\Phi^n$ is $c_1$-Lipschitz continuous and $c_2$-semi-concave with respect to $d_1$;
\item for any $n \geq 1$, $\Phi^n$ is continuously differentiable with respect to $m$
(for $d_1$) and 
its derivative $\delta \Phi^n/\delta m$ (defined on ${\mathcal P}({\mathbb T}^d) \times {\mathbb T}^d$)
is jointly continuous (with 
${\mathcal P}({\mathbb T}^d)$ being equipped with $d_1$)
and has jointly continuous derivatives of any order  in the $x$-variable.
\end{enumerate}
\end{lem}

The proof of the above result is a direct consequence of a mollification argument  explained in
Proposition 3.14
in \cite{cecchin2022weak}.
For completeness, we present a sketch of the proof in Appendix \ref{subse:mollif:arg}.

\subsection{Regularization by mollification of the measure argument}

\label{subsec:mollificationkernel}

\begin{lem}
\label{lem:convolution:rhodelta}
For two constants $c_1$ and $c_2$, let $\Phi : {\mathcal P}({\mathbb T}^d) \rightarrow {\mathbb R}$ be
$c_1$-Lipschitz continuous and $c_2$-semi-concave with respect to the $d_1$-Wasserstein distance. 
For a smooth symmetric density $\rho$ on ${\mathbb R}^d$ with compact support and for 
$\delta >0$, let
\begin{equation*}
\rho_\delta(x) := \frac1{\delta^d} \rho\bigl( \frac{x}{\delta} \bigr), \quad x \in {\mathbb R}^d, 
\end{equation*}
and then
\begin{equation*}
\Phi^{\delta}(m) := \Phi \bigl(m*\rho_\delta \bigr), \quad m \in {\mathcal P}({\mathbb T}^d). 
\end{equation*} 
Then, 
for any $s>0$, there exists a constant 
$\Gamma(d,\rho,s) \geq 0$ such that
\begin{enumerate}
\item 
$\sup_{m \in \pr(\T^d)} \vert \Phi^\delta(m) - \Phi(m) \vert \leq c_1 \Gamma(d,\rho,s) \delta$.
\item
$\Phi^{\delta}$ is $\Gamma(d,\rho,s) c_1 \delta^{-(s-1)}$-Lipschitz continuous and 
$\Gamma(d,\rho,s) c_2 \delta^{-2(s-1)}$-semi-concave with respect to 
$\| \cdot \|_{-s}$.
\end{enumerate}
Lastly, if $\Phi$ is continuous differentiable with respect to $m$, so is $\Phi^\delta$ and  the derivative is given by
\begin{equation*}
\frac{\delta \Phi^\delta}{\delta m}(m,x) = 
\Bigl( \frac{\delta \Phi}{\delta m}{\bigl(m*\rho_{\delta},\cdot\bigr)} * \rho_\delta\Bigr)(x), \quad m \in {\mathcal P}({\mathbb T}^d),
\ x \in {\mathbb T}^d.
\end{equation*}
 \end{lem}

 \begin{proof}
{ \ }
Let us first recall that there is a constant $\Gamma(d,\rho,s)$ such that for $f \in L^2$, 
\begin{align*}
    \| f * \rho_{\delta}\|_{s} \leq \Gamma(d,\rho,s) \delta^{-(s-1)} \|f\|_1.
\end{align*}
This fact can easily be proven via Young's convolution inequality when $s \in \N$, and can be proved in a straightforward way via Fourier analysis when $s \notin \N$, so we omit the details. We now proceed in several steps.
\newline \newline 
\underline{Step 1}. 
The first claim in the proof follows from the Lipschitz property of 
$\Phi$ under $d_1$ together with the fact that, for any 
$m \in \pr(\T^d)$, 
\begin{equation*}
d_1\bigl(m*\rho_\delta,m \bigr)
\leq \sup_{f 1-\textrm{Lip}}
\int_{{\mathbb T}^d} \bigl[ f *\rho_\delta(x) - f(x)\bigr] m(dx) 
\leq \delta. 
\end{equation*}

Similarly, the second claim in the proof follows from 
the fact that, for any two $m,m' \in {\mathcal P}({\mathbb T}^d)$, 
\begin{equation*} 
\begin{split}
d_1\bigl(m*\rho_\delta,m'*\rho_\delta\bigr)
&\leq 
\sup_{f 1-\textrm{Lip}}
\int_{{\mathbb T}^d} f *\rho_\delta(x)  \bigl( m-m'\bigr)(dx)
\\
&\leq 
\| m - m' \|_{-s}
\sup_{f 1-\textrm{Lip} : \widehat{f}^0=0}
\| f * \rho_\delta \|_{s} \\
&\leq  \| m - m' \|_{-s} \frac{\Gamma(d,\rho,s) }{\delta^{s-1}} 
\sup_{f 1-\textrm{Lip} : \widehat{f}^0=0}
\| f \|_{1} \\
&\leq \| m - m' \|_{-s} \frac{\Gamma(d,\rho,s) }{\delta^{s-1}}.
\end{split}
\end{equation*}
which completes the proof of the first part.
\newline 
\noindent \underline{Step 2}. The claim related with semi-concavity is shown in a somewhat similar manner. Indeed, we recall (once again)
 that
the semi-concavity of $\Phi$ may be expressed as 
\begin{equation*}
\Phi \bigl( \lambda m'+ (1-\lambda) m \bigr) 
\geq 
\lambda \Phi(m') +
 (1-\lambda)
 \Phi(m) - \frac{c_2}{2} \lambda (1-\lambda) d_1(m,m')^2. 
 \end{equation*}
 Replacing $m$ by $m*\rho_\delta$ and then 
 $m'$ by $m'*\rho_\delta$, we 
 see that the only difficulty is to upper bound 
 $d_1(m*\rho_\delta,m'*\rho_\delta)^2$ by $\|m-m' \|_{-s,2}^2$, up to a multiplicative constant, but this is exactly 
 what is done in the first step. 
\newline \newline
\underline{Step 3.} The formula for the derivative may be found in \cite[Chapter 5]{CarmonaDelarue_book_I}. 
\end{proof}

\subsection{Regularization by sup-convolution in $H^{-s}$}

\label{subsec:supconv}

We now turn to another regularization method that is key in our approach as it is shown to preserve the viscosity sub-solution property (see Proposition \ref{lem:uepsproperties} below).

\begin{prop} \label{prop:supconvproperties} Let $\Phi: \mathcal{P}(\T^d) \to \R$ 
be a function which is $C_L$-Lipschitz and $C_S$-semiconcave with respect to $H^{-s}$. For $\epsilon > 0$, let $\Phi^{\epsilon} : H^{-s} \to \R$ be defined by 
\begin{align}
\Phi^{\epsilon}(q) \coloneqq \sup_{m \in \pr(\T^d)} \Big\{ \Phi(m) - \frac{1}{2 \epsilon} \|q - m\|_{-s}^2 \Big\}.
\label{eq:sup:convolution:4.1}
\end{align}
Then for all $\epsilon < \frac{1}{2 C_S}$, we have 
\begin{enumerate}
    \item For all $m \in \mathcal{P}(\T^d)$,
$$ 0 \leq \Phi^{\epsilon}(m) - \Phi(m) \leq 2C_L^2 \epsilon.$$

    \item $\Phi^{\epsilon} \in \mathcal{C}^{1}(H^{-s})$, and we have (using the same notation as in \eqref{def:seminormc11})
    \begin{align} \label{c11bound}
    \big[\Phi^{\epsilon}\big]_{\cC^{1,1}(H^{-s})} \leq  \Big(\frac{1}{\epsilon} \vee (2C_S) \Big).
    \end{align}
    \item For $m \in \pr(\T^d)$, we have 
        (using the same notation as in \eqref{eq:D_{-s}:vs:nabla_{-s}})
    \begin{align*}
        \nabla_{-s} \Phi^{\epsilon}(m) = \frac{1}{\epsilon} \big(m_{\epsilon} - m\big), 
    \end{align*}
    where $m_{\epsilon}$ is the unique element of $\pr(\T^d)$ such that 
    \begin{align*}
        \Phi^{\epsilon}(m) = \Phi(m_{\epsilon}) - \frac{1}{2 \epsilon} \norm{m_{\epsilon} - m}_{-s}^2.
    \end{align*}
    Moreover, we have $\norm{m - m_{\epsilon}}_{-s} \leq 2 C_L \epsilon$.

     \item For all $m_1, m_2 \in \pr(\T^d)$, we have
$$ |\Phi^{\epsilon}(m_1)- \Phi^{\epsilon}(m_2) | \leq 2C_L \|m_1- m_2 \|_{-s}.$$

\end{enumerate}
\end{prop}

\begin{proof}
\underline{Claim \textit{(1)}.}
The fact that $\Phi^{\epsilon} \geq \Phi$ on $\pr(\T^d)$ is clear from the definition. For the other inequality, let us fix $m \in \mathcal{P}(\T^d)$ and let $m^{\epsilon} \in \mathcal{P}(\T^d)$ be a maximizer in the definition of $\Phi^{\epsilon}$, i.e.
\begin{align} \label{uepsdef}
    \Phi^{\epsilon}(m) = \Phi(m^{\epsilon}) - \frac{1}{2\epsilon}\|m-m^{\epsilon}\|_{-s}.
\end{align}
Notice that a maximizer does exist since ${\mathcal P}({\mathbb T}^d)$ is compact for the weak topology and the function 
$\Phi(\cdot)- \| \cdot - m\|_{-s}^2/2$ is, for a fixed $m$, continuous for the weak topology
(which is easy to check since the Fourier coefficients are continuous for the weak topology and $s$ is here assumed to be greater than $d/2$).  

Since $\Phi^{\epsilon}(m) \geq \Phi(m)$, rearranging \eqref{uepsdef} gives
\begin{align*}
\frac{1}{2\epsilon} \|m-m^{\epsilon} \|_{-s}^2 &\leq \Phi(m^{\epsilon}) - \Phi(m) \leq  C_L \|m-m^{\epsilon} \|_{-s}.
\end{align*}
We deduce that 
$$ \|m-m^{\epsilon} \|_{-s} \leq 2C_L\epsilon. $$
Coming back to the definition 
\eqref{uepsdef} of $m^{\epsilon}$ we get 
\begin{align*}
    \Phi^{\epsilon}(m) - \Phi(m) &= \Phi(m^{\epsilon}) - \Phi(m) - \frac{1}{2\epsilon} \|m-m^{\epsilon} \|_{-s}^2 \\
    &\leq \Phi(m^{\epsilon}) - \Phi(m) \leq C_L \|m - m^{\epsilon}\|_{-s} \leq 2 C_L^2 \epsilon. 
\end{align*}
\newline \newline
\underline{Claim \textit{(2)}.}
We first notice that, for every $m\in \mathcal{P}(\T^d)$, the map 
\begin{align*}
    q \mapsto \Phi(m) -\frac{1}{2\epsilon} \|q-m \|_{-s}^2 + \frac{1}{2\epsilon} \|q\|_{-s}^2 
\end{align*}
is linear in $q$ (which is precisely where the Hilbertian structure of $H^{-s}$ comes in) and therefore 
\begin{align*}
q \mapsto \Phi^{\epsilon}(q) + \frac{1}{2\epsilon} \|q\|_{-s}^2 = \sup_{m \in \pr(\T^d)} \Big\{ \Phi(m) -\frac{1}{2\epsilon} \|q-m \|_{-s}^2 + \frac{1}{2\epsilon} \|q\|_{-s}^2  \Big\}
\end{align*}
is convex. That is, $\Phi^{\epsilon}$ is $\frac{1}{\epsilon}$ semi-convex over $H^{-s}$ for any $\epsilon > 0$.

On the other hand, it is straightforward to check that 
$$ (m,q) \mapsto \Phi(m) - \frac{1}{2\epsilon} \|m-q\|^2_{-s} - C_S \|q\|^2_{-s} $$
is concave over $\mathcal{P}(\T^d) \times H^{-s}$ as soon as $\epsilon \leq \frac{1}{2C_S}$. In this case, a small adaptation of a lemma in  \cite{Lasry1986} is enough to deduce that $\Phi^{\epsilon} - C_S \|\cdot \|^2_{-s}$ is concave. That is, for $\epsilon \leq \frac{1}{2C_S}$, we have that $\Phi^{\epsilon}$ is $2C_S$ semi-concave As a consequence, and as explained in \cite{Lasry1986}, \eqref{c11bound} follows. In particular, $\Phi^{\epsilon}$ is continuously differentiable.
\newline \newline 
\underline{Claim \textit{(3)}.} The fact that $\nabla_{-s} \Phi^{\epsilon}(m) = \frac{1}{\epsilon}\big(m_{\epsilon} - m\big)$ is straightforward since, by definition of $m^{\epsilon}$, $m$ is a maximum of
$H^{-s} \ni q \mapsto \Phi(m^{\epsilon}) - \frac{1}{2\epsilon}\|q-m^{\epsilon} \|_{-s}^2 - \Phi^{\epsilon}(q)$ (the function is non-positive and matches $0$ at $m$). The estimate on $m_{\epsilon} - m$ was already obtained in the proof of Claim \textit{(1)}. 
\newline \newline 
\underline{Claim \textit{(4)}.} Lipschitz regularity for $\Phi^{\epsilon}$ follows easily from Claims \textit{(2)} and \textit{(3)}, since they show that $\Phi^{\epsilon}$ is $\cC^{1}$ on $H^{-s}$ and its gradient satisfies $\|\nabla_{-s} \Phi^{\epsilon}\|_{-s} = \frac{1}{\epsilon} \|m_{\epsilon} - m\|_{-s} \leq 2C_L$ on the convex subset $\pr(\T^d)$ of $H^{-s}$.

\end{proof}

One specificity of the sup-convolution operation introduced
in 
\eqref{eq:sup:convolution:4.1}
is that the 
supremum is just taken on a tiny subset of the Hilbert space 
$H^{-s}$. Indeed, it is easy to check that 
${\mathcal P}({\mathbb T}^d)$ has an empty interior w.r.t. to 
the $H^{-s}$ norm. Quite surprisingly, this creates substantial difficulty 
to study the regularity properties of $\Phi^\epsilon$ that may be possibly 
inherited from those of $\Phi$. When the maximization supporting the sup-convolution is taken over the whole space, there is no difficulty for, say, transferring any uniform continuity property from 
$\Phi$ (w.r.t. any arbitrary topology, possibly different from the $H^{-s}$-topology) 
to $\Phi^\epsilon$. Basically, a mere linear change of variable in the definition of the sup-convolution suffices. This fact is well reported in the paper \cite{Lasry1986}. 
However, things become more subtle when the argument in the sup-convolution is taken in 
a strict subset of the whole Hilbert space (here $H^{-s}$), as it is the case in 
\eqref{eq:sup:convolution:4.1}. Typically, issues may arise if maximizers are located at the boundary (here the boundary of ${\mathcal P}({\mathbb T}^d)$), which becomes very likely if the set over which the maximization is performed is tiny. 
This is the purpose of the next proposition to address this problem and to give sufficient conditions on the measure $m$ such that $m_\epsilon$ is `inside' the space of probability measures (for some topology) and to deduce subsequently further regularity properties on $\Phi^\varepsilon$ in the neighborhood of $m$. One additional subtlety is that we do so when $\Phi$ is in fact replaced by $\Phi^\delta$, as given by the previous 
Lemma \ref{lem:convolution:rhodelta}.

\begin{prop}
\label{prop:threshold}
For two constants $c_1$ and $c_2$, let $\Phi : {\mathcal P}({\mathbb T}^d) \rightarrow {\mathbb R}$ be
$c_1$-Lipschitz continuous and $c_2$-semi-concave with respect to the $d_1$-Wasserstein distance. 

For a given $s>0$ and for 
any $\delta,\epsilon>0$, let (with the same notation as in the statement of 
Lemma 
\ref{lem:convolution:rhodelta}
and, in particular, with the same choice for $\rho$), 
\begin{equation*}
\Phi^{\delta,\epsilon}(m_0) : = \sup_{ m \in {\mathcal P}({\mathbb T}^d) } 
\Big\{ \Phi^\delta(m) - \frac1{2 \epsilon} \| m - m_0 \|^2_{-s} \Bigr\}, 
\quad m_0 \in {\mathcal P}({\mathbb T}^d). 
\end{equation*} 
Then, 
for any $\eta>0$, there exists a constant $\gamma(d,\rho,s,\eta)>0$ such that, 
whenever 
$\delta^{2s-2}> 2c_2 \Gamma(d,\rho,s) \epsilon$ 
(with $\Gamma(d,\rho,s)$ as in the statement of Lemma 
\ref{lem:convolution:rhodelta})
and under the condition 
\begin{equation*} 
m_0  \geq  \frac{c_1  \epsilon \gamma(d,\rho,s,\eta)}{\delta^{(2s+d/2+\eta-1)}} \textrm{\rm \leb}, 
\end{equation*}
the (unique) maximizer 
$m_{\delta,\epsilon}$ in the definition of $\Phi^{\delta,\epsilon}$ satisfies 
\begin{equation*} 
\bigl\| {m}_{\delta,\epsilon} - {m}_0
\bigr\|_{L^\infty}  
\leq \frac{c_1 \gamma(d,\rho,s,\eta) \epsilon}{\delta^{(2s+d/2+\eta-1)}}.
\end{equation*} 
Moreover, 
\begin{equation*} 
 \Bigl\| D_x \frac{\delta \Phi^{\delta,\epsilon}}{\delta m}(m_{0},\cdot) \Bigr\|_{L^\infty} 
\leq c_1. 
\end{equation*}
\end{prop}

\begin{rmk} 
\label{rmk:measureboundedbelow}
We note that here and in the rest of the paper we use the notation $m \geq c \leb$ to mean that $m(A) \geq c \leb(A)$ for all Borel sets $A \subset \T^d$. Equivalently, the density of the absolutely continuous part of $m$ is lower bounded by $c$.
\end{rmk}

\begin{proof}
{ \ } 
The proof of Proposition \ref{prop:threshold} relies on a series of steps.
\newline
\underline{Step (1).}
The very preliminary one is to mollify the function $\Phi$ by means of the mollification procedure 
introduced in Lemma 
\ref{lem:regularisation}. In turn, we can define $\Phi^{n,\delta}$ and next $\Phi^{n,\delta,\epsilon}$ by replacing 
 $\Phi$ by $\Phi^n$ in the definitions of the latter functions (with the same choice of $\rho$ as in 
 Lemma \ref{lem:regularisation}). 
By Lemmas
\ref{lem:regularisation}
and
\ref{lem:convolution:rhodelta}, the function
$\Phi^{n,\delta}$ is continuously differentiable (continuity w.r.t. $m$ being understood for $d_1$) and, for any $m \in {\mathcal P}({\mathbb T}^d)$, 
\begin{equation}
\label{eq:bound:deltaVn,delta}
\Bigl\| \frac{\delta \Phi^{n,\delta}}{\delta m}(m,\cdot) \Bigr\|_{s} \leq \frac{c_1 \Gamma(d,\rho,s)}{ \delta^{s-1}}. 
\end{equation} 

The main idea is to show that, for any $n \geq 1$ and under the lower bound for $m_0$ prescribed in the statement, the equation 
\begin{equation}
\label{eq:fixed:point}
m - m_0 = \epsilon \Bigl( \frac{\delta \Phi^{n,\delta}}{\delta m}(m,\cdot) \Bigr)^*
\end{equation}
has a solution $m_{n}$ and that this solution identifies with the maximizer $m_{n,\delta,\epsilon}$ in the definition of $\Phi^{n,\delta,\epsilon}$. 
Using \eqref{dualformula}, we have that for any 
${\boldsymbol k} \in {\mathbb Z}^d$, 
\begin{equation}
\label{eq:fourier:dual}
\reallywidehat{ \displaystyle \Bigl( \frac{\delta \Phi^{n,\delta}}{\delta m}(m,\cdot) \Bigr)^*}^{\boldsymbol k} 
=  \big(1 + \sum_{i = 1}^d |k^i|^{2s} \big)
\reallywidehat{\displaystyle  \frac{\delta \Phi^{n,\delta}}{\delta m}(m,\cdot) }^{\boldsymbol k}.
\end{equation} 
When ${\boldsymbol k}=0$, the left-hand side is equal to $0$, courtesy of 
our choice of a normalization 
for the flat derivative. In particular, the dual element belongs to
\begin{equation*}
H_0^{-s}({\mathbb T}^d):=\Bigl\{q \in H^{-s}({\mathbb T}^d) : \widehat{q}^0 = 0 \Bigr\}. 
\end{equation*}
By 
\eqref{eq:fourier:dual}
and then 
\eqref{eq:bound:deltaVn,delta}, 
we obtain, for any 
real $r\geq 0$, 
\begin{equation}
\label{eq:fourier:dual:2}
\begin{split}
\biggl\|
\biggl( 
\displaystyle  \frac{\delta \Phi^{n,\delta}}{\delta m}(m,\cdot)
\biggr)^* 
\biggr\|_{r}^2
&= \sum_{{\boldsymbol k} \in {\mathbb Z}^d \setminus \{0\}}
\big(1 + \sum_{i = 1}^d |k^i|^{2r} \big)
\biggl( 
\biggl\vert 
\reallywidehat{\Bigl( \displaystyle  \frac{\delta \Phi^{n,\delta}}{\delta m}(m,\cdot) \Bigr)^*}^{\boldsymbol k}
\biggr\vert^2 
\biggr)
\\
&= \sum_{{\boldsymbol k} \in {\mathbb Z}^d \setminus \{0\}}
\biggl( 
\big(1 + \sum_{i = 1}^d |k^i|^{2s} \big)^2 \big(1 + \sum_{i = 1}^d |k^i|^{2r} \big)
\biggl\vert 
\reallywidehat{  \displaystyle  \frac{\delta \Phi^{n,\delta}}{\delta m}(m,\cdot) }^{\boldsymbol k}
\biggr\vert^2 
\biggr)
\\
& \leq C(d,r,s)
\biggl\|
  \frac{\delta \Phi^{n,\delta}}{\delta m}(m,\cdot)
\biggr\|_{2s+r}^2 \leq \frac{C(d,r,s) c_1^2 \Gamma^2(d,\rho,2s+r) }{\delta^{2(2s+r-1)}},
\end{split} 
\end{equation}
where $C(d,r,s)$ is a constant only depending on $d$, $r$ and $s$.

When $r=d/2+\eta$, Sobolev's embedding implies 
\begin{equation}
\label{eq:fourier:dual:3}
\begin{split}
\biggl\|
\biggl( 
\displaystyle  \frac{\delta \Phi^{n,\delta}}{\delta m}(m,\cdot)
\biggr)^* 
\biggr\|_{\linf}
\leq \gamma(d,\eta) 
\biggl\|
\biggl( 
\displaystyle  \frac{\delta \Phi^{n,\delta}}{\delta m}(m,\cdot)
\biggr)^* 
\biggr\|_{d/2+\eta}
\leq    \frac{c_1 \gamma(d,\rho,s,\eta) }{\delta^{2s+d/2+\eta-1}},
\end{split} 
\end{equation}
for two (new) 
constants $\gamma(d,\eta)$ and $\gamma(d,\rho,s,\eta)$. 
\vskip 4pt

\underline{Step (2).}
We now address the solvability of the equation 
\eqref{eq:fixed:point}.
The point is to apply Schauder's fixed point theorem in the subset 
\begin{equation*} 
{\mathcal C} = \biggl\{ m \in {\mathcal P}({\mathbb T}^d) : 
\|
m-m_0 \|_{L^\infty}
\leq \frac{c_1  \epsilon \gamma(d,\rho,s,\eta) }{\delta^{2s+d/2+\eta-1}} \biggr\},
\end{equation*} 
seen as a convex subset of $1+H^{-s}_0({\mathbb T}^d)$ (where $1$ is seen as the constant 
function, equal to 1).  
In the above right-hand side, neither $m$ nor $m_0$ is required to be in $L^\infty({\mathbb T}^d)$ (i.e., 
to have a bounded density) but the difference has to be.
The $L^\infty$ norm of $m-m_0$ can be expressed as 
\begin{equation} 
\label{eq:expression:l:infty}
\| m-m_0 \|_{\linf} = \sup_{f \in H^s({\mathbb T}^d) : \| f \|_{L^1} \leq 1} \int_{{\mathbb T}^d} 
f(x) d\bigl(m-m_0 \bigr)(x).
\end{equation}

Notably, ${\mathcal C}$ is closed under $\| \cdot \|_{-s}$. 
Indeed, 
${\mathcal P}({\mathbb T}^d)$ is stable under $\| \cdot \|_{-s}$
and convergence of a sequence of probability measures under 
$\| \cdot \|_{-s}$ is equivalent to weak convergence. The latter is a direct consequence of Prokhorov's theorem: any sequence of probability measures converging in $\| \cdot \|_{-s}$ has a weakly converging subsequence and the limits should obviously coincide. Moreover,
by a lower semi-continuity argument and by \eqref{eq:expression:l:infty}, 
the condition on $\| m - m_0\|_{L^{\infty}}$ in the definition of ${\mathcal C}$ can also be shown to be stable under 
$\| \cdot \|_{-s}$ (and also under weak convergence). 

The same reasoning shows that ${\mathcal C}$ is compact under $\| \cdot \|_{-s}$. 

Next, we observe that, whenever 
\begin{equation}
\label{eq:lower:bound} 
m_0 \geq   \frac{c_1   \epsilon \gamma(d,\rho,s,\eta) }{\delta^{2s+d/2+\eta-1}} \leb ,
\end{equation}
the set ${\mathcal C}$
can be merely described as  
\begin{equation}
\label{eq:application:bochner} 
 {\mathcal C} = \biggl\{ m \in 1 + H_0^{-s}({\mathbb T}^d) : 
\| m - m_0 \|_{\linf}
\leq  \frac{c_1   \epsilon \gamma(d,\rho,s,\eta) }{\delta^{2s+d/2+\eta-1}}  \biggr\}.
\end{equation} 
It suffices to check that any $m$ as in the right-hand side
is a probability measure. 
Obviously, $m$ is a finite signed measure with total mass one.  
The point is to prove that, for any non-negative function $f \in H^s({\mathbb T}^d)$, 
$\int_{{\mathbb T}^d} f(x) dm(x) \geq 0$. 

Clearly, 
by 
\eqref{eq:lower:bound} 
and by 
definition of ${\mathcal C}$, 
\begin{equation*} 
\begin{split}
\int_{{\mathbb T}^d} f(x) dm(x)
&= \int_{{\mathbb T}^d} f(x) dm_0(x)+ 
\int_{{\mathbb T}^d} f(x) d\bigl(m-m_0\bigr)(x)
\\
&\geq \frac{c_1  \epsilon \gamma(d,\rho,s,\eta) }{\delta^{2s+d/2+\eta-1}} \|f\|_{L^1} - \| f \|_{L^1} \| m - m_0 \|_{L^{\infty}}
\\
&\geq 
\frac{c_1  \epsilon \gamma(d,\rho,s,\eta) }{\delta^{2s+d/2+\eta-1}} 
\|f \|_{L^1} - 
\frac{c_1   \epsilon \gamma(d,\rho,s,\eta) }{\delta^{2s+d/2+\eta-1}} \| f \|_{L^1} \geq 0.
\end{split}
\end{equation*}

The next step is to consider the map
\begin{equation*}
\phi : m \in {\mathcal C} \mapsto m_0 + \epsilon 
 \Bigl( \frac{\delta \Phi^{n,\delta}}{\delta m}(m,\cdot) \Bigr)^*,
\end{equation*}
under the assumption 
\eqref{eq:lower:bound}.  
We claim that ${\mathcal C}$ is stable by $\phi$.
Indeed, by 
\eqref{eq:fourier:dual:3}, we get, 
for any $m \in {\mathcal P}({\mathbb T}^d)$ (which contains ${\mathcal C}$),  
\begin{equation*} 
\begin{split}
\bigl\| \phi(m) - m_0 
\bigr\|_\infty
&\leq \epsilon 
\biggl\|
\biggl( 
\displaystyle  \frac{\delta \Phi^{n,\delta}}{\delta m}(m,\cdot)
\biggr)^* 
\biggr\|_{\linf}
\leq
 \frac{c_1\epsilon \gamma(d,\rho,s,\eta)}{\delta^{(2s+d/2+\eta-1)}},
\end{split} 
\end{equation*} 
which, by 
\eqref{eq:application:bochner}, says that 
$\phi(m)$ indeed belongs to ${\mathcal C}$. 

In order to apply Schauder's theorem, it remains to show that 
$\phi$ is continuous from ${\mathcal C}$ into itself, when the latter 
is equipped with $\| \cdot \|_{-s}$. 
By 
\eqref{eq:fourier:dual:2} (with $r=d/2+\eta$),
the Fourier coefficients of 
$([\delta \Phi^{n,\delta}/\delta m](m,\cdot))^*$ are dominated by a summable sequence and, in turn,  
 it suffices to prove that each Fourier coefficient 
\begin{equation*}
\reallywidehat{\displaystyle  \frac{\delta \Phi^{n,\delta}}{\delta m}(m,\cdot) }^{\boldsymbol k}
\end{equation*}
is continuous with respect to $m$ for $\| \cdot \|_{-s}$. We already know that it is continuous 
with respect to $m$ for $d_1$ (since $\Phi^{n,\delta}$ is continuously differentiable). We then recall that convergence in 
$\| \cdot \|_{-s}$ implies weak 
convergence and, in turn, convergence in $d_1$ (by compactness of ${\mathbb T}^d$).
\vskip 4pt

\underline{Step (3).} By the second step, we know that 
\eqref{eq:fixed:point} has a fixed point $m_n$. We claim that this fixed point is in fact a strict maximizer of 
the function 
\begin{equation*}
\Psi : m \in {\mathcal P}({\mathbb T}^d) \mapsto \Phi^{n,\delta}(m) - \frac1{2 \epsilon} 
\| m - m_0 \|_{-s}^2.  
\end{equation*}
Indeed, 
we observe that, for any two $m,m' \in {\mathcal P}({\mathbb T}^d)$,   
\begin{equation*}
\begin{split}
\frac{d}{d\lambda} \Psi \Bigl( \lambda m' + (1-\lambda) m \Bigr){}_{\vert \lambda =0} 
&= 
\int_{{\mathbb T}^d} 
\frac{\delta \Phi^{n,\delta}}{\delta m}(m,x) d \bigl( m'-m \bigr)(x) 
- 
\frac1{\epsilon} \langle m'-m, m-m_0 \rangle_{-s}. 
\end{split}
\end{equation*}
By 
the definition of the dual element, 
\begin{equation*}
\begin{split}
\frac{d}{d\lambda} \Psi \Bigl( \lambda m' + (1-\lambda) m \Bigr){}_{\vert \lambda =0} 
&= 
\Bigl\langle m'-m, \Bigl( \frac{\delta \Phi^{n,\delta}}{\delta m}(m,\cdot) \Bigr)^*
\Bigr\rangle_{-s}
-
\frac1{\epsilon} \langle m'-m, m-m_0 \rangle_{-s}. 
\end{split}
\end{equation*}
In particular, choosing $m$ as $m_n$ (and then replacing $m'$ by $m$), we get 
\begin{equation}
\label{eq:critical:point}
\frac{d}{d\lambda} \Psi \Bigl( \lambda m + (1-\lambda) m_n \Bigr){}_{\vert \lambda =0} 
=0, \quad m \in {\mathcal P}({\mathbb T}^d), 
\end{equation}
which shows that 
$m_n$ is critical point of $\Psi$. 

We then argue by concavity to prove that $m_n$ is the unique maximizer of $\Psi$. 
Indeed, 
using the semi-concavity property of $\Phi^{n,\delta}$ (see Lemma
\ref{lem:convolution:rhodelta}),
we get, for any $m,m' \in {\mathcal P}({\mathbb T}^d)$,  
\begin{equation*}
\begin{split}
\Psi \Bigl( \lambda m' + (1-\lambda) m \Bigr)
&= \Phi^{n,\delta} \Bigl( \lambda m' + (1-\lambda) m \Bigr) - \frac1{2\epsilon} \bigl\| \lambda 
(m'-m_0) + (1-\lambda) (m-m_0) \bigr\|_{-s}^2
\\
&\geq \lambda \Phi^{n,\delta}(m') + (1-\lambda) \Phi^{n,\delta}(m) - \frac{\lambda}{2\epsilon} 
\| m' \|_{-s}^2 - \frac{1-\lambda}{2\epsilon} 
\| m \|_{-s}^2 
\\
&\hspace{15pt} - \frac{c_2 \Gamma(d,\rho,s)}{2 \delta^{2(s-1)}} \lambda (1-\lambda) \|m-m' \|_{-s}^2 
+ \frac{1}{2 \epsilon}   \lambda (1-\lambda) \|m-m' \|_{-s}^2
\\
&\geq \lambda \Psi(m') + (1-\lambda) \Psi( m ) + c   \lambda (1-\lambda) \|m-m' \|_{-s}^2,
\end{split}
\end{equation*} 
with $c:=1/\epsilon-c_2 \Gamma(d,\rho,s)/ {\delta^{2(s-1)}} >0$. 
And then,
choosing $m$ as $m_n$ (and then replacing $m'$ by $m$), we obtain, for any $m \in {\mathcal P}({\mathbb T}^d)$, 
\begin{equation*}
\begin{split}
\Psi(m) - \Psi(m_n) \leq 
\frac1{\lambda} 
\Bigl[ \Psi \Bigl( \lambda m + (1-\lambda) m_n \Bigr)
- \Psi(m_n) \Bigr] 
- c(1-\lambda) 
\|m-m_n \|_{-s}^2.
\end{split}
\end{equation*}
Letting $\lambda$ tend to $0$ and invoking \eqref{eq:critical:point}, we deduce that 
$m_n$ is indeed the unique maximizer of $\Psi$. Below, we thus write 
$m_{n,\delta,\epsilon}$ for $m_n$. 
The conclusion is that the maximizer of $\Psi$ satisfies 
\begin{equation}
\label{eq:mndeltaepsilon}
\bigl\| m_{n,\delta,\epsilon} - m_0 \bigr\|_{\linf}
\leq  \frac{c_1   \epsilon \gamma(d,\rho,s,\eta) }{\delta^{2s+d/2+\eta-1}},
\end{equation}
under the lower bound
\eqref{eq:lower:bound}. 
Also, from the fixed point identity 
\begin{equation*} 
\Bigl( 
\frac{\delta \Phi^{n,\delta}}{\delta m}(m_{n,\delta,\epsilon},\cdot) 
\Bigr)^* = \frac{m_{n,\delta,\epsilon}-m_0}{\epsilon},
\end{equation*} 
we deduce that 
\begin{equation}
\label{eq:dual:mndeltaepsilon} 
\Bigl( \frac{m_{n,\delta,\epsilon}-m_0}{\epsilon}
\Bigr)^{*} = \frac{\delta \Phi^{n,\delta}}{\delta m}(m_{n,\delta,\epsilon},\cdot). 
\end{equation} 

In particular, $[(m_{n,\delta,\epsilon}-m_0)/\epsilon]^{*}$
is a smooth function and 
\eqref{eq:dual:mndeltaepsilon} says that, as such, it is $c_1$-Lipschitz continuous (on ${\mathbb T}^d$) (the $x$-Lipschitz property of 
 $\frac{\delta \Phi^{n,\delta}}{\delta m}(m_{n,\delta,\epsilon},\cdot)$
 follows from the $\mu$-Lipschitz property of 
 $\Phi^{n,\delta}$, see for instance 
 \cite[Propositions 5.33 and 5.36]{CarmonaDelarue_book_I}). 
Importantly, $[(m_{n,\delta,\epsilon}-m_0)/\epsilon]^{*}$ has a zero mean (which can be shown by computing the zero Fourier mode) and is thus bounded by a constant that only depends on $c_1$ and $d$. 
\vskip 4pt

\underline{Step (4).} 
We now let $n$ tend to $\infty$. 
We recall from Lemmas 
\ref{lem:regularisation}
and 
\ref{lem:convolution:rhodelta} that $(\Phi^{n,\delta})_{n \geq 1}$ converges uniformly to 
$\Phi^{\delta}$. 

Obviously, the sequence $(m_{n,\delta,\epsilon})_{n \geq 1}$
is compact (for the weak convergence and also for $d_1$ and $\| \cdot \|_{-s}$) and passing to the limit in the inequality
\begin{equation*}
\Phi^{n,\delta}(m_{n,\delta,\epsilon}) - \frac1{2 \epsilon} 
\bigl\| m_{n,\delta,\epsilon} - m_0 \bigr\|_{-s}^2
\geq 
\Phi^{n,\delta}(m) - \frac1{2 \epsilon} 
\bigl\| m - m_0 \bigr\|_{-s}^2, \quad m \in {\mathcal P}({\mathbb T}^d), 
\end{equation*} 
we deduce that any weak limit point is a minimizer of the right-hand side. Recalling that $\delta^{2s-2}>2 c_2 \Gamma(d,\rho,s)\epsilon$, using the semi-concavity of $\Phi^{\delta}$ respect to $\|.\|_{-s}$ given by Lemma \eqref{lem:convolution:rhodelta} as well as Claim \textit{(3)} in Proposition \eqref{prop:supconvproperties} we deduce that the only possible limit is $m_{\delta,\epsilon}$ the unique maximum in the definition of $\Phi^{\delta,\epsilon}(m_0)$. In particular the whole sequence $(m_{n,\delta,\epsilon})_{ n \geq 0}$ converges to $m_{\delta,\epsilon}$.

Passing to the limit in 
\eqref{eq:mndeltaepsilon} (which is possible by the same argument as the one used in the first step of the proof), we get that 
\begin{equation}
\label{eq:mdeltaepsilon}
\bigl\| m_{\delta,\epsilon} - m_0 \bigr\|_{\linf}
\leq  \frac{c_1   \epsilon \gamma(d,\rho,s,\eta) }{\delta^{2s+d/2+\eta-1}},
\end{equation}
under the lower bound
\eqref{eq:lower:bound}. 

Lastly, we have that, for any $q \in H^{-s}({\mathbb T}^d)$, 
\begin{equation}
\Bigl\langle
\Bigl( \frac{m_{n,\delta,\epsilon}-m_0}{\epsilon}
\Bigr)^{*}, q \Bigr\rangle 
\underset{n \rightarrow \infty}{\longrightarrow} 
\Bigl\langle
\Bigl( \frac{m_{\delta,\epsilon}-m_0}{\epsilon}
\Bigr)^{*}, q \Bigr\rangle,  
\end{equation} 
as $n$ tends to $\infty$. As the functions 
$([(m_{n,\delta,\epsilon}-m_0)/\epsilon]^{*})_{n \geq 1}$ in the left-hand side
are uniformly bounded and $c_1$-Lipschitz-continuous, they do converge in sup norm, which shows that 
$[(m_{\delta,\epsilon}-m_0)/\epsilon]^{*}$ is $c_1$-Lipschitz continuous. It remains to recall from the analysis of the sup-convolution (see 
item \textit{(3)} in
Proposition \ref{prop:supconvproperties}) that 
\begin{equation*} 
\Bigl( 
\frac{\delta \Phi^{\delta,\epsilon}}{\delta m}(m_0,\cdot) 
\Bigr)^* = \frac{m_{\delta,\epsilon}-m_0}{\epsilon},
\end{equation*} 
which shows that 
\begin{equation*} 
\Bigl( 
\frac{m_{\delta,\epsilon}-m_0}{\epsilon}
\Bigr)^{*} = 
\frac{\delta \Phi^{\delta,\epsilon}}{\delta m}(m_0,\cdot).
\end{equation*} 
And then the right-hand side is $c_1$-Lipschitz continuous, which proves the last expected result. 
\end{proof}

\section{Proofs of the “hard inequalities"}

\label{sec:hardinequalities}

In this section we prove the “hard inequalities", i.e. the second inequality in Theorem \ref{thm:maind1} and the second inequality in Theorem \ref{thm:mainhs}. Recall that our main difficulty is that $U$ is not smooth, so simply plugging the projection $U^N(t,\bx) = U(t,m_{\bx}^N)$ into the PDE \eqref{hjbn} gives us no information. Our strategy is therefore to attempt to regularize $U$ while keeping track of the subsolution property. That is, we want to approximate $U$ by functions which are smooth enough and which are \textit{almost} subsolutions of \eqref{hjbinf}. In the $H^{-s}$-regular case, this is accomplished in just one step, since sup-convolution in $H^{-s}$ (as presented in Subsection \ref{subsec:supconv}) creates enough regularity for our purposes. In the $d_1$-regular case, the analysis is much more involved, since in this step we must start by “upgrading" $d_1$-regularity to $H^{-s}$-regularity by means of the regularization by mollification presented in Subsection \ref{subsec:mollificationkernel}, and only then apply a sup-convolution in $H^{-s}$. In each step the regularity of the approximation must be tracked, and the amount by which it fails to be a subsolution must be estimated.

\subsection{Analysis of $U^{\epsilon}$}

\label{subsec:uepsanalysis}

This subsection is concerned with the $H^{-s}$-regular case, and so throughout this subsection Assumption \ref{assump:Hs} holds, and in particular an integer $s > d/2 + 2$ is fixed. For $\epsilon >0$, we follow \eqref{eq:sup:convolution:4.1} in Proposition 
\ref{prop:supconvproperties} and define the sup-convolution $U^{\epsilon}: [0,T] \times H^{-s} \rightarrow \R$ by
\begin{equation} 
U^{\epsilon}(t,q) := \sup_{m \in \mathcal{P}(\T^d)} \Big\{ U(t,m) - \frac{1}{2\epsilon} \|q-m \|_{-s}^2 \Big\}.
\label{Definitionsupconvolution10/02/2023}
\end{equation}
 The goal of this subsection is to understand three questions: 
\begin{itemize}
    \item How close is $U^{\epsilon}$ to $U$? 
    \item What are the regularity properties of $U^{\epsilon}$? 
    \item By how much does $U^{\epsilon}$ fail to be a sub-solution of the PDE \eqref{hjbinf}?
\end{itemize}

Some properties of $U^{\epsilon}$ can be already be inferred directly from Propositions \ref{prop:upropshs} and \ref{prop:supconvproperties}, and the other relevant properties are summarized in the following proposition.

\begin{prop} \label{lem:uepsproperties} Let Assumption \ref{assump:Hs} hold and let $U^{\epsilon}$ be the sup-convolution defined in \eqref{Definitionsupconvolution10/02/2023}. Then, there is a constant $C > 0$ such that the following holds for all $\epsilon$ small enough:

\begin{enumerate}
    \item For all $t_1,t_2 \in [0,T]$ and $m_1,m_2 \in \pr(\T^d)$, we have 
    \begin{align*}
        |U^{\epsilon}(t_1,m_1) - U^{\epsilon}(t_2,m_2)| \leq C \big( |t_1 - t_2| + \|m_1 - m_2\|_{-s} \big).
    \end{align*}
    \item  The map $[0,T] \times H^{-s} \ni (t,q) \mapsto \nabla_{-s}U^{\epsilon}(t,q) \in H^{-s}$ is jointly continuous.
    \item $U^{\epsilon}$ satisfies
    \begin{align}
        - \partial_t U^{\epsilon}(t,m) - \int_{\T^d} \Delta_x \frac{\delta U^{\epsilon}}{\delta m}(t,m,x) m(dx) + \int_{\T^d} H\big(x, D_x \frac{\delta U^{\epsilon}}{\delta m}(t,m,x) \big)m(dx)  \leq \sF(m) + C\epsilon, 
        \label{InequationforUepsilon}
    \end{align}
    in the viscosity sense with test functions $\phi \in \cC^1((0,T) \times H^{-s})$.
\end{enumerate}
\end{prop}

\begin{rmk} \label{rmk:testfunctions}
    When we say that $U^{\epsilon}$ satisfies \eqref{InequationforUepsilon} in the viscosity sense with test functions in $\cC^{1}((0,T) \times H^{-s})$, we mean that for any $\phi \in \cC^1((0,T) \times H^{-s})$ and any $(t_0,m_0) \in (0,T) \times \mathcal{P}(\T^d)$  
    such that 
    \begin{align*}
        U^{\epsilon}(t_0,m_0) - \phi(t_0,m_0) = \sup_{(t,q) \in (0,T) \times H^{-s}} \Big\{ U^{\epsilon}(t,q) - \phi(t,q) \Big\},
    \end{align*}
    we have
    \begin{align}
        - \partial_t \phi(t_0,m_0) - \int_{\T^d} \Delta_x \frac{\delta \phi}{\delta m}(t_0,m_0,x) m_0(dy) + \int_{\T^d} H\big(x, D_x \frac{\delta \phi}{\delta m}(t_0,m_0,x) \big)m_0(dx)  \leq \sF(m_0) + C\epsilon.
    \end{align}
\end{rmk}

\begin{rmk}
Combining the bound on $\big[ U^{\epsilon} (t,\cdot) \big]_{\cC^{1,1}(H^{-s})}$ obtained from Proposition \ref{prop:supconvproperties} with Proposition \ref{prop:dmlipschitz}, we deduce that $m \mapsto D_m U^{\epsilon}(t,\cdot,x)$ is Lipschitz continuous over $\mathcal{P}(\T^d)$, uniformly in 
$x \in {\mathbb T}^d$ (but the Lipschitz constant depends on $\epsilon$).
\label{RemarkLipDm}
\end{rmk}

The rest of this subsection is devoted to proving Proposition \ref{lem:uepsproperties}. We start with a preliminary observation, which is that the inner product $\langle \cdot, \cdot \rangle_s$ (see Subsection \ref{subsec:notation})
can be re-written as
\begin{align} \label{hsinnerprod}
\langle f, g \rangle_s = \langle f, g \rangle_{L^2} + (2 \pi)^{-2s} \sum_{i = 1}^d \langle D_{x^i}^s f, D_{x^i}^s g \rangle_{L^2}, 
\end{align}
with $\langle f, g \rangle_{L^2} = \int_{\T^d} fg \, dx$ the usual $L^2$ inner product. Indeed, recalling that $s$ is an integer, \eqref{hsinnerprod} can be seen directly from the fact that for $\bm k \in \Z^d$, $\widehat{D_{x^i}^s f}^{\bm k} = (2 \pi \i k_i)^s \widehat{f}^{ \bm k}$. We now give two key lemmas.

\begin{lem} \label{lem:faadibruno}
    Let Assumption \ref{assump:Hs} hold, and let $f \in H^{s+1}$. Then there is a constant  $C = C(\|f\|_s)$ depending on $f$ only through the norm $\|f\|_s$ such that
    $$ | \langle H(\cdot,Df),f \rangle_s| \leq C.$$
\end{lem}

\begin{proof}
By \eqref{hsinnerprod}, we have 
\begin{align*}
    | \langle H(\cdot,Df),f \rangle_s| \leq |\langle H(\cdot, Df), f \rangle_{L^2}| + (2\pi)^{-2s} \sum_{i = 1}^d |\langle D_{x^i}^s \big[H(\cdot, Df)\big], D_{x^i}^s f \rangle_{L^2}|.
\end{align*}
Since $s > d/2 + 2$, $\|Df\|_{\linf} \leq C\|f\|_s$, and so clearly 
\begin{align*}
    |\langle H(\cdot, Df), f \rangle_{L^2}| \leq \|H(\cdot, Df)\|_{\linf}^{1/2} \|f\|_{L^2}^{1/2} \leq C(\|f\|_s){,}
\end{align*}
where $C(\|f\|_s)$ is a non-decreasing function of 
$\|f\|_s$.
To complete the proof, it suffices to fix $i \in \{1,...,d\}$ and show that
\begin{align*}
    \Big|\int_{\T^d} D_{x^i}^s [H(\cdot, Df)] D_{x^i}^s f dx \Big| \leq C, \quad C = C(\norm{f}_s).
\end{align*}
In order to do so, we can use (a generalization of) the Fa\`a Di Bruno formula to write $D_{x^i}^s \left[ H(\cdot,Df) \right]$ as a sum of terms of the form
\begin{equation}
 D_{x^i}^{ s-k} D_{p^{i_1}...p^{i_q}} H(\cdot, Df) D_{x^i}^{j_1} D_{x^{i_1}} f ... D_{x^i}^{j_q} D_{x^{i_q}} f,
\label{somanyderivatives}
\end{equation}
with $k,q \in \{0,...,s\}$, $i_l \in \{1,...,d\}$, $j_l \in \N$ satisfying $   \sum_{l = 1}^q j_l = k.$
So, to prove the Lemma, it in fact suffices to prove the following claim: 
\newline \newline 
\textbf{Claim:} For any $k,q$, $(i_l)_{l = 1,...,q}$ and $(j_l)_{l =1,...,q}$ as above, we have
\begin{align} \label{sufficientest}
    \Big| \int_{\T^d} \Big( D_{x^i}^{s-k} D_{p^{i_1}...p^{i_q}} H(\cdot, Df) D_{x^i}^{j_1} D_{x^{i_1}} f ... D_{x^i}^{j_q} D_{x^{i_q}} f \Big) D_{x^i}^s f dx \Big| \leq C, \quad C = C(\norm{f}_s).
\end{align}
We now prove the claim by considering several cases. 
\newline \newline 
\underline{Case 1} ($q = 0$): In the case $q = 0$ (and thus $k=0$), the left-hand side of \eqref{sufficientest} becomes 
\begin{align*}
    \Big| \int_{\T^d} D_{x^i}^{s} H(\cdot, Df)  D_{x^i}^{s} f dx \Big|, 
\end{align*}
and the estimate is easily proved using Cauchy-Schwarz and then Sobolev embedding (here we use $s > d/2 + 1$), and the fact that $H \in \cC^s$. 
Indeed, by Sobolev embedding, the argument $Df$ inside $D_{x^i}^s H$ is bounded by a constant $C(\| f \|_s)$ (depending on $\|f\|_s$), and by continuity, 
$D_{x^i}^s H$ is bounded on the compact subset $\{(x,p) \in {\mathbb T}^d \times 
{\mathbb R}^d : \vert p \vert \leq C(\| f \|_s)\}$.
\newline \newline
\underline{Case 2} ($q = 1$ and $k = s$): If $q = 1$ and $k = s$, then the left-hand side of \eqref{sufficientest} becomes 
\begin{align*}
   \biggl|\int_{\T^d} D_{p^{i_1}}H({\cdot}, Df) D_{x^i}^s D_{x^{i_1}} f D_{x^i}^s f dx\biggr|, 
\end{align*}
so we can use integration by parts to estimate
\begin{align*}
    \biggl|\int_{\T^d}& D_{p^{i_1}}H({\cdot}, Df) D_{x^i}^s D_{x^{i_1}} f D_{x^i}^s f dx    \biggl|= \frac{1}{2}  \biggr| \int_{\T^d} D_{p^{i_1}}H(\cdot, Df) D_{x^{i_1}} |D_{x_i}^s f|^2 dx    \biggr| \\
&=     \biggl| \frac{1}{2} \int_{\T^d} D_{x^{i_1}}\big[D_{p^{i_1}} H(\cdot, Df) \big] |D_{x^i}^s f|^2  dx   \biggr|
\leq \frac{1}{2} \norm{D_{x^{i_1}}\big[D_{p^{i_1}} H(\cdot, Df) \big]}_{\linf} \|f\|_{s}^2.
\end{align*}
To conclude \eqref{sufficientest} (and thus complete the proof), we need only 
to proceed as in Case 1 and
notice that because of the smoothness of $H$ and the fact that $s > d/2 + 2$, we have $\|D_{x^{i_1}}\big[D_{p^{i_1}} H(\cdot, Df) \big]\|_{\linf} \leq C(\|f\|_{s})$.
\newline \newline
\underline{Case 3} ($q = 1$ and $k < s$): In this case, the left-hand side of
\eqref{sufficientest}
 becomes 
\begin{align*}
  \biggl| \int_{\T^d} D_{x^i}^{s-k} D_{p^{i_1}} H(\cdot, Df) D_{x^i}^k D_{x^{i_1}} f D_{x^i}^s f dx \biggr|.
\end{align*}
Using Sobolev embedding and the smoothness of $H$ as in the Case 2, we see that 
the following estimate suffices: 
\begin{align*}
    \int_{\T^d} &| D_{x^i}^{k} D_{x^{i_1}} f| |  D^{s} f| dx \leq \| D_{x^i}^{k} D_{x^{i_1}} f \|_{L^2} \| D^{s} f \|_{L^2} \leq C\|f\|_s^2,
\end{align*}
where we use Cauchy-Schwarz and the fact that $k < s$. 
\newline \newline 
\underline{Case 4} ($q > 1$) In this case, we start
again by noticing that by Sobolev embedding,
\begin{align}
    \|D_{x^i}^{s-k} D_{p^{i_1}...p^{i_q}} H(\cdot, Df) \|_{\linf} \leq C(\|f\|_s).
\end{align}
Moreover, if for some $l \in \{1,...,q\}$, we have \begin{align} \label{sobolevthreshold}
\frac{1}{2} + \frac{j_l + 1}{d} < \frac{s}{d}, 
\end{align}
then by Sobolev embedding again, we have $\|D_{x^i}^{j_l} D_{x^{i_l}} f\|_{\linf} \leq C\|f\|_s$. Without loss of generality, we can arrange the terms in \eqref{sufficientest} so that for some $r \in \{0,1,...,q\}$ we have that $j_l$ satisfies \eqref{sobolevthreshold} for $l > r$, and $j_l$ fails to satisfy \eqref{sobolevthreshold} for $l \leq r$. So, ignoring the terms in \eqref{sufficientest} which are bounded in $\linf$ by constants depending on $\|f\|_s$, we find that in order to complete the proof in Case 4, it suffices to show that
\begin{align} \label{claim2}
    \int_{\T^d} \bigl| D_{x^i}^{j_1} D_{x^{i_1}} f ... D_{x^i}^{j_r} D_{x^{i_r}} f  D_{x^i}^s f \bigr| dx  \leq C, \quad C = C(\norm{f}_s),
\end{align}
where $q$,$r$, $(j_l)_{l = 1,...,q}$, $(i_l)_{l = 1,...,q}$ are as above and, importantly, 
\begin{align} \label{threshold2}
    \frac{1}{2} - \frac{s}{d} \geq - \frac{j_l + 1}{d}, \quad l = 1,...,r.
\end{align}
If $r = 0$, \eqref{claim2} obviously holds, while if $r = 1$, the same argument as in Case 3 applies (since $q > 1$ implies that $j_r < s$). Thus it suffices to prove \eqref{claim2} when $r > 1$. To this end, suppose that we can find a collection of positive reals $(p_l)_{l = 1,...,r}$ such that 
\begin{enumerate}
    \item $\sum_{l = 1}^r \frac{1}{p_l} = \frac{1}{2}$
    \item $\frac{1}{2} - \frac{s}{d} < \frac{1}{p_l} - \frac{j_l + 1}{d}, \quad $ for all $l = 1,...,r$. 
\end{enumerate}
Then we have 
\begin{align*}
   \int_{\T^d} |D_{x^i}^{j_1} D_{x^{i_1}} f| ... |D_{x^i}^{j_r} D_{x^{i_r}} f| |D_{x^i}^s f| dx &\leq  \| D_{x^i}^{j_1} D_{x^{i_1}} f \|_{L^{p_1}} ... \| D_{x^i}^{j_r} D_{x^{i_r}}f \|_{L^{p_r}} \|D_{x^i}^s f \|_{L^2} 
    \\ &\leq 
    \|f\|_{j_1 + 1, p_1}...\|f\|_{j_r + 1, p_r} \|f\|_{s} 
    \leq C\norm{f}_{s}^{q+1}, 
\end{align*}
where $\| \cdot \|_{j, p}$ denotes the usual Sobolev norm on the space $W^{j,p}$ of functions with $j$ weak derivatives each lying in $L^p$, and where the first inequality comes from the generalized H\"older's inequality and the last one comes from the fact that if $\frac{1}{2} - \frac{s}{d} < \frac{1}{p_j} - \frac{j_l + 1}{d}$, then $
    H^s \hookrightarrow W^{j_l +1, p_l}$
by Sobolev embedding. So, clearly to prove Claim 2, it suffices to prove that when $r > 1$, we can find $(p_l)_{l = 1,...,r}$ satisfying the two properties listed above. To prove this, it in turn suffices to establish that when $r > 1$, we have 
\begin{align} \label{pjprop}
    \sum_{l = 1}^r \Big(\frac{1}{2} - \frac{s}{d} + \frac{ j_l + 1}{d} \Big) < \frac{1}{2}, 
\end{align}
since then the $(p_l)_{l=1,...,r}$ defined by the equation
\begin{align} \label{pldef}
    \frac{1}{p_l} = \frac{1}{2} - \frac{s}{d} + \frac{ j_l + 1}{d} + \frac{1}{r} \bigg( \frac{1}{2} - \sum_{m = 1}^r \Big(\frac{1}{2} - \frac{s}{d} + \frac{ j_m + 1}{d} \Big) \bigg), \quad l=1,...,r,
\end{align}
will do the job. Indeed, the condition \eqref{threshold2} guarantees that the numbers $p_l$ defined in \eqref{pldef} are positive reals, and the condition (1) and (2) listed above are satisfied by construction. 
To prove \eqref{pjprop}, we use the fact that $1 < r \leq q$, $\sum_{l = 1}^r j_l \leq k \leq s$, and $s > \frac{d}{2} + 2$ to estimate 
\begin{align*}
    \sum_{l = 1}^r &\Big(\frac{1}{2} - \frac{s}{d} + \frac{ j_l + 1}{d} \Big)
    = r \big( \frac{1}{2} + \frac{1}{d} - \frac{s}{d} \big) + \frac{1}{d} \sum_{l = 1}^r j_l \\
    &\leq q \big( \frac{1}{2} + \frac{1}{d} - \frac{s}{d} \big) + \frac{s}{d} = (q-1) \big( \frac{1}{2} + \frac{1}{d} - \frac{s}{d} \big) + \frac{1}{2} + \frac{1}{d} \\
    &< (q-1) \big(\frac{1}{2} + \frac{1}{d} - \frac{2}{d} - \frac{1}{2} \big) + \frac{1}{2} + \frac{1}{d} \\
    &\leq (q-1)(- \frac{1}{d}) + \frac{1}{2} + \frac{1}{d} \leq \frac{1}{2}. 
\end{align*}
This completes the proof of the estimate \eqref{claim2} and hence of Case 4. \end{proof}

We next state a lemma which follows easily from the non-positivity of the Laplacian on $L^2$, together with the representation \eqref{hsinnerprod}.

\begin{lem} \label{lem:laplacenegative}
The Laplacian is non-positive on $H^s$. More precisely, if $f \in H^{s+2}$, then 
\begin{align*}
    \langle \Delta f, f \rangle_s \leq 0. 
\end{align*}
\end{lem}

We now use Lemmas \ref{lem:faadibruno} and \ref{lem:laplacenegative} to prove Lemma \ref{lem:uepsproperties}.

\begin{proof}[Proof of Proposition \ref{lem:uepsproperties}]
For Claim \textit{(1)}, the regularity of $U^{\epsilon}$ in $m$ follows directly from Proposition \ref{prop:supconvproperties}, and regularity time is inherited directly from the regularity of $U$ in $t$, since for any $t_1,t_2 \in [0,T]$ and $q \in H^{-s}$,
\begin{align*}
U^{\epsilon}(t_2,q) - U^{\epsilon}(t_1,q) \leq \sup_{m \in \mathcal{P}(\T^d)} 
\Bigl\{
U(t_1,m) - U(t_2,m)
\Bigr\}
 \leq C|t_1-t_2|,
\end{align*}
where we used the uniform in $m \in \mathcal{P}(\T^d)$ Lipschitz regularity of $U(\cdot,m)$.

For Claim \textit{(2)}, we first note by Proposition by Claim \textit{(3)} in Proposition \ref{prop:supconvproperties}, for all $\epsilon$ small enough and all $(t,q) \in [0,T] \times H^{-s}$ it holds 
{$$ \nabla_{-s}U^{\epsilon}(t,q) = \frac{1}{\epsilon}\bigl(m^{\epsilon}(t,q) - q \bigr).$$
Claim \textit{(2) }of the {lemma} will be proved if we can prove the continuity of the map $[0,T] \times H^{-s} \ni (t,q) \mapsto m^{\epsilon}(t,q) \in H^{-s}$. But for each fixed $t$, $q \mapsto  \nabla_{-s}U^{\epsilon}(t,q)$ is Lipschitz by Claim \textit{2} in Proposition \ref{prop:supconvproperties}, so in fact it suffices to show that for each fixed $q$, $t \mapsto m^{\epsilon}(t,q)$ is continuous. This can be established using the same compactness argument as in the fourth step in the proof of Proposition \ref{prop:threshold}, so we omit the details.}

We now turn to Claim \textit{(3)}, i.e., the subsolution estimate.
Let us fix a smooth function   $\phi = \phi(t,q) : (0,T) \times H^{-s} \to \R$ such that $U^{\epsilon} - \phi$ attains a maximum at $(t_0,m_0) \in (0,T) \times \pr(\T^d)$. Let $m_{\epsilon} \in \pr(\T^d)$ be the unique point such that 
\begin{align*}
    U^{\epsilon}(t_0,m_0) = U(t_0,m_{\epsilon}) - \frac{1}{2 \epsilon} \norm{m_0 - m_{\epsilon}}_{-s}^2, 
\end{align*}
and define $\psi(t,m) = \phi(t,m_0 - m_{\epsilon} + m)$. 
By a standard argument from the theory of viscosity solutions {which applies equally well in this infinite-dimensional setting} (see e.g. the ``alternative Proof of Theorem 5.8" in the notes \cite{Calder}), we have the following two facts :
\begin{enumerate}
    \item $U - \psi|_{(0,T) \times \pr(\T^d)}$ attains a maximum at $(t_0,m_{\epsilon})$.
    \item We have 
    \begin{align} \label{gradcharacterization}
        \nabla_{-s} \psi(t_0,m_{\epsilon}) = \nabla_{-s} \phi(t_0,m_0) = \frac{1}{\epsilon} (m_{\epsilon} - m_0),
    \end{align} and in particular 
    \begin{align*}
        \frac{ \delta \psi}{\delta m}(t_0,m_{\epsilon}) = \frac{\delta \phi}{\delta m}(t_0,m_0) = \frac{1}{\epsilon}  \big(m_{\epsilon} - m_0\big)^{\ast}.
    \end{align*}
    \end{enumerate}
Since $U - \psi$ has a maximum at $(t_0,m_{\epsilon})$ and $U$ is a subsolution of \eqref{hjbinf} (by Lemma \ref{lem:usubsol}), we deduce from 
  point (1) the inequality, 
    \begin{align*}
         &- \partial_t \phi(t_0,m_{0}) - \int_{\T^d} \Delta_x \frac{\delta \phi}{\delta m}(t_0,m_{0},x) m_{\epsilon}(dx) + \int_{\T^d} H\Bigl(x,D_x\frac{\delta \phi}{\delta m}(t_0,m_{0},x) \Bigr)m_{\epsilon}(dx) 
        \\  
        &= - \partial_t \psi(t_0,m_{\epsilon}) - \int_{\T^d} \Delta_x \frac{\delta \psi}{\delta m}(t_0,m_{\epsilon},x) m_{\epsilon}(dx) + \int_{\T^d} H\Bigl(y,D_x\frac{\delta \psi}{\delta m}(t_0,m_{\epsilon},x)\Bigr) m_{\epsilon}(dx)  \leq \sF(m_{\epsilon}). 
    \end{align*}
    Consequently, we have 
    \begin{align} \label{errors}
         &- \partial_t \phi(t_0,m_{0}) - \int_{\T^d} \Delta_x \frac{\delta \phi}{\delta m}(t_0,m_{0},x) m_{0}(dx) + \int_{\T^d} H\Bigl(y,D_x\frac{\delta \phi}{\delta m}(t_0,m_{0},x)\Bigr) m_{0}(dx)  \nonumber \\
         & \leq \sF(m_0) + E_1 + E_2 + E_3, 
    \end{align}
    where
    \begin{align*}
        E_1 &:= \int_{\T^d} \Delta_{x} \frac{\delta \phi}{\delta m}(t_0,m_0,x) (m_{\epsilon} - m_0)(dx), \\
        E_2 &:= \int_{\T^d} H\Bigl(x,D_x  \frac{\delta \phi}{\delta m}(t_0,m_0,x)\Bigr)(m_{0} - m_\epsilon)(dx), \\
        E_3 &:= \sF(m_{\epsilon}) - \sF(m_0). 
    \end{align*}
    Using point (2) above, we see that
    \begin{align} \label{e1bound}
        E_1 &= \epsilon \Bigl\langle \Delta_x \frac{\delta \phi}{\delta m}(t_0,m_0,\cdot),  \Big( \frac{\delta \phi}{\delta m}(t_0,m_0, \cdot) \Big)^{\ast} \Big\rangle_{s,-s} \nonumber \\
        &= \epsilon \Bigl\langle \Delta_x \frac{\delta \phi}{\delta m}(t_0,m_0,\cdot),  \frac{\delta \phi}{\delta m}(t_0,m_0, \cdot) \Bigr \rangle_{s} \leq 0,
    \end{align}
    by Lemma \ref{lem:laplacenegative} (which applies here since $m_{\epsilon} - m \in H^{-s+2}$, hence $\frac{\delta \phi}{\delta m}(t_0,m_0) = \epsilon (m_{\epsilon} - m)^* \in (H^{-s+2})^* = H^{s+2}$, as can be easily checked using the formula \eqref{dualformula}). For $E_2$, we use the same trick to write 
    \begin{align} \label{e2bound}
        E_2 &= \epsilon \Bigl\langle H \Bigl(\cdot, D_x \frac{\delta \phi}{\delta m}(t_0,m_0,\cdot) \Bigr), \Big(\frac{\delta \phi}{\delta m}(t_0,m_0, \cdot) \Big)^{\ast}
        \Bigr\rangle_{s,-s} \nonumber \\
        &= \epsilon \Bigl\langle H
        \Bigl(\cdot, D_x \frac{\delta \phi}{\delta m}(t_0,m_0,\cdot)\Bigr), \frac{\delta \phi}{\delta m}(t_0,m_0, \cdot)\Bigr\rangle_{s} \leq C\epsilon, 
    \end{align}
    where the last inequality follows from Lemma \ref{lem:faadibruno} together with the fact that Claim \textit{(3)} in Proposition \ref{prop:supconvproperties} shows that there is a constant $C$ independent of $\epsilon$ and $\phi$ such that
    \begin{align*}
        \norm{\frac{\delta \phi}{\delta m}(t_0,m_0,\cdot)}_s = \frac{1}{\epsilon} \norm{m_{\epsilon} - m_{0}}_{-s} \leq C.
    \end{align*}
    Finally for the third error term we simply use Lipschitz continuity of $\sF$ (see Remark \ref{rmk:hslip}) to estimate 
    \begin{align} \label{e3bound}
        E_3 \leq C \norm{m_{\epsilon} - m_{0}}_{-s} \leq C \epsilon. 
    \end{align}
    Combining \eqref{errors}, \eqref{e1bound}, \eqref{e2bound}, \eqref{e3bound} gives 
    \begin{align*}
         - \partial_t \phi&(t_0,m_{0}) - \int_{\T^d} \Delta_x \frac{\delta \phi}{\delta m}(t_0,m_{0},x) m_{0}(dx) + \int_{\T^d} H\Big(y,D_x\frac{\delta \phi}{\delta m}(t_0,m_{0},x)\Big) m_{0}(dx) \leq \sF(m_0) + C\epsilon,    
    \end{align*}
    which completes the proof.
\end{proof}

We close this subsection with a useful lemma, which, among other things, is used to show that the viscosity subsolution property in Claim \textit{(3)} of Proposition \ref{lem:uepsproperties} can be understood in a pointwise sense.
\begin{lem} \label{lem:fullderivative}
Let $\Phi = \Phi(t,m) : [0,T] \times H^{-s} \to \R$ such that $\sup_{0 \leq t \leq T} \big[\Phi(t,\cdot) \big]_{\cC^{1,1}(H^{-s})} \leq C$ for some constant $C$, and $\nabla_{-s} \Phi$ is jointly continuous. Then, for any $(t_0,m_0) \in [0,T] \times \mathcal{P}(\T^d)$ such that the derivative $\partial_t 
    \Phi(t_0,m_0)$ exists, the function $\Phi$ is (fully) differentiable at $(t_0,m_0)$.
    As a consequence, if $\Phi$ satisfies 
\begin{align}
        - \partial_t \Phi(t,m) - \int_{\T^d} \Delta_x \frac{\delta \Phi}{\delta m}(t,m,x) m(dx) + \int_{\T^d} H\Big(y, D_x \frac{\delta \Phi}{\delta m}(t,m,x) \Big)m(dx)  \leq \sF(m) + C, 
        \label{vinequality}
    \end{align}
    in the viscosity sense on $[0,T] \times \pr(\T^d)$, with test functions in $\cC^{1}((0,T) \times H^{-s})$ (see Remark \ref{rmk:testfunctions}), then \eqref{vinequality} is satisfied at any such point $(t_0,m_0)$.
\end{lem}

\begin{proof}
We first prove that $\Phi$ is differentiable at $(t_0,m_0)$. Using the bound on $\big[\Phi(t,\cdot) \big]_{\cC^{1,1}(H^{-s})}$, the existence of $\partial_t \Phi(t_0,m_0)$, and the joint continuity of $\nabla_{-s} \Phi$ (see Proposition \ref{lem:uepsproperties}), we can write
\begin{align*}
    \Phi(t,q) &- \Phi(t_0,m_0) = \Phi(t,q) - \Phi(t,m_0) + \Phi(t,m_0) - \Phi(t_0,m_0) \\
    &= \bigl\langle \nabla_{-s} \Phi(t,m_0), q - m_0 \bigr\rangle_{-s} + \partial_t \Phi(t_0,m_0)(t-t_0) + |t - t_0|e_1(t) + \norm{m - m_0} e_2(t,q) \\
    &= \bigl\langle \nabla_{-s} \Phi(t_0,m_0), q - m_0 \bigr\rangle_{-s} + \partial_t \Phi(t_0,m_0)(t-t_0) \\
    &\qquad +\bigl\langle \nabla_{-s} \Phi(t,m_0) - \nabla_{-s} \Phi(t_0,m_0), q - m_0 \bigr\rangle_{-s} +  |t - t_0|e_1(t) + \norm{q - m_0} e_2(t,m)  \\
    &= \bigl\langle \nabla_{-s} \Phi(t_0,m_0), m - m_0 \bigr\rangle_{-s} + \partial_t \Phi(t_0,m_0)(t-t_0) \\
    &\qquad + |t - t_0|e_1(t) + \norm{m - m_0} e_2(t,q)  + e_3(t)\norm{q - m_0},
\end{align*}
where $e_1, e_3 : [0,T] \to \R$ and $e_2  : [0,T] \times H^{-s} \to \R$ are functions such that 
\begin{align*}
    e_1(t_n) \to 0, \quad e_2(t_n, q_n) \to 0, \quad e_3(t_n) \to 0
\end{align*}
as $(t_n,q_n) \to (t_0,m_0)$. 
We conclude that that $\Phi(t_0,m_0)$ is differentiable at $(t_0,m_0)$. 

By classical arguments, see e.g. \cite{Barles1994} p.18 \footnote{In fact this reference only treats the finite dimensional case but the proof extends readily to any Hilbert space}, this implies that we can find $\Psi \in \mathcal{C}^1((0,T) \times H^{-s})$ which touches $\Phi$ from above at $(t_0,m_0)$. In particular, $\partial_t \Psi(t_0,m_0) = \partial_t \Phi(t_0,m_0)$ and $D_{-s} \Psi(t_0,m_0) = D_{-s} \Phi(t_0,m_0)$. By Proposition \ref{lem:uepsproperties}, this implies that Inequality \eqref{InequationforUepsilon} is satisfied at $(t_0,m_0)$.
\end{proof}

\subsection{Analysis of $U^{\delta}$, $U^{\delta,\epsilon}$, and $U^{\delta, \epsilon, \lambda}$}
\label{subsec:uepsdelta}

This section is concerned with the $d_1$-regular case. We fix throughout this section a real number $s$ with 
\begin{align*}
    s > d/2 + 1.
\end{align*}
For each $\delta, \epsilon > 0$, and $\lambda \in (0,1)$, we define functions $U^{\delta} : [0,T] \times \pr(\T^d) \to \R$, $U^{\delta, \epsilon}, U^{\delta, \epsilon, \lambda} :[0,T] \times \pr(\T^d) \to \R$,
by
\begin{align} \label{def:udelta}
    U^{\delta}(t,m) = U(t, m * \rho_{\delta}), 
\end{align}
with $\rho_{\delta}$ being defined as in Lemma \ref{lem:convolution:rhodelta},
\begin{align} \label{def:udeltaeps}
    U^{\delta, \epsilon}(t,q) = \sup_{m \in \pr(\T^d)} \Big\{U^{\delta}(t,m) - \frac{1}{2 \epsilon} \|m - q\|_{-s}^2 \Big\},  
\end{align}
and finally 

\begin{align}
    \label{def:udeltaepslambda}
    U^{\delta, \epsilon, \lambda}(t,q) = U^{\delta, \epsilon}(t, (1 - \lambda) q + \lambda \leb ).
\end{align}
Ultimately, we would like to answer the following questions: 
\begin{itemize}
    \item How close is $U^{\delta, \epsilon, \lambda}$ to $U$?
    \item What are the regularity properties of $U^{\delta, \epsilon, \lambda}$? 
    \item By how much does $U^{\delta, \epsilon, \lambda}$ fail to be a sub-solution of the PDE \eqref{hjbinf}?
\end{itemize}

Of course, a preliminary step will be to answer the same questions for $U^{\delta}$ and then $U^{\delta, \epsilon}$.


\begin{lem}  \label{lem:udeltaprops}
    There is a constant $C>0$ independent of $\delta >0$ such that, for any $\phi \in \cC^1((0,T) \times H^{-s})$ and any $(t_0,m_0) \in (0,T) \times \mathcal{P}(\T^d)$ such that $m_0 \geq c \text{Leb}$ for some $c > 0$ and
    \begin{align*}
        U^{\delta}(t_0,m_0) - \phi(t_0,m_0) = \sup_{(t,q) \in (0,T) \times \mathcal{P}(\T^d)} \Big\{ U^{\delta}(t,m) - \phi(t,m) \Big\},
    \end{align*}
    we have 
    \begin{align}
        - \partial_t \phi(t_0,m_0) - \int_{\T^d} \Delta_{x} \frac{\delta \phi}{\delta m}(t_0,m_0,x) m_0(dx) + \int_{\T^d} H\Big(x, D_x  \frac{\delta \phi}{\delta m}(t_0,m_0,x) \Big)m_0(dx)  \leq \sF(m_0) + C\delta.
    \end{align}
\end{lem}

\begin{proof}
First we take $\phi$, $t_0$ and $m_0$ as in the statement of the lemma as well as some $h>0$ such that $t_0 + h \leq T$. We let 
\begin{equation}
\alpha(t,x) = -D_pH \Bigl(x,D_x \frac{\delta \phi}{\delta m}(t_0,m_0,x) \Bigr), \hspace{20pt} (t,x) \in [t_0,t_0+h] \times \T^d,
\label{definitionalpha26avril2023}
\end{equation}
and $m$ be the associated trajectory starting from $(t_0,m_0)$, ie solution to 
\begin{equation}
\left \{
\begin{array}{ll}
    \partial_t m(t,x) +\div_{x}(\alpha(t,x) m(t,x) ) - \Delta_x m(t,x) = 0 & \mbox{in } (t_0,t_0+h) \times \T^d \\
    m_{t_0} = m_0.
\end{array}
\right.
\end{equation}
In particular, since $\alpha \in \mathcal{C}^1([t_0,t_0+h] \times \T^d)$, $m_t$ has a density for all $t \in (t_0,t_0+h]$. Now if we take such a pair $(\alpha,m)$, then  $(m_t^{\delta} := \rho_{\delta} * m_t )_{t \in [t_0,t_0+h]}$ solves 
\begin{equation}
\left \{
\begin{array}{ll}
    \partial_t m^{\delta}(t,x) +\div_{{x}}(\alpha^{\delta}(t,x) m^{\delta}(t,x) ) - \Delta_{x} m^{\delta}(t,x) = 0 & \mbox{in } (t_0,t_0+h) \times \T^d, \\
    m^{\delta}_{t_0} = m^{\delta}_0,
\end{array}
\right.
\end{equation}
with 
$$\alpha^{\delta}(t,x) := \frac{\rho_{\delta}*(\alpha(t,\cdot) m(t,\cdot))(x)}{\rho_{\delta}*m(t,\cdot)(x)}
\mathds{1}_{\{
\rho_{\delta}*m_t(x)>0\}},$$ 
noticing that 
$\rho_{\delta}*(\alpha(t,\cdot )m(t,\cdot))(x)=0$ if 
$\rho_{\delta}*m(t,\cdot )(x)=0$ and $\alpha(t,\cdot) \in L^1({\mathbb T}^d;m_t)$.
By dynamic programming it holds
\begin{align*}
U^{\delta}(t_0,m_0)&= U(t_0,\rho_{\delta} * m_0) \\
& \leq \int_{t_0}^{t_0+h} \int_{\T^d} L(x,\alpha^{\delta}(t,x))m^{\delta}_t(dx)dt + \int_{t_0}^{t_0+h} \mathcal{F}(m^{\delta}_t)dt + U(t_0+h,m^{\delta}_{t_0+h}), 
\end{align*}
and, using the definition of $U^{\delta}$ and the fact that $\phi$ touches $U^{\delta}$ from above at $m_0$ we have
\begin{align*}
\phi(t_0,m_0) &\leq  \int_{t_0}^{t_0+h} \int_{\T^d} L(x,\alpha^{\delta}(t,x))m^{\delta}_t(dx)dt + \int_{t_0}^{t_0+h} \mathcal{F}(m^{\delta}_t)dt + U^{\delta}(t_0+h,m_{t_0+h}) \\
&+ \phi(t_0,m_0) - U^{\delta}(t_0,m_0) \\
&\leq  \int_{t_0}^{t_0+h} \int_{\T^d} L(x,\alpha^{\delta}(t,x))m^{\delta}_t(dx)dt + \int_{t_0}^{t_0+h} \mathcal{F}(m^{\delta}_t)dt + \phi(t_0+h,m_{t_0+h}).
\end{align*}
Using the joint convexity of $\mathscr{L}: (a,b) \in \R^d \times \R^+ \mapsto L \bigl(x,\frac{a}{b} \bigr)b \in \R $ (set to be $+\infty$ if b=0), we can apply Jensen's inequality and deduce that, for all $(t,x) \in  (t_0,t_0+h] \times \T^d$,
\begin{equation}
\begin{split}
L\bigl(x,\alpha^{\delta}(t,x)\bigr)
m_t^{\delta}(x)
&=
{\mathscr L} \biggl( \bigl(\alpha(t,\cdot) m_t(\cdot),m_t(\cdot)\bigr)*\rho_\delta(x) 
\biggr) 
\\
&\leq 
{\mathscr L} \bigl( \alpha(t,\cdot) m_t(\cdot),m_t(\cdot)\bigr)*\rho_\delta(x) 
 =
 \rho_{\delta} * \left[ L(x, \alpha(t,\cdot))m_t(\cdot) \right] (x).
\label{Jensen22mars2023}
\end{split}
\end{equation}
As a consequence, we have
\begin{align*}
    \int_{t_0}^{t_0+h} &\int_{\T^d} L(x,\alpha^{\delta}(t,x)) m^{\delta}_t(x)dxdt  \leq \int_{t_0}^{t_0+h} \int_{T^d} \rho_{\delta} * [L(x,\alpha(t,\cdot)m_t(.)](x)dxdt \\
    &= \int_{t_0}^{t_0+h} \int_{\T^d} \int_{\T^d} \rho_{\delta}(x-y)L \bigl(y,\alpha(t,y)\bigr)m_t(dy)dxdt \\ 
    &\hspace{15pt} + \int_{t_0}^{t_0+h}\int_{\T^d} \int_{\T^d} \rho_{\delta}(x-y) \left[L\bigl(x,\alpha(t,y)\bigr) - L\bigl(y,\alpha(t,y)\bigr) \right]m_t(dy) dxdt \\ 
    &\leq \int_{t_0}^{t_0+h} \int_{\T^d} L\bigl(x,\alpha(t,x)\bigr)m_t(dx)dt + C(\|\alpha\|_{\linf}) \int_{t_0}^{t_0+h}\int_{\T^d} \rho_{\delta}(x-y) |x-y| m_t(dy)dxdt \\
    & \leq \int_{t_0}^{t_0+h} \int_{\T^d} L \bigl(x,\alpha(t,x) \bigr)m_t(dx)dt + C(\|\alpha\|_{\linf}) \delta h,
 \end{align*}
where $C(\|\alpha\|_{\linf})$ denotes the Lipschitz constant of $L$ in the variable $x$ over the set $\T^d \times B_{\|\alpha\|_{\linf}}$. On the other hand, $\mathcal{F}$ being Lipschitz continuous with respect to $d_1$, inequality $\bigl| \mathcal{F}(\rho_{\delta} *m) - \mathcal{F}(m) \bigr| \leq C \delta$ holds for all $m\in \mathcal{P}(\T^d)$ and some $C>0$ depending only on the Lipschitz constant of $\mathcal{F}$. Therefore we have,
\begin{align}
\notag \phi(t_0,m_0) - \phi(t_0+h,m_{t_0+h}) &\leq  \int_{t_0}^{t_0+h} \int_{\T^d} L \bigl(x,\alpha(t,x) \bigr)m_t(dx)dt + \int_{t_0}^{t_0+h} \mathcal{F}(m_t)dt \\
&+ C(\|\alpha\|_{\linf}) \delta h + C \delta h.
\label{beforedividingbyh26avril}
\end{align}
Since $\alpha$ is given by \eqref{definitionalpha26avril2023} it holds
$$L\bigl(x,\alpha(t,x)\bigr) = -\alpha(t,x).D_x \frac{\delta \phi}{\delta m}(t_0,m_0,x) - H \Bigl(x,D_x \frac{\delta \phi}{\delta m}(t_0,m_0,x) \Bigr), \hspace{20pt} \forall (t,x) \in [t_0,t_0+h]\times \T^d,$$
and therefore, using the equation satisfied by $m$ we easily deduce, after dividing \eqref{beforedividingbyh26avril} by $h$ and letting $h$ tends to $0$ that
\begin{align}
\notag &- \partial_t \phi(t_0,m_0) - \int_{\T^d} \Delta_{x} \frac{\delta \phi}{\delta m}(t_0,m_0,x) m_0(dx) + \int_{\T^d} H\Big(x, D_x  \frac{\delta \phi}{\delta m}(t_0,m_0,x) \Big)m_0(dx)  \\
&\leq \mathcal{F}(m_0) + C(\|\alpha\|_{\linf}) \delta  + C \delta.
\label{Inequalitybeforeboundsonalpha27avril2023}
\end{align}
To conclude, it remains to prove that $C(\|\alpha\|_{\linf})$ is bounded independently from $\delta >0$. We are going to show that
\begin{equation} 
\bigl\| D_x \frac{\delta \phi}{\delta m}(t_0,m_0,.) \bigr \|_{\linf} \leq L_U,
\label{thedesiredinequality4avril2023}
\end{equation}
where $L_U$ is a Lipschitz constant for $U$ in the argument $m$ with respect to $d_1$. By Lemma \ref{lem:convolution:rhodelta}, $L_U$ also provides a bound for the Lipschitz constant of $U^\delta$ (in $m$ w.r.t. $d_1$). To see that indeed we have \eqref{thedesiredinequality4avril2023} we note that because $\phi$ touches $U^{\delta}$ from above at $(t_0,m_0)$, we must have 
\begin{align}
\notag \phi(t_0,m_0) &- \phi\bigl(t_0,m_0 + \epsilon D_{\vec{e}}g \leb \bigr)  \leq U^{\delta}(t_0,m_0) - U^{\delta}
\bigl(t_0,m_0+ \epsilon D_{\vec{e}}g \leb \bigr) \\
 &\leq L_U |\epsilon| \sup_{\|Df \|_{\linf} \leq 1} \int_{{\mathbb T}^d} D_{\vec{e}} g f dx = L_U |\epsilon |\sup_{\|Df\|_{\linf} \leq 1} 
 \biggl( - \int_{\T^d} g Df dx 
 \biggr) = L_U |\epsilon |\norm{g}_{L^1},
 \label{towardLipschitzestimate28mars2023}
\end{align}
for each smooth $g: {\mathbb T}^d \rightarrow {\mathbb R}$ --- where $D_{\vec{e}} g:= D_x g \cdot \vec{e}$ denotes the derivative of $g$ in the direction of some unit vector $\vec{e} \in {\mathbb R}^d$--- and each $\epsilon \in (-1,1)$ with $|\epsilon|$ small enough so that $m_0 + \epsilon D_{\vec{e}}g \leb \in \mathcal{P}({\mathbb T}^d)$ (which is indeed possible since 
$m_0 \geq c \leb$). 
Being $\phi$ in $C^1([0,T] \times H^{-s})$, it holds
\begin{equation}
\begin{split}
\phi(t_0,m_0+ \epsilon D_{\vec{e}}g \leb ) - \phi(t_0,m_0) &= \epsilon \int_{\T^d} \frac{\delta \phi}{\delta m}(t_0,m_0,x) D_{\vec{e}}g(x)dx +o(\epsilon)
\\
&= -\epsilon \int_{\T^d} \Bigl( D_x \frac{\delta \phi}{\delta m}(t_0,m_0,x)\cdot \vec{e} \Bigr) g(x)dx +o(\epsilon),
\end{split}
\end{equation}
and therefore, by arbitrariness of $\epsilon$ small enough, we deduce from \eqref{towardLipschitzestimate28mars2023} that   
\begin{align*}
\biggl|\int_{\T^d} \Bigl(D_x\frac{\delta \phi}{\delta m} (t_0,m_0,x) \cdot \vec{e} \Bigr) g(x) dx \biggr|  \leq L_U \norm{g}_{L^1},
\end{align*}
for all smooth $g$, which easily gives \eqref{thedesiredinequality4avril2023}.
\end{proof}

We now turn to the relevant properties of $U^{\delta, \epsilon}$, as defined in 
\eqref{def:udelta}.

\begin{lem} \label{lem:udelteepsproperties}
    For any $\eta > 0$ there are constants $c, C > 0$ (which can depend on $\eta$), such that for each $\delta > 0$, each $0 < \epsilon < c \delta^{2(s-1)}$ and each $(t_0,m_0) \in (0,T) \times \pr(\T^d)$ such that $m_0 \geq C \epsilon \delta^{-(2s + d/2 + \eta - 1)} \textrm{\rm \leb}$ and $U^{\delta,\epsilon}$ is differentiable at $(t_0,m_0)$, we have 
    \begin{align} \label{thresholdsub}
    &- \partial_t U^{\delta, \epsilon}(t_0,m_0) - \int_{\T^d} \Delta_x \frac{\delta U^{\delta,\epsilon}}{\delta m}(t_0,m_0,x) m_0(dx) + \int_{\T^d}H\Bigl(x, D_x \frac{\delta U^{\delta,\epsilon}}{\delta m}(t_0,m_0,x)\Bigr) m_0
    (dx) \nonumber  \\
   & \leq \sF(m_0) + C \big(\delta +  \epsilon \delta^{(2s + d/2 + \eta - 1)}\big).
    \end{align}
    Moreover, $\displaystyle \|D_{x} \frac{\delta U^{\delta, \epsilon}}{\delta m}(t_0,m_0,\cdot)\|_{\linf} \leq C$.
\end{lem}

\begin{proof}

Being $U^{\delta,\epsilon}$ differentiable at $(t_0,m_0)$,  there exists $\phi \in \cC^1((0,T) \times H^{-s})$ which touches $U^{\delta, \epsilon}$ from above at $(t_0,m_0)$
and which satisfies
\begin{align*}
\partial_t \phi(t_0,m_0) = \partial_t U^{\delta, \epsilon}(t_0,m_0), \quad \nabla_{-s} \phi(t_0,m_0) = \nabla_{-s} U^{\delta,\epsilon}(t_0,m_0).
\end{align*}

Now we let $m_{\delta,\epsilon} \in \mathcal{P}(\T^d)$ be the unique point such that 
\begin{align*}
    U^{\delta,\epsilon}(t_0,m_0) = U^{\delta}(t_0,m_{\delta,\epsilon}) - \frac{1}{2 \epsilon} \|m_0 - m_{\delta,\epsilon}\|_{-s}^2, 
\end{align*}
and we follow the reasoning in Proposition \ref{lem:uepsproperties} to conclude that because $\phi$ touches $U^{\delta,\epsilon}$ from above at $(t_0,m_0)$, the function $\psi$ defined by
\begin{align*}
\psi(t,m) = \phi(t, m_0 - m_{\delta,\epsilon} + m)
\end{align*}
touches $U^{\delta}$ from above at $(t_0,m_{\delta,\epsilon})$. In particular assuming that the constant $C$ appearing in the statement of the lemma satisfies $C>c_1$ where $c_1$ is a Lipschitz constant for $U$ with respect to $d_1$ we can infer from Proposition \ref{prop:threshold} that $m^{\epsilon,\delta}\geq c\leb$ for some $c>0$. As a consequence we can apply Lemma \ref{lem:udeltaprops} and deduce that
\begin{align*}
 &- \partial_t U^{\delta, \epsilon}(t_0,m_0) - \int_{\T^d} \Delta_{x} \frac{\delta U^{\delta,\epsilon}}{\delta m}\bigl(t_0,m_0,x\bigr) m_{\delta,\epsilon}(dx) + \int_{\T^d}H\Bigl(x, D_x \frac{\delta U^{\delta,\epsilon}}{\delta m}(t_0,m_{0},x)\Bigr) m_{\delta,\epsilon}(dx) \\
& = - \partial_t \psi(t_0,m_{\delta,\epsilon}) - \int_{\T^d} \Delta_{x} \frac{\delta \psi}{\delta m}\bigl(t_0,m_{\delta,\epsilon},x\bigr) m_{\delta,\epsilon}(dx) + \int_{\T^d}H
\Bigl(x, D_x \frac{\delta \psi}{\delta m}(t_0,m_{\delta,\epsilon},x)\Bigr) m_{\delta,\epsilon} (dx)
\\
&  \leq \sF(m_{\epsilon}) + C \delta.
\end{align*}
Thus 
\begin{align} \label{udeltaepserror} \nonumber 
 &- \partial_t U^{\delta, \epsilon}(t_0,m_0) - \int_{\T^d} \Delta_{x} \frac{\delta U^{\delta,\epsilon}}{\delta m} \bigl(t_0,m_0,x\bigr) m_{0}(dx) + \int_{\T^d}H\Bigl(x, D_{x} \frac{\delta U^{\delta,\epsilon}}{\delta m}(t_0,m_{0},x)\Bigr) m_{0}(dx) \\
 & \leq \sF(m_{0}) + C' \delta + E_1 + E_2 + E_3, 
\end{align}
where $C'>0$ depends only on the Lipschitz constants of $U$ and $\mathcal{F}$ with respect to $d_1$ and 
\begin{align*}
        E_1 &:= \int_{\T^d} \Delta_{x} \frac{\delta U^{\delta, \epsilon}}{\delta m}\bigl(t_0,m_0,x\bigr) d(m_{\delta,\epsilon} - m_0)(x), 
        \\
        E_2 &:= \int_{\T^d} H\Bigl(x,D_{x}  \frac{\delta U^{\delta, \epsilon}}{\delta m}(t_0,m_0,x)\Bigr) d(m_{0} - m_{\delta,\epsilon})(x), \\
        E_3 &:= \sF(m_{\delta,\epsilon}) - \sF(m_0). 
    \end{align*}
Lemma \ref{lem:laplacenegative} 
(together with 
Claim \textit{(3)} in Proposition 
\ref{prop:supconvproperties})
shows that $E_1 \leq 0$ exactly as in the proof of Proposition \ref{lem:uepsproperties}. To bound $E_2$, we use Proposition \ref{prop:threshold} to see that because $m_0 \geq C \epsilon \delta^{-(2s + d/2 + \eta - 1)} \leb$ (see Remark \ref{rmk:measureboundedbelow}), we have 
\begin{align*}
\|m_0 - m_{\delta,\epsilon}\|_{\linf} \leq C \delta^{-(2s + d/2 + \eta - 1)}\epsilon, 
\end{align*}
and also, as already mentioned, we have
(from Proposition \ref{prop:threshold})
$\|D_x \frac{\delta U^{\delta,\epsilon}}{\delta m}(t_0,m_0,\cdot)\|_{\linf} \leq C$. Combining these two inequalities easily gives 
\begin{align*}
    E_2 \leq C \epsilon \delta^{-(2s + d/2 + \eta - 1)}. 
\end{align*}
To bound $E_3$, we use $d_1$-Lipschitz continuity of $\sF$ to get 
\begin{align*}
E_3 \leq C d_1(m_0,m_{\delta,\epsilon}) \leq C \|m_0 - m_{\delta,\epsilon}\|_{\linf} \leq C \epsilon \delta^{-(2s + d/2 + \eta -1)}.
\end{align*}
Combining the estimates of $E_1$, $E_2$, and $E_3$ with \eqref{udeltaepserror} completes the proof.
\end{proof}

Finally, we have the properties of $U^{\delta, \epsilon, \lambda}$ defined in \eqref{def:udeltaepslambda}.

\begin{lem} \label{lem:udelteepslambda}
    For any $\eta > 0$ there are constants $c, C > 0$ which can depend on $\eta$, such that for each $\delta > 0$, each $0 < \epsilon < c \delta^{2(s-1)}$, and each $\lambda \in (0,1/2)$ such that 
    \begin{align*}
    \lambda \geq C \epsilon \delta^{-(2s + d/2 + \eta - 1)}, 
    \end{align*}
    we have 
    \begin{enumerate}
    \item $U^{\delta, \epsilon,\lambda}$ is $C \delta^{-(s-1)}$-Lipschitz over $\pr(\T^d)$ with respect to $\| \cdot \|_{-s}$ and  $C \delta^{-2(s-1)}$-semi-concave (again with respect to $\|\cdot\|_{-s}$) over all of $H^{-s}$. Moreover, for each $S \in (0,T)$, 
    the function 
    $U^{\delta,\epsilon,\lambda}$ is Lipschitz continuous in time, uniformly in 
    $(\delta,\epsilon,\lambda)$, on $[0,S] \times \pr(\T^d)$.    
    \item For each $(t,m) \in [0,T] \times {\mathcal P}({\mathbb T}^d)$, 
    \begin{align*}
        |U^{\delta, \epsilon, \lambda}(t,m) - U(t,m)| \leq C \big(\lambda + \delta + \epsilon \delta^{-2(s-1)} \big).
    \end{align*}
    \item  The map $[0,T] \times H^{-s} \ni (t,q) \mapsto \nabla_{-s}U^{\delta, \epsilon, \lambda}(t,q) \in H^{-s}$ is jointly continuous.
    \item For each $t \in [0,T]$, $U^{\delta, \epsilon}(t, \cdot)$ lies in $\cC^{1,1}(H^{-s})$ and satisfies 
    \begin{align*}
    \sup_{0 \leq t \leq T} \big[U^{\delta, \epsilon, \lambda}(t,\cdot) \big]_{\cC^{1,1}(H^{-s})} \leq \frac{C}{\epsilon}.
    \end{align*}
    \item If $(t_0,m_0) \in (0,T) \times \pr(\T^d)$ is such that $t \mapsto U^{\delta, \epsilon,\lambda}(t,m_0)$ is differentiable at $t_0$ then $U^{\delta,\epsilon,\lambda}$ has a (full) derivative at $(t_0,m_0)$ and we have
    \begin{align*}
    &- \partial_t U^{\delta, \epsilon, \lambda}(t_0,m_0) - \int_{\T^d} \Delta_{x} \frac{\delta U^{\delta,\epsilon,\lambda}}{\delta m}(t_0,m_0,x) m_0(dx) + \int_{\T^d}H\Bigl( x, D_{x} \frac{\delta U^{\delta,\epsilon,\lambda}}{\delta m}(t_0,m_0,x)\Bigr) m_0(dx) \\
&    \leq \sF(m_0) + C \big(\delta +   \lambda + \epsilon \delta^{(2s + d/2 + \eta - 1)}\big).
    \end{align*}
    \end{enumerate}
\end{lem}

\begin{proof}
Claims \textit{(1)--(4)} are easily proved using the regularity of $U$ from Proposition \eqref{prop:uregularityd1} as well as the properties of the various regularization in Lemma \eqref{lem:convolution:rhodelta} and Proposition \eqref{prop:supconvproperties}. For \textit{(5)}, we fix $\eta>0$ and we choose $C,c>0$ so that Lemma \ref{lem:udelteepsproperties} holds with the corresponding constants.

We assume that $(t_0,m_0) \in (0,T) \times \pr(\T^d)$ is such that $t \mapsto U^{\delta,\epsilon,\lambda}$ is differentiable at $t_0$. Then, by Lemma \eqref{lem:fullderivative},  $U^{\delta, \epsilon, \lambda}$ has a full derivative at $(t_0,m_0)$. Set $m_{\lambda} = \lambda \leb+ (1 - \lambda) m_0$. Then $U^{\delta, \epsilon}$ has a full derivative at $(t_0,m_{\lambda})$. If $\lambda \geq C \epsilon \delta^{- (2s + d/2 + \eta - 1)}$, then we can apply Lemma \eqref{lem:udelteepsproperties} and deduce that  \eqref{thresholdsub} holds at $(t_0,m_{\lambda})$ (even though the absolutely continuous part of $m_0$ is not lower bounded). Since 
\begin{align*}
\partial_t U^{\delta, \epsilon, \lambda}(t_0,m_0) = \partial_t U^{\delta, \epsilon}(t_0,m_{\lambda}), \quad \frac{\delta U^{\delta, \epsilon, \lambda}}{\delta m} (t_0,m_0,\cdot) = (1 - \lambda)\frac{ \delta U^{\delta, \epsilon}}{\delta m}(t_0,m_{\lambda},\cdot), 
\end{align*}
we have 
\begin{align}  \label{udelcomp}
    \nonumber &- \partial_t U^{\delta, \epsilon, \lambda}(t_0,m_0) - \int_{\T^d} \Delta_{x} \frac{\delta U^{\delta, \epsilon, \lambda}}{\delta m} (t_0,m_0,x) m_0(dx) + \int_{\T^d} H\Bigl(x, D_{x} \frac{\delta U^{\delta, \epsilon, \lambda}}{\delta m} (t_0,m_0,x)\Bigr)  m_0(dx) \\
    &= - \partial_t U^{\delta, \epsilon}(t_0,m_{\lambda}) - (1 - \lambda) \int_{\T^d} \Delta_{x} \frac{\delta U^{\delta, \epsilon}}{\delta m}(t_0,m_{\lambda},x) m_0(dx) 
    \\
    &\hspace{15pt} + \int_{\T^d} 
    H\Bigl(x, (1 - \lambda) D_{x} \frac{\delta U^{\delta, \epsilon}}{\delta m} (t_0,m_{\lambda},x)\Bigr) m_0(dx) 
    \nonumber \\
    &= - \partial_t U^{\delta, \epsilon}(t_0,m_{\lambda}) - \int_{\T^d} \Delta_{x} \frac{\delta U^{\delta, \epsilon}}{ \delta m}(t_0,m_{\lambda},x) m_\lambda(dx) - \lambda 
     \int_{\T^d} \Delta_{x} \frac{\delta U^{\delta, \epsilon}}{\delta m}(t_0,m_{\lambda},x) \leb(dx) \nonumber 
     \\
    &\quad + \frac{1}{1 - \lambda}\int_{\T^d} H\Bigl(x, (1 - \lambda) D_{x} \frac{\delta U^{\delta, \epsilon}}{\delta m} (t_0,m_{\lambda},x)\Bigr) 
    m_{\lambda} (dx) \nonumber
    \\
    &\quad + \frac{\lambda}{1 - \lambda} \int_{\T^d} H\Bigl(x, (1 - \lambda) D_{x}
    \frac{\delta U^{\delta, \epsilon}}{\delta m} (t_0,m_{\lambda},x)\Bigr) \leb(dx).\nonumber
\end{align}
To go further, we first use integration by parts to conclude that
\begin{align} \label{ibp}
\int_{\T^d} \Delta_{x} \frac{\delta U^{\delta, \epsilon}}{\delta m}(t_0,m_{\lambda},x) \leb(dx) = 0.
\end{align}
Next, we use convexity of the Hamiltonian in the second argument to get 
\begin{align*}
H\Bigl(x, (1 - \lambda) D_{x} \frac{\delta U^{\delta,\epsilon}}{\delta	 m}(t_0,m_{\lambda},x)\Bigr) &\leq (1 - \lambda) H\Bigl(x, D_{x} \frac{\delta U^{\delta,\epsilon}}{\delta m}(t_0,m_{\lambda},x)\Bigr) + \lambda H(x,0)  \\
&\leq (1 - \lambda) H\Bigl( x, D_{x} 
\frac{\delta U^{\delta,\epsilon}}{\delta m}(t_0,m_{\lambda},x)\Bigr) + C \lambda, 
\end{align*}
and hence 
\begin{align}  \label{hconvex}
&\frac{1}{1 - \lambda}\int_{\T^d} H\Bigr(x, (1 - \lambda) D_{x} \frac{\delta U^{\delta, \epsilon}}{\delta m} (t_0,m_{\lambda},x)\Bigr) m_{\lambda}(dx) 
\\
&\quad \leq \int_{\T^d} H\Bigl(x,D_{x} \frac{\delta U^{\delta, \epsilon}}{\delta m} (t_0,m_{\lambda},x)\Bigr) m_{\lambda}(dx) 
 + C \lambda. \nonumber
\end{align}
Similarly, using the bound on $ \|D_{x} \frac{\delta U^{\delta,\epsilon}}{\delta m}(t,m_{\lambda},\cdot)\|_{\linf}$ from Lemma \ref{lem:udelteepsproperties}, we have 
\begin{align} \label{derivlinf}
    \frac{\lambda}{1 - \lambda} \int_{\T^d} H\Bigl(x, (1 - \lambda) D_{x} \frac{\delta U^{\delta, \epsilon}}{\delta m} (t_0,m_{\lambda},x)\Bigr) \leb(dx) \leq C \lambda. 
\end{align}
Combining \eqref{ibp}, \eqref{hconvex}, and \eqref{derivlinf} with \eqref{udelcomp} and using Lemma \ref{lem:udelteepsproperties}, we arrive at
\begin{align}  \label{udelcomp2}
    \nonumber &- \partial_t U^{\delta, \epsilon, \lambda}(t_0,m_0) - \int_{\T^d} \Delta_{x} \frac{\delta U^{\delta, \epsilon, \lambda}}{\delta m}(t_0,m_0,x) m_0(dx) + \int_{\T^d} 
    H\Bigl(x, D_{x} \frac{\delta U^{\delta, \epsilon, \lambda}}{\delta m}(t_0,m_0,x)\Bigr)  m_0(dx) \\
    &\leq  - \partial_t U^{\delta, \epsilon}(t_0,m_{\lambda}) - \int_{\T^d} \Delta_{x} \frac{\delta U^{\delta, \epsilon}}{\delta m}(t_0,m_{\lambda},x) m_\lambda(dx) + \int_{\T^d} H
    \Bigr(x, D_{x} \frac{\delta U^{\delta, \epsilon}}{\delta m} (t_0,m_{\lambda},x)\Bigr) m_{\lambda}(dx) + C\lambda \nonumber \\
    &\leq \sF(m_{\lambda}) + C \Bigl( \lambda+
    \delta +  \epsilon \delta^{(2s + d/2 + \eta - 1)}\Bigr)
     \leq \sF(m_0) + C \Bigl( \lambda+
    \delta +  \epsilon \delta^{(2s + d/2 + \eta - 1)}\Bigr), 
\end{align}
where the last line uses the fact that $\sF$ is $d_1$-Lipschitz. This completes the proof.

\end{proof}

\subsection{Projections and rates of convergence}

The goal of this section is to use the properties of $U^{\epsilon}$ and $U^{\delta, \epsilon, \lambda}$ proved in the previous section to establish estimates on the distance between their projections and the finite-dimensional value functions $\Phi^N$. The main technical tool we use to do this is the following proposition.

\begin{prop} \label{prop:projections}
    Suppose that $\Phi = \Phi(t,q) : [0,T] \times H^{-s} \to \R$ is a continuous map such that 
    \begin{itemize}
        \item $\Phi$ is uniformly continuous on $[0,T] \times \pr(\T^d)$ ($\pr(\T^d)$ being endowed with the metric induced by $\| \cdot \|_{-s}$) and Lipschitz (for the same metric) on each set of the form $[0,t] \times \pr(\T^d)$ with $t < T$. 
        
        \item $\sup_{0 \leq t \leq T} \big[\Phi(t,\cdot)\big]_{\cC^{1,1}(H^{-s})} < \infty$.
        \item The derivative $\nabla_{-s} \Phi$ is jointly continuous on $[0,T] \times \pr(\T^d)$.
        \item The inequality
        \begin{align}
        - \partial_t \Phi(t,m) - \int_{\T^d} \Delta_{x} \frac{\delta \Phi}{\delta m}(t,m,x) m(dx) + \int_{\T^d} H\Big(x, D_{x} \frac{\delta \Phi}{\delta m}(t,m,x) \Big)m(dx)  \leq \sF(m) + c_1, 
        \label{vinequality2}
    \end{align}        
    is satisfied at any $(t,m) \in [0,T] \times \pr(\T^d)$ where $\Phi$ has a full derivative, and $\Phi(T,m) \leq \sG(m) + c_2$, for two positive constants $c_1$ and $c_2$.
    \end{itemize}
    Define for each $N \in \N$ a function $\Phi^N(t,\bx) : [0,T] \times (\T^d)^N \rightarrow {\mathbb R}$ by 
    \begin{align*}
        \Phi^N(t,\bx) = \Phi(t,m_{\bx}^N). 
    \end{align*}
    Then we have 
    \begin{align*}
    \Phi^{N}(t,\bx) \leq V^N(t,\bx) + T\Big(c_1 + \frac{C}{N} \sup_r \big[\Phi(r,\cdot) \big]_{\cC^{1,1}(H^{-s})} \Big) + c_2, 
    \end{align*}
    for each $N \in \N$ and each $(t,\bx) \in [0,T] \times (\T^d)^N$, and a constant $C$ which depends only on $d$ and $s$.
\end{prop}

Before moving to the proof of Proposition \ref{prop:projections}, we explain how, together with the results of the previous section, it can be used to easily establish the “hard inequalities", from Theorems \ref{thm:mainhs} and \ref{thm:maind1}.

\begin{prop}[Upper bound in Theorem \ref{thm:mainhs}] \label{prop:hardinequality}
    Let Assumption \ref{assump:Hs} hold. Then there is a constant $C$ such that for each $N \in \N$, 
    \begin{align*}
        U(t,m_{\bx}^N) \leq V^N(t,m_{\bx}^N) + \frac{C}{\sqrt{N}},
    \end{align*}
    for all $(t,\bx) \in [0,T] \times (\T^d)^N$. 
\end{prop}

\begin{proof}
    For each $\epsilon$ and $N$, we combine Proposition \ref{prop:projections}, Proposition \ref{lem:uepsproperties}, and Proposition \ref{prop:supconvproperties} to conclude that
    \begin{align*}
        U(t,m_{\bx}^N) - V^N(t,\bx) &=  U(t,m_{\bx}^N) - U^{\epsilon}(t,m_{\bx}^N) + U^{\epsilon}(t,m_{\bx}^N) - V^N(t,\bx), \\
        &\leq C \Big(\epsilon + \frac{1}{N\epsilon} \Big).
    \end{align*}
    We choose $\epsilon = \frac{1}{\sqrt{N}}$ to get the result.
\end{proof}

\begin{prop}[Upper bound in Theorem \ref{thm:maind1}]
 \label{prop:hardinequalityd1}
    Let Assumption \ref{assump:Hs} hold. Then for each $\eta > 0$, there is a constant $C$ such that for each $N \in \N$ and $(t,\bx) \in [0,T] \times (\T^d)^N$,
    \begin{align*}
        U(t,m_{\bx}^N) \leq V^N(t,m_{\bx}^N) + C N^{ - \beta(d) + \eta}, \quad \beta(d) := \frac{2}{3d + 6}.
    \end{align*}
\end{prop}

\begin{proof}[Proof of Proposition \ref{prop:hardinequalityd1}]
We use Proposition \ref{prop:projections} and Lemma \ref{lem:udelteepslambda} to conclude that for any $\eta > 0$, there are constants $c$ and $C$ such that for any $\delta > 0$, $0 < \epsilon < c \delta^{2(s-1)}$, and any $\lambda \geq C \epsilon \delta^{- (2s + d/2 + \eta - 1)}$, we have 
\begin{align*}
    U(t,m_{\bx}^N) - V^N(t,\bx) &= U(t,m_{\bx}^N) - U^{\delta, \epsilon, \lambda}(t,m_{\bx}^N) + U^{\delta, \epsilon, \lambda}(t,m_{\bx}^N) - V^N(t, \bx) \\
    &\leq C\Big(\lambda + \delta + \epsilon\big( \delta^{ - 2(s-1)} + \delta^{-(2s + d/2 + \eta - 1)} \big) + \frac{1}{N \epsilon} \Big) \\ 
    &\leq  C\Big(\lambda + \delta + \epsilon \delta^{-(2s + d/2 + \eta - 1)} + \frac{1}{N \epsilon} \Big). 
\end{align*}
We choose $\lambda = C \epsilon \delta^{- (2s + d/2 + \eta - 1)}$ (with the $C$ here greater than $C$ in the statement of Lemma 
 \ref{lem:udelteepslambda})
 to find 
\begin{align*}
U(t,m_{\bx}^N) - V^N(t,\bx) \leq  C\Big(\delta + \epsilon \delta^{-(2s + d/2 + \eta - 1)} + \frac{1}{N \epsilon} \Big), 
\end{align*}
and then choose 
\begin{align*}
 \delta = N^{-1/(2s + d/2 + \eta  +1)}, \quad \epsilon = N^{-1} \delta^{-1}. 
\end{align*}
This yields 
\begin{align*}
U(t,m_{\bx}^N) - V^N(t,\bx) \leq  C N^{- 1/ (2s + d/2 + \eta + 1)}.
\end{align*}
Since $s$ was chosen to be an arbitrary real number with $s > d/2 + 1$, we find that for any $\alpha > 0$, there is a constant $C$ such that 
\begin{align*}
U(t,m_{\bx}^N) \leq V^N(t,\bx) + C N^{ -\frac{1}{(\frac{3d}{2} + 3 + \alpha)}}, 
\end{align*}
for all $(t,\bx) \in [0,T] \times (\T^d)^N$, as claimed. 
\end{proof}

The rest of this subsection is devoted to a proof of Proposition \ref{prop:projections}. We start with two lemmas.

\begin{lem} \label{lem:phinreg}
Let $\Phi$ satisfy the hypotheses of Proposition \ref{prop:projections} and $\Phi^N$ be as in the statement of the same proposition. Then we have
\begin{enumerate}
    \item $\Phi^N$ is uniformly continuous on $[0,T] \times (\T^d)^N$ and $\Phi^N \in \rm Lip\big([0,t] \times (\T^d)^N\big)$ for each $ t < T$. 
    \item $\Phi^N$ is continuously differentiable in $\bx$, with 
\begin{align*}
    D_{x^i} \Phi^N(t,\bx) = \frac{1}{N} D_{m} \Phi(t,m_{\bx}^N,x^i) = \frac{1}{N} D_{y} \big[D_{-s} \Phi(t,m_{\bx}^N)\big](x^{i}),\quad 
    i \in \{1,\cdots,N\},
\end{align*}
where we recall that $x^i$ denotes the $i^{\it th}$ $d$-dimensional entry of the $N$-tuple $\bx$.
\item The derivative $D_{x^i}\Phi^N$ is Lipschitz continuous with respect to $\bx$, uniformly in $t$.
\item The derivatives $\partial_t \Phi^N$, $D_{x^i} \Phi^N$, $D^2_{x^ix^j} \Phi^N$ exist almost everywhere, and define versions of the corresponding weak derivatives of $\Phi^N$. {Moreover the spatial derivatives $D_{x^i} \Phi^N$ and $D^2_{x^ix^j} \Phi^N$ are uniformly bounded on $[0,T] \times \T^d$ and the time derivative $\partial_t \Phi^N$ is uniformly bounded on $[0,t] \times \T^d$ for each $t < T$,} so that in particular 
\begin{align*}
    \Phi^N \in W^{1,2}_{\linf}([0,t] \times \T^d),
\end{align*}
for each $t < T$, with $\Phi^N$, with $W^{1,2}_{\linf}\big([0,t] \times (\T^d)^N\big)$.
\end{enumerate}
\end{lem}
{We note that in Claim \textit{(4)} we use the notation $W^{1,2}_{\linf}\big([0,t] \times (\T^d)^N\big)$ for the set of functions $u : [0,t] \times  (\T^d)^N \to \R$ with uniformly bounded weak derivatives $\partial_t u$, $D_{x_i} u$, $i = 1,...,N$, and $D_{x_ix_j}^{2} u$, $i,j = 1,...,N$}

\begin{proof}
  Claim \it (1) follows from the fact that $\|m_{\bx}^N - m_{\by}^N\|_{-s} \leq Cd_1(m_{\bx}^N, m_{\by}^N) \leq \frac{C}{N} \sum_{i = 1}^N |x^i - y^i|$. 
  
  Claim \it (2) follows easily from Lemma \ref{lem:derivativeconnection}. For \it (3), we use \it (2), Proposition \ref{prop:dmlipschitz}, the bound on $[\Phi(t,\cdot)]_{C^{1,1}(H^{-s})}$, and also the fact that $D_y D_m \Phi$ is bounded to conclude that 
  \begin{align*}
      |D_{x^i} \Phi^{N}(t,\bx) &- D_{x^i}\Phi^{N}(t,\by)| = \frac{1}{N}|D_m \Phi(
      t,
      m_{\bx}^N,x^i) - D_m \Phi(t,m_{\by}^N,y^i)| \\
      &\leq C |D_m \Phi(t,m_{\bx}^N,x^i) - D_m \Phi(t,m_{\by}^N,x^i)| +  C |D_m \Phi(t,m_{\by}^N,x^i) - D_m \Phi(t,m_{\by}^N,y^i)| \\
      &\leq C |\bx - \by|. 
  \end{align*}
Finally claim \it (4) follows \it (1)  and \it (3) together with Rademacher's Theorem.
\end{proof}

\begin{lem} 
\label{lem:dxixierror}
Let $\Phi$ and $\Phi^N$ be as in the statement of Proposition \ref{prop:projections}. Then for almost every $(t,\bx) \in [0,T] \times (\T^d)^N$,
the function
\begin{align} \label{dmumap}
    \T^d \ni x \mapsto D_m \Phi\Bigl(t,\frac{1}{N} \sum_{j \neq i} \delta_{x^j} + \frac{1}{N} \delta_x,x^i\Bigr)
\end{align}
is differentiable at $x^i$. At any point $(t,\bx)$ such that the map in \eqref{dmumap} is differentiable at $x^i$, the second derivative $D^{2}_{x^ix^i} \Phi^N$ exists (in the classical sense) and satisfies
\begin{align*}
    D^2_{x^i x^i} \Phi^N(t,\bx) = \frac{1}{N} D_yD_m \Phi(t,m_{\bx}^N, x^i) + \frac{1}{N^2} R^{N,i}(t,\bx),
\end{align*}
where $R^{N,i}$ is defined explicitly as 
\begin{align} \label{rdef}
    R^{ N, i}(t,\bx) = N D_{y}\Big[ \T^d \ni y \mapsto D_m\Phi\Bigl(
    t,\frac{1}{N} \sum_{j \neq i} \delta_{x^j} + \frac{1}{N} \delta_y,x^i \Bigr)
    \Big](x^i)
\end{align}
and satisfies $\norm{R^{N,i}}_{\linf} \leq 
 C \sup_{0 \leq t \leq T} [\Phi(t,\cdot)]_{\cC^{1,1}(H^{-s})}$, for $C$ depending only on $s$ and $d$.
\end{lem}

\begin{rmk}
\label{rmk:RNi:D2m}
Notice that if $\Phi$ in the above statement were twice differentiable (in the variable $m$), then 
$R^{N,i}(t,\bx)$ would be equal to 
\begin{equation*}
R^{N,i}(t,\bx) = D^2_{mm} \Phi(t,m_{\bx}^N, x^i,x^i). 
\end{equation*} 
The very much spirit of Lemma 
\ref{lem:dxixierror} is that
we are still able to provide an almost everywhere formula for $R^{N,i}(t,\bx)$ even though 
$\Phi$ just satisfies the assumptions of
Proposition 
\ref{prop:projections} (under which $\Phi$ may not be twice differentiable).
\end{rmk}

\begin{proof}
   The first claim follows easily from the Lipschitz continuity of $m \mapsto D_m \Phi(t,m,y)$ (uniformly in $t$ and $y$), see Remark \ref{RemarkLipDm}.
    Recall from Lemma \ref{lem:phinreg} that
    \begin{align*}
        D_{x^i} \Phi^N(t,\bx) = \frac{1}{N} D_{m} \Phi(t,m_{\bx}^N,x^i). 
    \end{align*}
 However, $y \mapsto D_{m} \Phi(t,m_{\bx}^N,y) = D_{y}\left[D_{-s}\Phi(t,m_{\bx}^N) \right](y)$ is continuously differentiable since 
 $D_{-s}\Phi(t,m_{\bx}^N)(\cdot) \in H^s \subset C^2(\T^d)$. And so it is clear that if $y \mapsto D_m \Phi(t,\frac{1}{N} \sum_{j \neq i} \delta_{x^j} + \frac{1}{N} \delta_y,x^i)$ is differentiable at $x^i$, then $D^{2}_{x^ix^i} \Phi(t,\bx)$ exists and is given by
    \begin{align*}
        D_{x^ix^i}^2 \Phi(t,\bx) = \frac{1}{N} D_y D_m \Phi(t, m_{\bx}^N, x^i) + \frac{1}{N} D_{y} \Big[  y \mapsto D_m \Phi \Bigl(t,\frac{1}{N} \sum_{j \neq i} \delta_{x^j} + \frac{1}{N} \delta_y,x^i \Bigr) \Big](x^i), 
    \end{align*}
    as required. The $\linf$ bound on $R^{N,i}$ follows easily from the fact that by Proposition \ref{prop:dmlipschitz}, for all $\bx, \bar{\bx} \in (\T^d)^N$,
    \begin{align*}
        |D_m \Phi(t_0, m_{\bx}^N,y) - D_m \Phi(t_0,m_{\overline{\bx}}^N,y)| \leq C [\Phi(t_0,\cdot)]_{C^{1,1}} d_1(m_{\bx}^{N}, m_{\overline{\bx}}^N) \leq C[\Phi(t_0,\cdot)]_{C^{1,1}} \frac{1}{N} \sum_{i = 1}^N |x^i - \overline{x}^i|. 
    \end{align*}
    
\end{proof}

We are now ready to combine the previous two lemmas to prove Proposition \ref{prop:projections}

\begin{proof}[Proof of Proposition \ref{prop:projections}]
We start by proving that for each $N$, we have, for almost every 
$(t,\bx) \in [0,T] \times ({\mathbb T}^d)^N$,
\begin{align} \label{uepsnsubsol}
    &- \partial_t \Phi^N(t,\bx) - \sum_{i = 1}^N \Delta_{x^i} \Phi^N(t,\bx) + \frac{1}{N} \sum_{i = 1}^N H\Bigl(x^i, ND_{x^i}\Phi^N(t,\bx)\Bigr) 
    \\
    &\quad \leq \sF(m_{\bx}^N) + \frac{C}{N} \sup_{r \in [0,T]} [\Phi(r, \cdot)]_{\cC^{1,1}(H^{-s})}.
    \nonumber
\end{align} 
To do this, we start by observing that for almost every $(t,\bx)$, we have 
    \begin{enumerate}
        \item the time derivative $\partial_t \Phi^N$ exists at $(t,\bx)$, 
        \item the map
\begin{align}
    \T^d \ni y \mapsto D_m \Phi\Bigl(t,\frac{1}{N} \sum_{j \neq i} \delta_{x^j} + \frac{1}{N} \delta_y,x^i\Bigr)
\end{align}
is differentiable at $x^i$.
    \end{enumerate}
    But at any point $(t,\bx)$ satisfying both of these conditions, Lemma \ref{lem:fullderivative} guarantees that $\Phi$ is differentiable at $(t,\bx)$, and so we can use Lemma \ref{lem:dxixierror} to get
\begin{align*}
        &- \partial_t \Phi^{N}(t,\bx) - \sum_{i = 1}^N \Delta_{x^i} \Phi^{N}(t,\bx) + \frac{1}{N} \sum_{i = 1}^N H\Bigl(x^i, N D_{x^i} \Phi^{N}(t,\bx)\Bigr) \\
        &= - \partial_t \Phi(t,m_{\bx}^N) - \sum_{i = 1}^N \Big(\frac{1}{N} \Delta_x \frac{\delta \Phi}{\delta m}(t,m_{\bx}^N,x^{i}) 
        \\
        &\qquad \qquad + \frac{1}{N^2} R^{N,i}(t,\bx) \Big) + \frac{1}{N} \sum_{i = 1}^N H\bigl(x^i, D_m\Phi(t,m_{\bx}^N,x^{i})\bigr)
        \\
        &= - \partial_t \Phi(t,m_{\bx}^N) - \int_{\T^d} \Delta_x \frac{\delta \Phi}{\delta m} (t,m_{\bx}^N,y) m_{\bx}^N(dy) \\
        &\qquad \qquad + \int_{\T^d} H\bigl(y,D_m\Phi(t,m_{\bx}^N,y)\bigr)m_{\bx}^N(dy) -  \frac{1}{N^2} \sum_{i=1}^{N} R^{N,i}(t,\bx) \\
        &\leq \sF(m_{\bx}^N) + \frac{C}{N} \sup_{r \in [0,T]} [\Phi(r,\cdot)]_{\cC^{1,1}(H^{-s})}, 
    \end{align*}
    where we used the bound on $\norm{R^{N,i}}_{\linf}$ from Lemma \ref{lem:dxixierror}. This completes the proof of \eqref{uepsnsubsol}. 
    
    At this point, we know that for each $t < T$, $\Phi^{N}$ lies in the space $W^{1,2}_{\infty}\big([0,t) \times (\T^d)^N\big)$
of functions with $\linf$ weak derivatives $\partial_t \Phi^N$, $D \Phi^N$, $D^2 \Phi^N$, and also that $\Phi^N$ satisfies the inequality \eqref{uepsnsubsol} at almost every $(t,\bx)$. Moreover, we have by assumption
    \begin{align*}
        \Phi^N(T,\bx) \leq V^{N}(T,\bx) + c_2 = \sG(m^N_\bx) + c_2.
    \end{align*}
    The rest of the proof is a standard “verification"-type argument. Indeed, we fix $(t_0,\bx_0)$ and an admissible control
    $\bm \alpha = (\alpha^1,...,\alpha^N)$ in feedback form, and we define $\bX = (X^1,...,X^N)$ by
    \begin{align*}
        dX_t^i = \alpha^i(t,X_t^i) dt + dW_t^i, \quad t \in [t_0,T]; \quad  X_{t_0}^i = x_0^i.
    \end{align*}
    We use the It\^o-Krylov formula, see \cite{Krylov1980} Section 2.10, to see that 
    \begin{align*}
        \frac{d}{dt} \E\big[\Phi^{N}(t,\bX_t)\big] &= \E\Big[\partial_t \Phi^{N}(t,\bX_t) + \sum_{i = 1}^N \Delta_{x^i} \Phi^{N}(t,\bX_t) + \sum_{i = 1}^N D_{x^i} \Phi^N(t,\bX_t) \cdot \alpha^i(t,\bX_t) \Big] \\
        &\geq 
        \E\bigg[\frac{1}{N} \sum_{i = 1}^N \Big(H\bigl(X_t^i, ND_{x^i}\Phi^N(t,\bX_t) \bigr)+ ND_{x^i} \Phi^N(t,\bX_t) \cdot \alpha^i(t,\bX_t) \Big) \\ &\qquad  - \sF(m_{\bX_t}^N) - c_1 - C \sup_{r \in [0,T]} [\Phi(r,\cdot)]_{\cC^{1,1}(H^{-s})} \bigg]
        \\
        &\geq -\E\bigg[\frac{1}{N} \sum_{i = 1}^N L(X_t^i, \alpha_t^i)  + \sF(m_{\bX_t}^N)- c_1 - C \sup_{r \in [0,T]} [\Phi(r,\cdot)]_{\cC^{1,1}(H^{-s})} \bigg].
    \end{align*}
    Thus, integrating between $t_0$ and $T$ the above inequality we get  
    \begin{align*}
        \Phi^N(t_0,\bx_0) &\leq \E\bigg[\Phi^N(T,\bX_T) + \int_{t_0}^T \Big(\frac{1}{N} \sum_{i = 1}^N L(X_t^i, \alpha_t^i)  + \sF(m_{\bX_t}^N) \Big)dt \bigg] +  T \Big(c_1 + C \sup_r [\Phi(r,\cdot)]_{\cC^{1,1}(H^{-s})} \Big) \\
        &\leq \E\bigg[\sG(\bX_T) + c_2 + \int_{t_0}^T \Big(\frac{1}{N} \sum_{i = 1}^N L(X_t^i, \alpha_t^i)  + \sF(m_{\bX_t}^N) \Big) dt\bigg] + T \Big(  c_1 + C \sup_r [\Phi(r,\cdot)]_{\cC^{1,1}(H^{-s})} \Big).
    \end{align*}
    Taking the infimum over $\alpha$ gives the result.
\end{proof}

\section{Proofs of the “easy inequalities"} \label{sec:easyinequalities}

\subsection{The “easy inequality" in the $d_1$-regular case.}

We start with a Lemma stating that $V^N$ satisfies uniform in $N$ Lipschitz estimates. The proof is essentially the same as Lemma 3.1 in \cite{cardaliaguet2023algebraic}, and so is omitted.

\begin{lem} \label{lem:uniforminn}
Let Assumption \ref{assump:d1} hold, except for possibly $d_1$-semiconcavity of $\sF$ and $\sG$. Then there is a constant $C$ such that for each $N \in \N$, $t \in [0,T]$ and $\bx, \by \in (\T^d)^N$,
\begin{align*}
|V^N(t,\bx) - V^N(t,\by)| \leq \frac{C}{N} \sum_{i = 1}^N |x^i - y^i|.
\end{align*}
\end{lem}

\begin{prop} \label{prop:easyinequalityd1}
    Let Assumption \ref{assump:d1} hold, except for possibly $d_1$-semiconcavity of $\sF$ and $\sG$. Then there is a constant $C$ such that for each $N \in \N$, 
    \begin{align*}
        V^N(t,\bx) \leq U(t,m_{\bx}^N) + CR_{N,d},
    \end{align*}
    for all $(t,\bx) \in [0,T] \times (\T^d)^N$, where $R_{N,d}$ is as defined in \eqref{fournierguillinrate}.
\end{prop}

\begin{proof}
    We first define a lift $\hat{V}^N = \hat{V}^N(t,m) : [0,T] \times \pr(\T^d) \to \R$ by
    \begin{align*}
        \hat{V}^N(t,m) = \int_{(\T^d)^N} V^N(t,\bx) 
       m^{\otimes N}(d\bx). 
    \end{align*}
    Using the Markov structure of the control problem,
    it is standard to show that $\hat{V}^N$ satisfies 
    \begin{align} \label{vhatchar}
        \hat{V}^N(t,m) = \inf_{\bm \alpha \in \sA^N} \E \bigg[ \int_{t}^T \Big( \frac{1}{N} \sum_{i = 1}^N L(X_s^i, \alpha^i_s) + \sF(m_{\bX_s}^N) \Big) ds + \sG(m_{\bX_T}^N) \bigg]
    \end{align}
    with
    \begin{align*}
        dX_s^i = \alpha_s^i ds + dW_s^i, \quad s \in [t,T]; \quad X_t^i = \xi^i, 
    \end{align*}
    with $(\xi^i)_{i = 1,...,N}$ i.i.d. with common law $m$. Now let $(\alpha,m)$ be a candidate for the optimization problem defining $U(t,m)$, and suppose that $\alpha$ is bounded so that 
     $m_s = \sL(Y_s)$, for $s \in [t,T]$, where $Y$ is the unique strong solution of 
     \begin{align*}
     dY_s = \alpha(s,Y_s) ds + dW_s, \quad s \in [t,T]; \quad X_t^i = \xi \sim m.
     \end{align*}
     Next define $\bY = (Y^1,...,Y^N)$ by 
    \begin{align*}
        &dY_s^i = \alpha(s,Y^i_s) ds + dW^i_s, \quad s \in [t,T]; \quad X_t^i = \xi^i, 
    \end{align*}
    where $(\xi^i)_{i = 1,...,N}$ are as above. Then using \eqref{vhatchar}, we have
    \begin{align*}
        \hat{V}^N(t,m) &\leq \E \bigg[ \int_{t}^T \Big( \frac{1}{N} \sum_{i = 1}^N L\bigl(Y_s^i, \alpha(s,Y_s^i)\bigr) + \sF(m_{\bY_s}^N) \Big) ds + \sG(m_{\bY_T}^N) \bigg] \\
        &= \E \bigg[ \int_{t}^T \Bigl( L
        \bigl(Y_s, \alpha(s,Y_s)\bigr) + \sF(m_{\bY_s}^N) \Big) ds + \sG(m_{\bY_T}^N) \bigg] \\
        &\leq \E \bigg[ \int_{t}^T \Bigl( L(Y_s, \alpha(s,Y_s)) + \sF\bigl(\sL(Y_s)\bigr) \Big) ds + \sG\bigl(\sL(Y_T)\bigr) \bigg] + CR_{N,d},
    \end{align*}
    with the last inequality following from \cite{FournierGuillin}. Taking an infimum in $\alpha$, we conclude that 
    \begin{align} \label{vhatlower}
        \hat{V}^N(t,m) \leq U(t,m) + CR_{N,d}.
    \end{align}
    Next, we use the fact that $V^N$ is symmetric in $(x^1,...,x^N)$ (the entries of $\bx$) and satisfies the estimate in Lemma \ref{lem:uniforminn} to get
    \begin{align}  \label{vhatvest} 
        \Bigl\vert        \hat{V}^N(t,m^N_{\bx}) -V^N(t,\bx) \Bigr\vert &=
        \biggl\vert \int_{(\T^d)^N} 
        \Bigl[ V^N(t,\by) - V^N(t,\bx) \Bigr] 
        m^{\otimes N}_{\bx}(d\by)\biggr\vert
       \nonumber  \\
        &\leq C      
        \int_{({\mathbb T}^d)^N} \Bigl[
         \frac{1}N
         \inf_{\varsigma \in {\mathcal S}_N}  
            \sum_{i=1}^N \vert y^{\varsigma(i)} -   x^i \vert 
            \Bigr]
          m^{\otimes N}_{\bx}(d\by)
          \nonumber   \\
           &= C 
\int_{({\mathbb T}^d)^N} d_1(m^N_{\by},m^N_{\bx})
        m^{\otimes N}_{\bx}(d\by) \leq  C R_{N,d}, 
    \end{align}
    with the last ineqality coming from \cite{FournierGuillin}.
    Combining \eqref{vhatvest} with \eqref{vhatlower} completes the proof.
\end{proof}

\begin{proof}[Proof of Theorem \ref{thm:maind1}]
    Combine Propositions \ref{prop:easyinequalityd1} and \ref{prop:hardinequalityd1}.
\end{proof}

\subsection{An auxiliary estimate for linear PDEs on the Wasserstein space}
In order to prove the lower bound in Theorem \ref{thm:mainhs}, we need to first prove a similar bound for a class of linear PDEs on the Wasserstein space. We fix two functions
\begin{align*}
    \alpha = \alpha(t,x) : [0,T] \times \T^d \to \R, \quad \phi = \phi(m) : \mathcal{P}(\T^d) \to \R.
\end{align*}
The following Assumption will be in force throughout the Subsection.
\begin{assumption} \label{assump:section6}
    We assume that $s$ is a real number with $s > d/2 + 1$, and that $\alpha = \alpha(t,x) : [0,T] \times \T^d \to \R$ is jointly continuous and satisfies 
    \begin{align*}
        \sup_{0 \leq t \leq T} \|\alpha(t,\cdot)\|_{s-1} < \infty.
    \end{align*}
\end{assumption}
We will consider the maps  
\begin{align}
\label{eq:PhiN:Phiinfty:1}
    \Phi^N = \Phi^N(t,\bx) : [0,T] \times (\T^d)^N \to \R, \quad \Phi^{\infty} = \Phi^{\infty}(t,\mu) : [0,T] \times \pr(\T^d) \to \R,
\end{align}
defined by 
\begin{align}
\label{eq:PhiN:Phiinfty:2}
\Phi^N(t,\bx) = \E\big[\phi(m^N_{\bX^{t,\bx}_T})\big], 
\quad 
\text{and} \quad \Phi^{\infty}(t,\mu) = \phi(m_T^{t,\mu}),
\end{align}
where $(\bX^{t,\bx}_s)_{t \leq s \leq T}$ denotes the $N$-tuple 
$(X^1_s,...,X^N_s)_{t \leq s \leq T}$ solving the system of SDEs
\begin{align*}
dX_s^i = \alpha(s,X_s^{i}) ds + \sqrt{2} dW_s^i, \quad s\in [t,T]; \quad X_{t}^i = x^i;
\quad i=1,...,N,
\end{align*}
and $(m_s)_{t \leq s \leq T}$ denotes the solution to the Fokker Planck equation
\begin{align*}
\partial_s m_s = \Delta_{x} m_s - \text{div}_{x}\bigl(m_s \alpha(s,\cdot)\bigr),
\quad  s \in [t,T];
 \quad m_t = \mu.
\end{align*}
Formally, we expect (by the method of characteristics) that $\Phi^N$ solves the PDE
\begin{align} \label{phinpde}
    \begin{cases} 
   \displaystyle - \partial_t \Phi^N(t,\bx) - \sum_{i=1}^{N} \Delta_{x^i} \Phi^N(t,\bx) -\sum_{i=1}^{N}\alpha(s,x^{i}) \cdot D_{x^{i}} \Phi^N(t,\bx)= 0, \quad (t,\bx) \in [0,T) \times (\T^d)^N, \vspace{.1cm} \\
    \Phi^N(T,\bx) = \phi(m_{\bx}^N), \quad \bx \in (\T^d)^N,
    \end{cases}
\end{align}
and that $\Phi^{\infty}$ solves the PDE 
\begin{align} \label{phiinfpde}
    \begin{cases} 
\displaystyle    - \partial_t \Phi^{\infty}(t,m) - \int_{\T^d} \Delta_{y} \frac{\delta \Phi^{\infty}}{\delta m}(t,m,y) m(dy) - \int_{\R^d} \alpha(t,y) \cdot D_y \frac{\delta \Phi^{\infty}}{\delta m}(t,m,y)m(dy) = 0, \vspace{.2cm} \\
    \qquad \qquad (t,m) \in [0,T) \times \pr(\T^d), \vspace{.2cm} \\
    \Phi^{\infty} (T,m) = \phi(m), \quad m \in \pr(\T^d).
    \end{cases}
\end{align}
{The existence and uniqueness of classical solutions to \eqref{phiinfpde} under appropriate technical conditions is by now well-understood, and we refer to \cite{buckdahn2017, CHAUDRUDERAYNAL2022, Tse2021} for some results in this direction. It is also well-known that when $\Phi^{\infty}$ is smooth enough, the functions $\Phi^N$ converge to $\Phi^{\infty}$ with the rate $\frac{1}{N}$, see e.g. \cite{CST, CHAUDRUDERAYNAL2021, DelarueTse}. In the proof of the lower bound in Theorem \ref{thm:mainhs}, we need a similar result, which is stated here.}

\begin{prop} \label{prop:lineqn}
    Suppose that Assumption \ref{assump:section6} holds, and that $\phi = \phi(m) : \pr(\T^d) \to \R$ is Lipschitz and semi-concave with respect to $\| \cdot \|_{-s}$.
    Then we have 
    \begin{align*}
        \Phi^N(t,\bx) \leq \Phi^{\infty}(t,m_{\bx}^N) + \frac{C}{N}, 
    \end{align*}
    for all $(t,\bx) \in [0,T] \times (\T^d)^N$ and a constant $C$ independent of $N$.
\end{prop}

{We emphasize that the estimate appearing in Proposition \ref{prop:lineqn} is by now fairly standard. Compared to existing results, the main novelty is that we do not require $\phi$ to be twice differentiable, we only require a one-sided bound in the form of semi-concavity. This means that $\Phi^{\infty}$ may not in fact be a classical solution to the PDE \eqref{phiinfpde}, and some mollification arguments are required. Lemma \ref{lem:semiconcavebound} below is the key technical result which shows how semi-concavity with respect to $\|\cdot\|_{-s}$ can lead to a one-sided estimate between $\Phi^N$ and $\Phi^{\infty}$.}

Before proving Proposition \ref{prop:lineqn} we need several preliminary results.

\begin{lem} \label{lem:semiconcavebound}
Suppose that Assumption \ref{assump:section6} holds and that $\psi = \psi(q) : H^{-s} \to \R$ is $\cC^{1}$ with $\big[\psi \big]_{\cC^{1,1}(H^{-s})} < \infty$ and $C_S$ semi-concave, i.e. the map
\begin{align*}
    q \mapsto \psi(q) - \frac{C_S}{2} \norm{q}_{-s}^2
\end{align*}
is concave. Then there is a constant $C$ depending on $\psi$ only through $C_S$ such that the inequality
\begin{align} \label{onesidedbound}
   \sum_{j = 1}^d \Big( D_{y_j} \Big[\T^d \ni y=(y_1,...,y_d) \mapsto  (D_m \psi)^j\Bigl(\frac{1}{N} \sum_{k \neq i} \delta_{x^k} + \frac{1}{N} \delta_{y}, x^i\Bigr)\Big](x^i) \Big) \leq \frac{C}{N}
\end{align}
holds at any point $\bx =(x^1,...,x^N) \in (\R^d)^N$ and any index $i \in \{1,...,N\}$ such that the derivatives appearing in \eqref{onesidedbound} exist.
\end{lem}

\begin{proof}
The starting point is to notice that by semi-concavity, we have, for all $q_1,q_2 \in H^{-s}$,
\begin{align} \label{semiconcaveequiv}
    \bigl\langle D_{-s} \psi(q^1) - D_{-s} \psi(q^2), q^1 - q^2 
  \bigr\rangle_{s,-s} \leq C_S \norm{q_1 - q_2}_{-s}^2. 
\end{align}
We fix $\bx = (x^1,...,x^N) \in (\T^d)^N$, $\lambda >0$, $i \in \{1,...,N\}$ and $j \in \{1,...,d\}$ and test the inequality \eqref{semiconcaveequiv} with 
\begin{align*}
    q^1 = m_{\bx}^N, \quad q^2 = m_{\bx}^N + \frac{\lambda}{N} D_{y_j}( \delta_{x^i}),
\end{align*}
to find (using the fact that $s>d/2+1$ and 
 that $H^s$ embeds into $\cC^2(\T^d)$) that
\begin{align*}
    -\Bigl\langle D_{-s} \psi(m_{\bx}^N) - D_{-s} \psi\Bigl(m_{\bx}^N + \frac{\lambda}{N} D_{y_j}( \delta_{x^i} ) \Bigr), \frac{\lambda}{N} D_{y_j}( \delta_{x^i} )
   \Bigr\rangle_{s,-s} \leq C_S \frac{\lambda^2}{N^2} \norm{D_{y_j}(\delta_{x^i})}_{-s}^2 \leq C \frac{\lambda^2}{N^2}, 
\end{align*}
or equivalently 
\begin{align} \label{anothercomp}
    D_{y_j} \big[D_{-s} \psi(m_{\bx}^N)\big](x^i) - D_{y_j} \Big[D_{-s} \psi \Bigl(
    m_{\bx}^N + \frac{\lambda}{N} D_{y_j}( \delta_{x^i} ) 
    \Bigr) \Big](x^i) \leq  \frac{C\lambda}{N}.
\end{align}
Next, we define $\bx^{\lambda} = (x^1,...,x^i - \lambda e^j, x^{i+1},...,x^N)$, and we note that
\begin{align*}
     \norm{m_{\bx}^N + \frac{\lambda}{N} D_{y_j}(\delta_{x^i}) - m_{\bx^{\lambda}}^N}_{-s}  = \frac{1}{N} \sup_{\norm{f}_s = 1} \Big( f(x^i) - f(x^i - \lambda e^j) - \lambda D_{y_j} f(x^i) \Big) \leq \frac{C \lambda^2}{N}, 
\end{align*}
where we used again the Sobolev embedding of $H^s$ into $\cC^2(\T^d)$. Combining this with \eqref{anothercomp}
and Lemma \ref{lem:derivativeconnection}, we have
\begin{align*}
    (D_m \psi)^j\bigl(m_{\bx}^N,x^i\bigr) - (D_m \psi)^j
    \bigl(m_{\bx^{\lambda}}^N, x^i\bigr) &= D_{y_j}\Big[\frac{\delta \psi}{\delta m}(m_{\bx}^N)
    \Big](x^i) - D_{y_j} \Big[
    \frac{\delta \psi}{\delta m}(m_{\bx^{\lambda}}^N)\Big](x^i) \\
    &= D_{y_j} \big[D_{-s} \psi (m_{\bx}^N)\big](x^i) -  D_{y_j}\big[D_{-s} \psi (m_{\bx^{\lambda}}^N)\big](x^i) \\
    &= D_{y_j} \big[D_{-s} \psi (m_{\bx}^N)\big](x^i) 
    - D_{y_j} \Big[D_{-s} \psi\bigl(m_{\bx}^N + \frac{\lambda}{N} D_{y_j}( \delta_{x^i} ) \bigr) \Big](x^i) 
    \\
    &\quad +  D_{y_j} \Big[D_{-s} \psi \bigl(m_{\bx}^N + \frac{\lambda}{N} D_{j}( \delta_{x^i} ) \bigr) \Bigr](x^i) - D_{y_j}\big[D_{-s} \psi (m_{\bx^{\lambda}}^N)\big](x^i) \\
    &\leq \frac{C \lambda }{N} + \frac{C \big[\psi\big]_{\cC^{1,1}(H^{-s})} \lambda^2}{N}.
\end{align*}
 Thus, if $\bx$ is such that the derivatives appearing in \eqref{onesidedbound} exist, then we must have
\begin{align*}
    &D_{y_j} \Big[y \mapsto (D_m \psi)^j
    \Big( \frac{1}{N} \sum_{k \neq i} \delta_{x^k} + \frac{1}{N} \delta_{y}, x^i \Big) \Big](x^i) 
    \\
    &= \lim_{\lambda \to 0^+} \frac{1}{\lambda} \Big((D_m \psi)^j(m_{\bx}^N,x^i) - (D_m \psi)^j(m_{\bx^{\lambda}}^N, x^i)\Big) \leq \frac{C}{N},
\end{align*}
and summing over $j$ gives the result.
\end{proof}

\begin{lem} \label{lem:lineqnwellposed}
If $\phi = \phi(q) : H^{-s} \rightarrow \R$ is $\cC^{1,1}_{\text{loc}}(H^{-s})$ and $C_S$-semi-concave with respect to $\|.\|_{-s}$, and Assumption \ref{assump:section6} is in force, then $\Phi^{\infty}$ is a classical solution of \eqref{phiinfpde}. Moreover, $\Phi^{\infty}$ can be extended to a function
\begin{align*}
    \Phi^{\infty} = \Phi^{\infty}(t,q) : [0,T] \times H^{-s} \to \R.
\end{align*}
such that $\Phi$ is $CC_S$-semi-concave, with $C$ independent of $\phi$.
\end{lem}

We remark that in the statement of Lemma \ref{lem:lineqnwellposed} we use the standard notation $\cC^{1,1}_{\text{loc}}(H^{-s})$ for the set of $\cC^1$ functions $\phi$
such that $\nabla_{-s} \phi$ is Lispchitz on each bounded subset of $H^{-s}$. 

\begin{proof}



 For any $q_0 \in \mathcal{C}^{\infty}(\T^d)$ we can find a unique classical solution to the equation
 \begin{equation*}
     \partial_t q_{t} +\div_{x}(\alpha(t,\cdot) q_{t}) - \Delta_{x} q_{t} = 0 \mbox{ in } (t_0,T) \times \T^d, \hspace{30pt} q(t_0)=q_0,
 \end{equation*}
and moreover, following the proof of Lemma \ref{stabilityLemmaHs27Mars} given in Appendix 
\ref{sec:LinearPDEsappendix}, there is some $C :=C(\sup_{t \in [0,T]} \|\alpha(t,\cdot) \|_{s-1})>0$ such that, for any two solutions $(q^{t_0,q_1}_t)_{t\in [t_0,T]},(q^{t_0,q_2}_t)_{t\in [t_0,T]}$ starting respectively from $q_1$ and $q_2$ it holds 
\begin{equation}
\sup_{t \in [t_0,T]} \|q^{t_0,q_1}_t -q^{t_0,q_2}_t \|_{-s} \leq C \|q_1 - q_2 \|_{-s}.
\label{stabilitygeneralinitialdata29mars}
\end{equation}
As a consequence, we can extend by density the map $q_0 \mapsto (q^{t_0,q_0}_t)_{t\in [t_0,T]} \in \mathcal{C}([t_0,T], H^{-s})$ to the whole of $H^{-s}$ so that \eqref{stabilitygeneralinitialdata29mars} is satisfied for any $q_0 \in H^{-s}$. We define $\Phi^{\infty}$ over the whole of $[0,T] \times H^{-s}(\T^d)$ by
$$ \Phi^{\infty}(t_0,q_0) = \phi(q^{t_0,q_0}_T).$$
The map $\Phi^{\infty}(t_0,\cdot)$ is the composition map $\phi \in \cC^{1,1}_{loc}(H^{-s})$ and the bounded linear map $H^{-s} \ni q_0 \mapsto q_T^{t_0,q_0} \in H^{-s}$. Therefore $\Phi^{\infty}(t_0,\cdot) \in \cC^{1,1}_{loc}(H^{-s})$. 
In particular, the derivatives in the measure argument appearing 
in equation \eqref{phiinfpde} exist.
Moreover, if $\phi$ is $C_S$ semi-concave then $\Phi^{\infty}(t_0,\cdot)$ is $C_SC$ semi-concave for some $C;= C( \sup_{t \in [0,T]} \|  \alpha(t,\cdot) \|_{s-1})>0$. It is then plain to check, by expanding $\Phi^{\infty}$ along solutions to the Fokker-Planck equation that $\Phi^{\infty}(\cdot,q_0)$ is differentiable whenever $q_0 \in \mathcal{P}(\T^d)$ and that $\Phi^{\infty}$ solves \eqref{phiinfpde} over $[0,T] \times \mathcal{P}(\T^d).$

\end{proof}

\begin{proof}[Proof of Proposition \ref{prop:lineqn}]
    We start by defining, for each $\epsilon > 0$, the sup-convolution $\phi^{\epsilon}(q):H^{-s} \to \R$ as in \eqref{eq:sup:convolution:4.1}, i.e. 
    \begin{align*}
        \phi^{\epsilon}(q) = \sup_{m \in \pr(\T^d)} \Big\{\phi(q) - \frac{1}{2\epsilon} \norm{q - m}_{-s}^2 \Big\}. 
    \end{align*}
    By Proposition \ref{prop:supconvproperties}, we know that $\phi^{\epsilon} \in \cC^{1,1}_{\text{loc}}(H^{-s})$. By the proof of the same proposition, it is clear that for all $\epsilon < \frac{1}{2 C_S}$, $\phi^{\epsilon}$ is $2C_S$-semi-concave, where $C_S$ denotes the semi-concavity constant of $\phi$. Thus we can apply Lemma \ref{lem:lineqnwellposed} to see that the map $\Phi^{\infty,\epsilon}$ given by
    \eqref{eq:PhiN:Phiinfty:1}--\eqref{eq:PhiN:Phiinfty:2} with 
    $\phi^\epsilon$ instead of $\phi$
    is a classical solution of \eqref{phiinfpde}(with 
    $\phi^\epsilon$ as boundary condition), and extends to a map $\Phi^{\infty,\epsilon} : [0,T] \times H^{-s} \to \R$ such that for each $t \in [0,T]$, $\Phi^{\infty,\epsilon}(t,\cdot)$ is semi-concave with constant $C_S'$ depending only on $\phi$ only through $C_S$. Moreover, since $\Phi^{\infty,\epsilon} \in \cC^1\big([0,T] \times \pr(\T^d)\big)$, the projections $\Phi^{\infty,\epsilon,N} : [0,T] \times (\T^d)^N \to \R$ given by
    \begin{align*}
        \Phi^{\infty,\epsilon,N}(t,\bx) = \Phi^{\infty,\epsilon}(t,m_{\bx}^N), \quad (t,\bx) \in [0,T] \times \pr(\T^d),
    \end{align*}
    have derivatives 
    \begin{align*}
        \partial_t \Phi^{\infty, \epsilon, N}(t,\bx) = \partial_t \Phi^{\infty, \epsilon}(t,m_{\bx}^N), \quad  D_{x^i} \Phi^{\infty, \epsilon, N}(t,\bx) = \frac{1}{N} D_m \Phi^{\infty, \epsilon}(t,m_{\bx}^N, x^i),
    \end{align*}
    which are continuous and bounded on $[0,T] \times (\T^d)^N$. In addition, we can apply (the reasoning of) Lemma \ref{lem:phinreg} to conclude that
    $\Phi^{\infty, \epsilon, N}$ belongs to $W^{1,2}_{\linf}([0,T] \times {\mathbb T}^d)$, with derivatives $D_{x^i x^j} \Phi^{\infty,\epsilon,N}$ existing almost everywhere. Next, following the proof of Lemma 
     \ref{lem:dxixierror},
     we note that for almost every $(t,\bx)$, the derivatives
    \begin{align}
    D_{y_j} \Bigl[
    y \mapsto  (D_m \Phi^{\infty,\epsilon})^j \Bigl(
    \frac{1}{N} \sum_{k \neq i} \delta_{x^k} + \frac{1}{N} \delta_{y}, x^i
    \Bigr)\Bigr](x^i)
\end{align}
exist for each $j =1,...,d$, and that at any such point we have, for any $i=1,...,N$,
\begin{align*}
    \Delta_{x^i} \Phi^{\infty, \epsilon ,N}(t,\bx) &= \frac{1}{N} \tr\big(D_y D_m \Phi^{\infty, \epsilon}(t,m_{\bx}^N, x^i)\big) 
    \\
&\quad    + \frac{1}{N} \sum_{j = 1}^d \Big( D_{y_j} \Big[
y \mapsto  (D_m \Phi^{\infty,\epsilon})^j
\Bigl(\frac{1}{N} \sum_{k \neq i} \delta_{x^k} + \frac{1}{N} \delta_{y}, x^i\Bigr)\Bigr](x^i) \Big), 
\end{align*}
and so at such a point
\begin{align*}
    - \partial_t& \Phi^{\infty,\epsilon,N}(t,\bx) - \sum_{i=1}^{N} \Delta_{x^i} \Phi^{\infty,\epsilon,N}(t,\bx) -\sum_{i=1}^{N}\alpha(t,x^{i}) \cdot D_{x^{i}} \Phi^{\infty,\epsilon,N}(t,\bx) \\
    &= - \partial_t \Phi^{\infty, \epsilon}(t,m_{\bx}^N) - \frac{1}{N} \sum_{i=1}^{N} \tr\big(D_y D_m \Phi^{\infty, \epsilon}(t,m_{\bx}^N, x^i)\big) -\frac{1}{N} \sum_{i=1}^{N}\alpha(t,x^{i}) \cdot D_m \Phi^{\infty,\epsilon}(t,m_{\bx}^N)   \\
    &\quad - \frac{1}{N} \sum_{i = 1}^N \sum_{j = 1}^d \Big( D_{y_j} 
    \Big[ y \mapsto  (D_m \Phi^{\infty,\epsilon})^j
    \Bigr(\frac{1}{N} \sum_{k \neq i} \delta_{x^k} + \frac{1}{N} \delta_{y}, x^i
    \Bigr) \Big](x^i) \Big) \\
    &= - \frac{1}{N} \sum_{i = 1}^N \sum_{j = 1}^d \Big( D_j 
    \Big[ y \mapsto  (D_m \Phi^{\infty,\epsilon})^j \Bigr(\frac{1}{N} \sum_{k \neq i} \delta_{x^k} + \frac{1}{N} \delta_{y}, x^i\Bigr) \Big](x^i) \Big) \geq -\frac{C}{N}, 
\end{align*}
where $C$ depends on $\phi$ only through $C_S$ and we have used both the PDE for $\Phi^{\infty,\epsilon}$ and Lemma \ref{lem:semiconcavebound}. To summarize, at this stage we know that $\Phi^{\infty, \epsilon, N} \in W^{1,2}_{\infty}([0,T] \times {\mathbb T}^d)$
satisfies 
\begin{align*}
- \partial_t \Phi^{\infty,\epsilon,N}(t,\bx) - \sum_{i=1}^{N} \Delta_{x^i} \Phi^{\infty,\epsilon,N}(t,\bx) -\sum_{i=1}^{N}\alpha(t,x^{i}) \cdot D_{x^{i}} \Phi^{\infty,\epsilon,N}(t,\bx) \geq - \frac{C}N,
\end{align*}
for almost every $(t,\bx) \in [0,T) \times (\T^d)^N$, while it is clear that the function $\Phi^{N,\epsilon}$ defined by
\eqref{eq:PhiN:Phiinfty:1}--\eqref{eq:PhiN:Phiinfty:2}
with $\phi^\epsilon$ instead of $\phi$
satisfies the PDE 
\eqref{phinpde} (with $\bx \mapsto \phi^\epsilon(m^N_\bx)$ as terminal condition)
in a classical sense on $[0,T) \times (\T^d)^N$. This is enough to conclude by the comparison principle (which can be justified in this setting by the It\^o-Krylov formula, as in the proof of Proposition \ref{prop:projections}) that 
\begin{align*}
\Phi^{N,\epsilon}(t,\bx) \leq \Phi^{\infty, \epsilon}(t,m_{\bx}^N) + \frac{C}{N}, 
\end{align*}
for each $(t,\bx) \in [0,T] \times (\T^d)^N$ and a constant $C$ depending on $\phi$ only through $C_S$ (and in partuclar $C$ is independent of $\epsilon$). Finally, it is clear that $\Phi^{N,\epsilon} \to \Phi^N$, $\Phi^{\infty, \epsilon} \to \Phi^{\infty}$ uniformly as $\epsilon \downarrow 0$, and this completes the proof.
\end{proof}

\subsection{The “easy inequality" in the $H^{-s}$-regular case}

 \begin{prop} \label{prop:easyinequality}
    Let Assumption \ref{assump:Hs} hold. Then there is a constant $C$ such that for each $N \in \N$, 
    \begin{align*}
        V^N(t,\bx) \leq U(t,m_{\bx}^N) + \frac{C}{N},
    \end{align*}
    for all $(t,\bx) \in [0,T] \times (\T^d)^N$. 
\end{prop}

\begin{proof}

Fix $(t_0,\bx_0) \in [0,T) \times (\T^d)^N$, and let $\alpha$ be an optimizer for the mean field control problem starting from $(t_0, m_{\bx_0}^N)$. We define 
 $\bX = (X^i)_{i = 1,...,N}$  by 
\begin{align*}
    dX^i_t = \alpha(t,X_t^i) dt + \sqrt{2} dW_t^i, \quad t \in [t_0,T]; \quad X_{t_0}^i = x_0^i, 
\end{align*}
and define $m$ by 
\begin{align*}
    \partial_t m_{t} = \Delta_{x} m_{t} - \text{div}_{x}\bigl(
    m_{t}\alpha(t,\cdot) \bigr), \quad 
    t \in [t_0,T]; \quad m_{t_0} = m_{\bx_0}^N. 
\end{align*}
Notice that by linearity of the Fokker-Planck equation, $m_t = \frac{1}{N} \sum_{i = 1}^N \sL(X_t^i)$, so that 
\begin{align*}
    \E \biggl[ \frac{1}{N} \sum_{i = 1}^N L(X_t^i, \alpha(t,X_t^i))
    \biggr] = \frac{1}{N} \sum_{i=1}^{N} \int_{\T^d} L(x,\alpha(t,x)) \sL(X_t^i)(dx) = \int_{\T^d} L(x,\alpha(t,x)) m_t(dx).
\end{align*}
Using this, we can write
\begin{align*}
    V^N(t_0,\bx_0) &\leq \E\bigg[\int_{t_0}^T \Big( \frac{1}{N} \sum_{i = 1}^N L(X_t^i, \alpha(t,X_t^i) \big) + \sF(m_{\bX_t}^N) \Big)dt + \sG(m_{\bX_T}^N) \bigg] \\
    &= \int_{t_0}^T\int_{\R^d} L(x, \alpha(t,x)) m_t(dx) dt  + \E\bigg[\int_{t_0}^T \sF(m_{\bX_t}^N) dt + \sG(m_{\bX_T}^N) \bigg]  \\
    &= U(t_0,m_{\bx_0}) + \int_{t_0}^T E^1_t dt + E^2,  
\end{align*}
where 
\begin{align*}
    E^1_t := \E[\sF(m_{\bX_t}^N)] - \sF(m_t), \quad E^2 := \E[\sG(m_{\bX_T}^N)] - \sG(m_T). 
\end{align*}
But using the Lipschitz and semi-concavity of $\sF$ and $\sG$ in $H^{-s}$ and applying Proposition \ref{prop:lineqn} (with $\phi = \sF(t,\cdot)$ and $\phi = \sG)$, we see that
\begin{align*}
E^1_t \leq \frac{C}{N}, \quad E^2 \leq \frac{C}{N}, 
\end{align*}
and so 
\begin{align*}
    V^N(t_0,x_0) \leq U(t_0, m_{\bx_0}^N) + \frac{C}{N},
\end{align*}
which completes the proof. 

\end{proof}

\begin{proof}[Proof of Theorem \ref{thm:mainhs}]
    Combine Propositions \ref{prop:easyinequality} and \ref{prop:hardinequality}.
\end{proof}

\section{Proofs for the examples and the convex case}
\label{sec:convexandexproofs}

\begin{lem} \label{lem:mimic}
Suppose that $\sF$ and $\sG$ are convex and $d_1$-Lipschitz, and $H$ satisfies the conditions of Assumption \ref{assump:d1}. Then
\begin{equation}
    U(t,m_{\bx}^N) \leq V^N(t,\bx),
\label{InequalityConvexcase11Avril2023}
\end{equation}
for all $N \in \N$ and $(t,\bx) \in [0,T] \times (\T^d)^N$.
\end{lem}

\begin{proof}
Let us first start with the observation that by Jensen's inequality, we have 
\begin{align} \label{jensens}
    \sF \biggl(\frac{1}{N} \sum_{i = 1}^N \sL(X^i)\biggr) \leq \E\big[ \sF(m_{\bX}^N) \big], \quad 
    \sG\biggl(\frac{1}{N} \sum_{i = 1}^N \sL(X^i)\biggl) \leq \E\big[ \sG(m_{\bX}^N) \big],
\end{align}
for any $N \in \N$ and any $\bX = (X^1,...,X^N)$ a random vector taking values in $(\T^d)^N$.  





Now fix $(t_0,\bx_0) \in [0,T] \times (\T^d)^N$, and let $\bm \alpha = (\alpha^1,...,\alpha^N)$ denote an admissible control for the $N$-particle control problem started from $(t_0,\bx_0)$. We aim to build a control for the mean field control problem out of the control $\bm \alpha$. To do this, we start by using a mimicking result from \cite{Brunick2013} to conclude that we can find for each $i \in \{1,...,N\}$ a bounded and measurable function 
\begin{align*}
    \hat{\alpha}^i = \hat{\alpha}^i(t,x) : [t_0, T] \times \T^d \to \R^d
\end{align*}
such that 
\begin{itemize}
    \item for almost every $t \in [0,T]$, $\hat{\alpha}^i(t,\cdot)$ is a version of $\E[\alpha_t^i | X_t^i = \cdot]$
    \item for each $t \in [0,T]$, we have 
    \begin{align*}
        \sL(X_t^i) = \sL(\hat{X}_t^i), 
    \end{align*}
    where $X^i$ denotes the solution to the SDE 
    \begin{align*}
        d\hat{X}_t^i = \hat{\alpha}^i(t,\hat{X}_t^i) dt + dW_t^i, \quad 
        t \in [t_0,T]; \quad
        \hat{X}_{t_0}^i = x_0^i. 
    \end{align*}
\end{itemize}
Next, we fix an $\sF_{t_0}$-measurable random variable $\xi$ with $\sL(\xi) = m_{\bx_0}^N$. We choose $\xi$ in such a way that we can find a partition of $\Omega$ by $\sF_{t_0}$-measurable sets $A^1,...,A^N$ such that $\bP[A^i] = \frac{1}{N}$ and $\xi = x_0^i$ on $A^i$ for each $i \in \{1,...,N\}$. (For instance, if
we require $\xi$ to take values in $\{x_0^1,...,x_0^N\}$, then 
$A^i$ is just taken as the pre-image of $x_0^i$ by $\xi$.) Now, we build a control for the mean field problem (started from $(t_0,m_{\bx_0}^N)$) as follows. First, we define a function 
\begin{align*}
    \hat{\alpha} = \hat{\alpha}(\omega, t, x) : \Omega \times [t_0,T] \times \T^d \to \R^d
\end{align*}
by the formula 
\begin{align*}
    \hat{\alpha}(\omega, t, x) = \sum_{i = 1}^N 1_{A^i}(\omega) \hat{\alpha}^i(t,x). 
\end{align*}
We denote by $\hat{X}$ the solution of the SDE 
\begin{align} \label{xhatdef}
    d\hat{X}_t = \hat{\alpha}(\omega, t, \hat{X}_t) dt + dW_t, \quad t \in [t_0,T]; \quad \hat{X}_{t_0} = \xi, 
\end{align}
and we define $\hat{\alpha}_t := \hat{\alpha}(\omega, t, X_t)$. It is easy to argue that in fact \eqref{xhatdef} has a unique strong solution thanks to boundedness of the functions $\hat{\alpha}^i$, and it can be built for instance by setting $\hat{X}_t = \sum_{i = 1}^N \hat{X}^i_t 1_{A^i}$, where $\hat{X}_t^i$ solves 
\begin{align} \label{xhatidef}
    d\hat{X}^i_t = \hat{\alpha}^i(t, \hat{X}_t) dt + dW_t, \quad t \in [t_0,T]; \quad  \hat{X}_{t_0} = x_0^i.
\end{align}
We note that we clearly have
\begin{align*}
    \sL(\hat{X}_t) = \frac{1}{N} \sum_{i = 1}^N \sL(\hat{X}_t^i) = \frac{1}{N} \sum_{i = 1}^N \sL(X_t^i). 
\end{align*}
Moreover, we notice that for almost every $t$,
\begin{align} \label{lagrangianest}
    \E\big[L(\hat{X}_t,\hat{\alpha}_t) \big] &= \frac{1}{N} \sum_{i = 1}^N \E \big[L(\hat{X}_t^i, \hat{\alpha}^i(t,\hat{X}^i_t)) \big] 
    = \frac{1}{N} \sum_{i = 1}^N \E \big[L(X_t^i, \hat{\alpha}^i(t,X_t^i)) \big] \nonumber  \\
    &= \frac{1}{N} \sum_{i = 1}^N \E \big[L(X_t^i, \E[\alpha_t^i | X_t^i]) \big] \leq \frac{1}{N} \sum_{i = 1}^N \E\big[ L(X_t^i, \alpha_t^i) \big].
\end{align}
We now conclude by using \eqref{jensens} and \eqref{lagrangianest} to estimate 
\begin{align*}
    U(t_,m_{\bx_0}^N) &\leq \E\bigg[\int_{t_0}^T  \Big( L(\hat{X}_t, \hat{\alpha}_t) + \sF(\sL(\hat{X}_t)) \Big) dt + \sG(\sL(\hat{X}_T)) \bigg] \\
    &\leq \E\bigg[ \int_{t_0}^T \Big(\frac{1}{N} \sum_{i = 1}^N L(X_t^i, \alpha_t^i)  + \sF(\sL(\hat{X}_t)) \Big) dt + \sG(\sL(\hat{X}_{T})) \bigg]  \\
    &= \E\bigg[ \int_{t_0}^T \Big\{ \frac{1}{N} \sum_{i = 1}^N L(X_t^i, \alpha_t^i)  + \sF \Bigl( \frac{1}{N} \sum_{i = 1}^N \sL(X_t^i)) \Bigr) \Big\} dt + \sG \Bigl(\frac{1}{N} \sum_{i = 1}^N \sL(X_{T}^i)\Bigr) \bigg] \\
    &\leq \E\bigg[ \int_{t_0}^T \Big(\frac{1}{N} \sum_{i = 1}^N L(X_t^i, \alpha_t^i)  + \sF(m_{\bX_t}^N) \Big) dt + \sG(m_{\bX_T}^N) \bigg].  
\end{align*}
Taking the infimum over the admissible controls leads to \eqref{InequalityConvexcase11Avril2023}.
\end{proof}

\begin{proof}[Proof of Proposition \ref{prop:convexlip}]
    One inequality is provided by Lemma \ref{lem:mimic}, and the other is given by Proposition \ref{prop:easyinequalityd1}.
\end{proof}

\begin{proof}[Proof of Proposition \ref{prop:ex1}]
Since $\sF$ and $\sG$ depend only on $\bar{m}$ and $L$ is convex, we can conclude that the infimum in the definition of $U$ can be restricted to deterministic controls, i.e. 
\begin{align*}
    U(t,m) = \inf_{\alpha = \alpha(s) : [t,T] \to \R^d}
   \bigg\{ \int_{t}^T \Big(L(\alpha(s)) + \sF(\sL(X_s))  \Big)ds + \sG(\sL(X_T)) \bigg\},
\end{align*}
subject to 
\begin{align*}
    dX_s = \alpha(s) ds + \sqrt{2} dW_s, \quad 
    s \in [t,T]; \quad 
 {\mathcal L}(X_t) = m_0,
\end{align*}
with the infimum taken over square-integrable functions $\alpha$.
Indeed, for any admissible control $\alpha_t$ (in open-loop formulation) we can get a better reward by taking the deterministic control $ t \mapsto \E\left[\alpha_t \right]$.
But if we set $\bar{m}_s = \bar{\sL(X_s)} = \E[X_s]$, then we have the ODE 
\begin{align} \label{ode}
    \frac{d}{dt} \bar{m}_s = \alpha(s), \quad 
    s \in [t,T]; \quad
    \bar{m}_t = \bar{m}, 
\end{align}
and so we can further rewrite $U$ as 
\begin{align*}
    U(t,m) = \inf_{\alpha = \alpha(s) : [t,T] \to \R^d}
   \bigg\{ \int_{t}^T \Big(L(\alpha(s)) + F(\bar{m}_s)  \Big)ds + \sG(\bar{m}_T) \bigg\},
\end{align*}
subject to \eqref{ode},
and it now follows from dynamic programming that $U(t,m) = u(t,\bar{m})$. 

We next argue that $V^N(t,\bx) = v^N(t, \frac{1}{N} \sum_{i = 1}^N x^i)$. For this, we note the convexity of $L$ together with the form of $\sF$ and $\sG$ easily implies that the infimum in the definition of $V^N$ can be restricted to open-loop feedbacks such that $\alpha^i_s = \alpha_s$, for each $i$ and some square-integrable process $\alpha$. That is, we have  
\begin{align*}
    V^N(t,\bx) = \inf_{\alpha \in \sA} \E\bigg[\int_{t}^T \Big(L(\alpha_s) +\sF(m_{\bX_{s}}^N)  \Big)ds + \sG(m_{\bX_T}^N) \bigg] 
\end{align*}
with 
\begin{align*}
    dX_s^i = \alpha_s ds + \sqrt{2} dW_s^i, \quad 
    s \in [t,T]; \quad  X_t^i = x^i. 
\end{align*}
Next, we notice that if we set $\bar{X}_{s} = \frac{1}{N} \sum_{i = 1}^N X^i_{s}$ for $s \in [t,T]$, then we have the SDE 
\begin{align} \label{sdevanishing}
    d\bar{X}_s = \alpha_s ds + \frac{1}{\sqrt{N}} d\bar{W}_{s}^N, \quad X_t = \bar{m_{\bx}^N} := \frac1N \sum_{i=1}^N x_i, 
\end{align}
where $(\bar{W}^N_{s}:= \frac{1}{\sqrt{N}} \sum_{i = 1}^N W^i_{s})_{s \in [t,T]} $ is a Brownian motion. Thus we can write 
\begin{align*} 
    V^N(t,\bx) = \inf_{\alpha} \E\bigg[\int_{t}^T \Big(L(\alpha_s) +F(\bar{X}_s)  \Big)ds + G(\bar{X}_T) \bigg], 
\end{align*}
subject to \eqref{sdevanishing}. 
It again follows from classical representation result that $V^N(t,\bx) = v^N(t,\frac{1}{N} \sum_{i = 1}^N x^i)$. 

\end{proof}

\begin{proof}[Proof of Proposition \ref{prop:ex2}]
On the one hand, it is obvious that 
\begin{align*}
    U(0,\delta_0) = 0, 
\end{align*}
since this is exactly the cost incurred by playing the control $\alpha = 0$. 
Assume now that $T$ is large enough that the density of $\cN_T$ is bounded above by $1$, i.e. $T \geq \frac{1}{2 \pi}$. Using the Cole-Hopf transform and Feynman-Kac formula to solve the PDE \eqref{hjbn} for $V^N$, we have a simple representation formula for $V^N(0, 0)$, which reads 
\begin{align*}
    V^N(0,0) = - \frac{1}{N} \log \int_{(\R^d)^N} e^{- N d_1(m_{\bx}^N, \cN_T)}  \cN_T^{\otimes N}(d\bx) = - \frac{1}{N} \log \E\Big[ e^{-N d_1(m_{\bm \xi}^N, \cN_T)} \Big],
\end{align*}
with $(\xi^i)_{i \in \N}$ being i.d.d. with common law $\cN_T$, and with $m_{\bm \xi}^N = \frac{1}{N} \sum_{i = 1}^N \delta_{\xi^i}$. 
Our goal is to use this formula to prove that for infinitely many $N$, $V^N(0,0) \geq c N^{-1/d}$, with $c > 0$ being independent of $N$. Clearly it suffices to show that there are constants $C$ and $c$ such that
\begin{align} \label{suffest}
    \E\Big[ e^{-N d_1(m_{\bm \xi}^N, \cN_T)} \Big] \leq C e^{-cN^{1 - 1/d}},
\end{align}
for infinitely many $N$. In order to establish \eqref{suffest}, we fix $N \in \N$ such that $N = M^d$ for some $M \in \N$, and we consider a partition of the unit cube $\cC = [0,1]^{d}$ into $N$ cubes $(\cC^{N,i})_{i = 1,...,N}$ of side length $N^{-1/d} = 1/M$. We next consider, for $p \in (0,1)$, the event $A_{p,N}$ defined by 
\begin{align*}
    A_{p,N} := \Big\{ \Big| \{ i : m_{\bm \xi}^N(\cC^{N,i}) > 0\} \Big| > (1-p)N \Big\} = \Big\{ \Big| \{ i : \bm \xi_j \in \cC^{N,i} \text{ for some } 1 \leq j \leq N \} \Big| > (1-p)N \Big\}.
\end{align*}
In other words, $A_{p,N}$ is the even that at least $(1-p) N$ of the $N$ cubes $(\cC^{N,i})_{i = 1,...,N}$ are charged by the random measure $m_{\bm \xi}^N$.
We will now state two lemmas about the event $A_{p,N}$, both of which will be proved below.

\begin{lem} \label{lem:d1lowerbound}
On the event $(A_{p,N})^{\complement}$, we have 
\begin{align*}
    d_1 \bigl(m_{\bm \xi}^N, \cN_T\bigr) \geq c p N^{-1/d}, 
\end{align*}
for some $c > 0$ independent of $p$ and $N$.
\end{lem}

\begin{lem} \label{lem:couponcollector}
For $p > 0$ small enough that $p < \frac{(1-p)^2}{8}$ we have 
\begin{align*}
    \bP[A_{p,N}] \leq C e^{-c(p)N},  
\end{align*}
where $C$ is some constant independent of $N$ and 
$c(p) =  \frac{(1-p)^2}{8} - p$.
\end{lem}

With these lemmas in hand, we can easily conclude the proof, since for $p > 0$ small enough we have
\begin{align*}
    \E \Bigl[e^{-Nd_1(m_{\bm \xi}^N, \cN_T)}\Bigr] &\leq \bP[A_{p,N}] + \bigl(1 - \bP[A_{p,N}]\bigr) e^{- c N^{1-1/d}} \\
    &\leq C e^{-cN} + e^{-c N^{1 - 1/d}} 
    \leq C e^{-c N^{1 - 1/d}},
\end{align*}
for some $c$, $C$ independent of $N$, which clearly implies \eqref{suffest}.
\end{proof}

We now prove the two lemmas used in the above example.

\begin{proof}[Proof of Lemma \ref{lem:d1lowerbound}]
    Fix $\omega \in (A_{p,N})^{\complement}$. Define a set $I \subset \{1,...,N\}$ by 
    \begin{align*}
        I := \Big\{i \in \{1,...,N\} : m_{\bm \xi(\omega)}(\cC^{N,i}) = 0 \Big\},
    \end{align*}
    and then set $\cC^I = \cup_{i \in I} \cC^{N,i}$. 
    Then by the definition of $A_{p,N}$, we must have $|I| \geq pN$. Now let $\gamma$ be any coupling between $m_{\bm \xi(\omega)}^N$ and $\cN_T$. That is,  $\gamma \in \pr(\R^d \times \R^d)$, and we have 
    \begin{align*}
        \gamma( A \times \R^d) = m_{\bm \xi(\omega)}(A), \quad \gamma(\R^d \times A) = \cN_T(A),
    \end{align*}
    for any Borel set $A \subset \R^d$. Finally, for $i = 1,...,N$, let $\cC^{N,i}_{1/2}$ denote the cube with the same center and half the side-length of $\cC^{N,i}$. Then we have 
    \begin{align*}
        \int_{\R^d \times \R^d} |x-y| d \gamma(x,y) &= \int_{\R^d \times \R^d} |x-y| 1_{x \in (\cC^I)^{\complement} }d\gamma(x,y)  \\
        &\geq \int_{\R^d \times \R^d} \Big(  |x-y| 1_{x \in (\cC^I)^{\complement}} \sum_{j \in I} 1_{y \in \cC^{N,j}_{1/2}} \Big) d\gamma(x,y) \\
        &\geq \int_{\R^d \times \R^d} \Big(c N^{-1/d} 1_{x \in (\cC^I)^{\complement}} \sum_{j \in I} 1_{y \in \cC^j_{1/2}} \Big) d\gamma(x,y) \\
        &= cN^{-1/d} \sum_{j \in I} \cN_T(\cC^j_{1/2}) 
         \geq c p N^{-1/d}, 
    \end{align*}
    where in the last step we used the fact that the density of $\cN_T$ is bounded 
    from below on $\cC$.
    Taking an infimum over the couplings $\gamma$ shows that 
    \begin{align*}
        d_1(m_{\bm \xi}(\omega), \cN_T) \geq cpN^{-1/d},
    \end{align*}
    for each $\omega \in (A_{p,N})^c$, which completes the proof.
\end{proof}

\begin{proof}[Proof of Lemma \ref{lem:couponcollector}]
    First, we make a simplifying observation. Let $(U^i)_{i \in \N}$ denote a sequence of i.i.d. random variables each of which is uniformly distributed on the unit cube $\cC$, i.e. $\sL(U^i)$ is the Lebesgue measure restricted to $\cC$.  
    For ${\bm U}:=(U^1,...,U^N)$,
    let
    $m_{\bm U}^N = \frac{1}{N} \sum_{i = 1}^N \delta_{U^i}$, and define 
    \begin{align*}
    B_{p,N} = \Big\{ \Big| \{ i : m_{\bm U}^N(\cC^{N,i}) > 0\} \Big| > (1-p)N \Big\} = \Big\{ \Big| \{ i :  U_j \in \cC^{N,i} \text{ for some } 1 \leq j \leq N \} \Big| > (1-p)N \Big\}.
\end{align*}
In other words, $B_{p,N}$ is defined exactly as $A_{p,N}$, but with the i.i.d. uniform random variables $(U^i)_{i = 1,...,N}$ replacing the i.i.d. Gaussian random variables $(\xi^i)_{i = 1,...,N}$. Thanks to the fact that we assumed $T$ is large enough that the density of $\cN_T$ is bounded above by $1$, it is clear that $\bP[A_{p,N}] \leq \bP[B_{p,N}]$,
and so it suffices to prove the estimate 
\begin{align} \label{bpnest}
    \bP[B_{p,N}] \leq C e^{-cN},
\end{align}
for $p$ small enough and $C, c$ independent of $N$. We mention that $\bP[B_{p,N}]$ is connected to the well-known ``coupon collector problem" - in the language of this problem, $\bP[B_{p,N}]$ is the probability that the coupon collector collects at least $(1-p)N$ distinct coupons in the first $N$ draws from $N$ distinct coupons. 

To analyze the probability $\bP[B_{p,N}]$, we first define an increasing sequence of random times $(T_m)_{m \in \N}$ by setting $T_0 :=0$ and  
\begin{align*}
    T_m := \inf \Big\{ j \in \N : \Big|\{i : m_{\bm U}^j(\cC^{N,i}) > 0 \}  \Big| = m  \Big\}, \quad m \in \N,
\end{align*}
and then we set $\tau_m := T_m - T_{m-1}$ for $m \in \N\setminus \{0\}$.
So in the language of the coupon collector problem, $T_m$ is the number of draws needed to find the $m^{th}$ distinct coupon, and $\tau_m$ is the number of draws between the discovery of the $m^{th}$ distinct coupon and that of the $(m-1)^{th}$ distinct coupon. It is easy to check that (using the notation $\lceil \cdot \rceil$ for the ceil part)
\begin{align*}
    \bP[B_{p,N}] \leq  \bP[\tau_1 + ... + \tau_{\lceil(1-p) N \rceil} \leq  N], 
\end{align*}
and that $(\tau_i)_{i \in \N}$ are independent, with $\tau_i \sim \text{Geom}(1 - \frac{i-1}{N})$
i.e. for $k \in \N \setminus \{0\}$, 
\begin{align*}
\bP[\tau^i = k] = \Bigl(\frac{i-1}{N} \Bigr)^{k-1} 
\Bigl(1 - \frac{i-1}{N} \Bigr).
\end{align*}
Using the well-known formula for the moment generating function of a geometric random variable, we find that
\begin{align*}
  \E[e^{- \lambda \tau_i}] =  \frac{e^{ - \lambda} \Big(1 - \frac{i-1}{N} \Big)}{1 - e^{- \lambda } \big(\frac{i-1}{N}\big)},  
\end{align*}
for each $\lambda > 0$. Thus, we see that for $\lambda > 0$, we have
\begin{align} \label{mgfcomp}
    \bP[\tau_1 + ... &+ \tau_{\lceil(1-p) N \rceil} \leq N] = \bP\Big[e^{- \lambda (\tau_1 + ... + \tau_{\lceil(1-p) N \rceil})} \geq e^{-\lambda N}\Big] \nonumber  \\
    &\leq e^{\lambda N} \E\Big[ e^{- \lambda (\tau_1 + ... + \tau_{\lceil(1-p) N \rceil})}  \Big] \nonumber  \\
&= e^{\lambda N} \prod_{i = 1}^{\lceil(1-p) N \rceil} \E\Big[e^{- \lambda \tau_i} \Big]   \nonumber \\
&=
e^{\lambda N} \prod_{i = 1}^{\lceil(1-p) N \rceil} \frac{e^{ - \lambda} \Big(1 - \frac{i-1}{N} \Big)}{1 - e^{- \lambda } \big(\frac{i-1}{N}\big)} \nonumber \\
& = \exp\Big\{ \lfloor pN \rfloor \lambda  + \sum_{i = 1}^{\lceil(1-p)N \rceil} \Big( \ln \big(1 - \frac{i-1}{N}\big) - \ln \big(1 - e^{-\lambda} (\frac{i-1}{N}) \big) \Big) \Big \}. 
\end{align}
To go further, we note that 
\begin{align} \label{logcomp}
\log \Bigl(1 - \frac{i-1}{N}\Bigr) - \log
\Bigl(1 - e^{-\lambda}\frac{i-1}{N}\Bigr) &= \sum_{n = 1}^{\infty} n^{-1} e^{- \lambda n}\big(\frac{i - 1}{N} \big)^n - \sum_{n = 1}^{\infty} n^{-1} \big(\frac{i - 1}{N} \big)^n 
\\
&\leq (e^{- \lambda} - 1) \big(\frac{i - 1}{n}\big).\nonumber
\end{align}
Notice that for any $0< \lambda \leq 1$, we have $\exp (- \lambda) - 1 \leq - \frac{\lambda}{4}$, and so combining \eqref{mgfcomp} with \eqref{logcomp} shows that for $0 < \lambda \leq 1$,  
\begin{align*}
 \bP[\tau_1 + ... &+ \tau_{\lceil(1-p) N \rceil} \leq N]  \leq 
\exp\Big\{ \lfloor pN \rfloor \lambda  - (1 - e^{- \lambda}) \sum_{i = 1}^{\lceil(1-p)N \rceil} \frac{i - 1}{N}  \Big \} \\
&\leq \exp \Big\{p N \lambda - (1 - e^{- \lambda}) (1-p)^2 \frac{N}{2}  \Big\} \\
&\leq \exp\Big\{ p N \lambda - \frac{1}{8} \lambda N (1 - p)^2 \Big\} = \exp\Bigl\{ \lambda N  \Bigl(p - \frac{(1-p)^2}{8} \Bigr) \Bigr\}.
\end{align*}
Fixing $\lambda = 1${,} we see that if $p$ small enough that $p < \frac{(1-p)^2}{8}$, then we have
\begin{align*}
    \bP[B_{p,N}] \leq \exp(- c N). 
\end{align*}
with $c = \big(p - \frac{(1-p)^2}{8} \big)$.
This concludes the proof.
\end{proof}

\appendix

\section{Some auxiliary estimates for finite-dimensional PDEs}

In this appendix, we first give some results about linear and semi-linear parabolic equations that were used in Section \ref{sec: PropertiesoftheValuefunction}. 
We also give a sketch of proof of Lemma \ref{lem:regularisation} in 
Subsection \ref{subse:mollif:arg}.

The following form of Grönwall's Lemma will be useful (see \cite{YeGaoDing}):

\begin{lem}
    Assume that $l \in L^{\infty}([0,T])$ satisfies, for some non-decreasing measurable $C_1 :[0,T] \rightarrow \R^+ $, some $C_2 >0$ and all $t \in [0,T]$,
    $$l(t) \leq C_1(t) + C_2 \int_t^T \frac{l(s)}{\sqrt{s-t}}ds,$$
then, there exists a constant $C(C_2,T)$ (depending on $C_2$ and $T$) such that, for all $t \in [0,T]$, 
$$l(t) \leq C_1(t) C(C_2,T).$$
\label{GronwallLemma}
\end{lem}

\subsection{Linear PDEs}

\label{sec:LinearPDEsappendix}

\begin{lem}
Suppose that $\alpha:[t_0,T] \times \T^d \rightarrow \R^d$ is Lipschitz continuous in $x$ uniformly in $t$, take $g \in \mathcal{C}^1(\T^d)$, $t_1 \in (t_0,T]$ and let $v$ be the solution to
\begin{equation} 
-\partial_t v(t,x) - \Delta_{x} v(t,x) -\alpha(t,x) \cdot D_{x}v(t,x) =0 \mbox{ in } (t_0,t_1) \times \T^d, \hspace{15pt} v(t_1,x)=g(x) \mbox{ in }\T^d.
\label{equationbacward20mars2023}
\end{equation} 
Then, 
\begin{itemize}
    \item If $\alpha$ satisfies
\begin{equation}
\sup_{t_0 \leq t  \leq T} \|\alpha(t,\cdot) \|_{\linf} +\sup_{t_0 \leq t \leq T} \sqrt{T-t} \|D_{x} \alpha(t,\cdot) \|_{\linf} \leq C_1,
\label{regularityforalpha27avril2023}
\end{equation}
for some $C_1>0$. Then, there is $C_1'$ depending only on $C_1$ such that
\begin{align}
  \sup_{t \in [t_0,t_1]} \|D_{x}v(t,\cdot) \|_{\linf} +\sup_{t \in [t_0,t_1]} \sqrt{t_1-t} \|D^2_{xx}v(t,\cdot) \|_{\linf} \leq C_1' \|D_{x}g \|_{\linf}.
\label{Mixedestimated1reg28mars2023}
\end{align}
     \item If we further assume that $g \in H^s$ and $\sup_{t_0 \leq t \leq T} \| \alpha (t,\cdot) \|_{s-1} \leq C_2$ for some $C_2>0$, then

\begin{equation}
    \sup_{t \in [t_0,t_1] } \|v(t,\cdot ) \|_{H^s} \leq C_2' \|g \|_{H^s},
\label{estimateHs28mars2023}
\end{equation}
for some some $C_2'$ depending only on $C_2$.
\end{itemize}
\label{lem:LinearBackward28Mars}
\end{lem}

\begin{proof}
 We start with the proof of \eqref{Mixedestimated1reg28mars2023} Let us write $P_t$ for the heat semi-group on $\T^d$, i.e. given $\phi : \T^d \to \R$,
    \begin{align*}
    \big(P_t \phi \big)(x) = \psi(t,x), \quad \text{where $\psi$ solves }  \partial_t \psi = \Delta_{x} \psi, \,\, \psi(0,x) = \phi(x). 
    \end{align*}
We note for later use that we have the classical estimates for $P$:
    \begin{align} \label{smoothingSchauder}
        \norm{P_t f}_{\linf} \leq \norm{f}_{\linf}, \quad 
        \norm{D_{x}P_t f}_{\linf} \leq \frac{C}{\sqrt{t}} \norm{f}_{\linf}, 
    \end{align}
    \begin{align} \label{smoothing}
        \norm{P_t f}_s \leq \norm{f}_s, \quad 
        \norm{P_t f}_{s} \leq \frac{C}{\sqrt{t}} \norm{f}_{s-1}. 
    \end{align}
Notice that the solution $v$ to \eqref{equationbacward20mars2023} must satisfy  for $(t,x) \in [t_0,t_1] \times \T^d$
\begin{equation}
    v(t,x) = P_{t_1-t}g(x) +\int_t^{t_1} P_{s-t} \left[\alpha(s,\cdot)\cdot D_{x}v(s,\cdot)\right] (x) ds.
    \label{equationforv21Mars2023}
\end{equation}
And therefore, differentiating \eqref{equationforv21Mars2023} in $x$, it holds that
$$D_{x}v(t,x) = P_{t_1-t}D_{x}g(x) + \int_t^{t_1} D_{x}P_{s-t} \left[ \alpha(s,\cdot) \cdot D_{x}v(s,\cdot) \right](x)ds,$$
which, in turn, using the smoothing properties \eqref{smoothingSchauder}, leads to
\begin{align*}
    \|D_{x}v(t,\cdot) \|_{\linf} \leq \|D_{x}g \|_{\linf} + \sup_{t \in [t_0,t_1]}\|\alpha (t,\cdot) \|_{\linf} \int_t^{t_1} \frac{\|D_{x}v(s,\cdot) \|_{\linf}}{\sqrt{s-t}}ds,
\end{align*}
and we conclude with Lemma \ref{GronwallLemma} that inequality
\begin{equation}
    \sup_{t\in [t_0,t_1]} \|D_{x}v(t,\cdot) \|_{\linf} \leq C \|D_x g \|_{\linf}
\label{Lipestimateforv27avril2023}
\end{equation}
holds for some $C=C(C_1,T)$.
For the estimate on $\|D_{xx}^2v(t,\cdot) \|_{\infty}$, we differentiate \eqref{equationforv21Mars2023} as follows
\begin{equation*}
D_{xx}^2v(t,x) = D_{x}P_{t_1-t}D_{x}g (x) + \int_t^{t_1} D_{x}P_{s-t} \left[D_{x}\alpha(s,\cdot)D_{x}v(s,\cdot) + D_{xx}^2v(s,\cdot)\alpha(s,\cdot) \right](x) ds.
\end{equation*}
And therefore, using the smoothing properties \eqref{smoothingSchauder} and the regularity condition \eqref{regularityforalpha27avril2023} on $\alpha$ we have
\begin{align*}
\|D^2_{xx}v(t,\cdot) \|_{\linf} &\leq \frac{\|D_{x}g\|_{\linf}}{\sqrt{t_1-t}} + \sup_{t\in [t_0,t_1]}\|D_{x}v(t,\cdot) \|_{\linf} \int_t^{t_1} \frac{\|D_{x}\alpha(s,\cdot) \|_{\linf}}{\sqrt{s-t}} ds \\
&+ \sup_{t\in [t_0,t_1]} \|\alpha(t,\cdot) \|_{\linf} \int_t^{t_1} \frac{\|D_{xx}^2v(s,\cdot)\|_{\linf}}{\sqrt{t-s}}ds \\
&\leq  \frac{\|D_{x}g\|_{\linf}}{\sqrt{t_1-t}} + C_1 \sup_{t\in [t_0,t_1]} \|D_{x}v(t,.) \|_{\linf} \int_t^{t_1} \frac{1}{\sqrt{T-s}\sqrt{s-t}} ds \\
&+ C_1 \int_t^{t_1} \frac{\|D_{xx}^2v(s,\cdot)\|_{\linf}}{\sqrt{t-s}}ds. 
\end{align*} 
Noticing that 
$$\int_t^{t_1} \frac{1}{\sqrt{T-s}\sqrt{s-t}} ds \leq \int_t^T \frac{1}{\sqrt{T-s}\sqrt{s-t}} ds =\pi,$$
we get, applying Gronwall's Lemma \ref{GronwallLemma},
$$ \|D^2_{xx}v(t,\cdot) \|_{\linf} \leq C \Bigl( \sup_{t\in [t_0,t_1]}\|D_{x}v(t,\cdot)\|_{\linf} + \frac{\|D_{x}g\|_{\linf}}{\sqrt{t_1-t}} \Bigr),$$
for some $C=C(C_1,T)$. We conclude using the bound \eqref{Lipestimateforv27avril2023} on $D_xv$. 
\newline \newline
We continue with the proof of \eqref{estimateHs28mars2023} under the corresponding assumption. Being $s-1>d/2$, there is some $C_s$ depending only on $s$ such that $\|fg \|_{H^{s-1}} \leq C_s \|f\|_{H^{s-1}}  \|g \|_{H^{s-1}}$ whenever $f,g \in H^{s-1}$. Therefore, using the smoothness assumed on $\alpha$, we have for each fixed $r \in [t_0,t_1]$ 
    \begin{align*}
        \norm{\alpha(r,\cdot) \cdot D_{x}v(r,\cdot)}_{H^{s-1}} \leq C \norm{D_{x}v(r,.)}_{H^{s-1}} \leq C \norm{v}_{H^s}.
    \end{align*}
Combining this with the smoothing property \eqref{smoothing}, we have 
    \begin{align}
    \norm{P_{r-t} \big(\alpha(r,\cdot) \cdot D_{x}v(r,\cdot)\big)}_{H^s} \leq \frac{C}{\sqrt{r-t}} \norm{v(r,\cdot)}_{H^s}.
    \end{align}
But now using \eqref{equationforv21Mars2023}, we have, for all $t\in [t_0,t_1]$, 
    \begin{align*}
        \norm{v(t,\cdot)}_{H^s} &\leq \norm{g}_{H^s} + \int_t^{t_1} \norm{P_{r-t} \big(b(r,\cdot) \cdot D_{x}v(r,\cdot)\big)}_{H^s} dr \\
        &\leq \norm{g}_{H^s} + C  \int_t^{t_1} \frac{\|v(r,\cdot)\|_{H^s}}{\sqrt{r-t}}dr.
    \end{align*}
Using Grönwall Lemma \ref{GronwallLemma}, we get \eqref{estimateHs28mars2023}.

\end{proof}







Now we can prove the stability Lemmas 
\ref{stabilityLemmaHs27Mars} and \ref{lem:stabilityEstimate20Mars2023}.

\begin{proof}[Proof of Lemma \ref{stabilityLemmaHs27Mars}]
We argue by duality using Lemma 
\ref{lem:LinearBackward28Mars} above. Indeed, for any $m^1, m^2$ satisfying 
\begin{align*}
           \partial_t m^i_{t} = \Delta m^i_{t} - \div \bigl(m^i_{t} \alpha (t,\cdot) \bigr),
    \end{align*}
for $(t,x) \in [t_0, T] \times \T^d$, any $g \in H^s$ and any $t_1 \in [t_0,T]$, we have 
    \begin{align*}
    \int_{\T^d} g(x) (m^1_{t_1} - m^2_{t_1})(dx) = \int_{\T^d} v(t_0, x) (m_{t_0}^1 - m_{t_0}^2)(dx),
    \end{align*}
where $v$ solves \eqref{equationbacward20mars2023}. As a consequence, 
    \begin{align*}
    \int_{\T^d} g(x) (m^1_{t_1} - m^2_{t_1})(dx) \leq \norm{v(t_0,\cdot )}_{H^s} \norm{m_{t_0}^1 - m_{t_0}^2}_{H^{-s}} \leq C \norm{g}_{H^s} \norm{m_{t_0}^1 - m_{t_0}^2}_{H^{-s}}, 
    \end{align*}
    and so taking a supremum over $g$ with $\norm{g}_{H^s} \leq 1$ gives the result.  
\end{proof}

\begin{proof}[Proof of Lemma \ref{lem:stabilityEstimate20Mars2023}]
We argue once again by duality to deduce that, for all $t_1 \in [t_0,T]$ and all $1$-Lipschitz function $g$,
\begin{align*}
    \int_{\T^d} g(x)(m^1_{t_1}-m^2_{t_1})(dx) &= \int_{\mathbb T^d} v(t_0,x)(m^1_{t_0}-m^2_{t_0})(dx) \\
    &\leq \|D_xv(t_0,.) \|_{\linf} d_1(m^1_{t_0},m^2_{t_0})
    \leq C_1' d_1(m^1_{t_0},m^2_{t_0}),
\end{align*}
where $v$ is solution to \eqref{equationbacward20mars2023} with terminal data $g$ and $C_1'$ is the constant appearing in \eqref{Mixedestimated1reg28mars2023}. Taking the supremum over $1$-Lipschitz functions $g$ gives the first statement in Lemma \ref{lem:stabilityEstimate20Mars2023}. For the second statement we proceed similarly. Once again, let $t_1 \in (t_0,T]$, let $g$ be a $1$-Lipschitz function and $v$ the solution to \eqref{equationbacward20mars2023} with terminal data $g$. Using \eqref{Mixedestimated1reg28mars2023} we infer that
\begin{align*}
\int_{\T^d} g(x) \bigl( m_{t_1}^1 - m_{t_1}^2 \bigr)(dx) &= \int_{\mathbb T^d} v(t_0,x)(m^1_{t_0}-m^2_{t_0})(dx) \\
    &\leq \Bigl( \|D_xv(t_0,.) \|_{\linf} +  \|D^2_{xx}v(t_0,.) \|_{\linf} \Bigr) \norm{m^1_{t_0}
-m^2_{t_0} }_{-2,\infty} \\
&\leq \frac{C_1'}{\sqrt{t_1-t_0}} \norm{m^1_{t_0}
-m^2_{t_0} }_{-2,\infty}.
\end{align*}
Taking the supremum over $1$-Lipschitz functions $g$ gives the result.
\end{proof}


\subsection{HJB equation}

\begin{lem}
    Assume that $H$ satisfies (1) and (2) in Assumption \eqref{assump:d1}. Take $f\in L^{\infty}([t_0,T],\mathcal{C}^2(\T^d))$ and $g \in \mathcal{C}^2(\T^d))$. Let $u \in C^{1,2}([t_0,T) \times \T^d)$ be the solution to 
\begin{equation}
    -\partial_t u(t,x) + H\bigl(x,
    D_{x} u(t,x)\bigr) - \Delta_{x} u(t,x) = f(t,x)  \text{ in } (t_0,T)\times \T^d, \quad 
    u(T,x) = g(x) \text{ in } \T^d.
\end{equation}
Then $ \|D_{x}u \|_{\linf} + \sup_{t \in [t_0,T]} \sqrt{T-t} \|D_{xx}^2 u(t,.) \|_{\linf} \leq C$
for some $C= C(\|D_{x}f \|_{\linf}, \|D_{x}g \|_{\linf})>0.$
If we further assume that $H \in \cC^s(\T^d \times \R^d)$, $f \in L^{\infty}([0,T], \mathcal{C}^s(\T^d))$ and $g \in \mathcal{C}^s(\T^d)$ then
$$\sup_{t \in [t_0,T]} \|u(t,\cdot) \|_{\cC^s} + \sup_{t \in [t_0,T]} \sqrt{T-t}\|u(t,\cdot) \|_{\cC^{s+1}} \leq C'$$
for some $C'=C'( \sup_{t\in [t_0,T]} \|f(t,\cdot) \|_{\cC^s}, \|g \|_{\cC^s})$.
\end{lem}

\begin{proof}
    The Lipschitz estimate is standard and follows from the classical Bernstein method and the condition on the growth of $D_xH$ in Assumption \eqref{assump:d1}. Therefore, there is $C = C(\|D_{x}f \|_{\linf}, \|D_{x}g \|_{\linf})$ such that 
\begin{equation}
    \|D_{x}u \|_{\linf} \leq C.
    \label{Lipschitzestimate27Mars2023}
\end{equation}
Being a classical solution, $u$ must satisfy, for all $(t,x) \in [t_0,T] \times \T^d$,

    $$u(t,x) = P_{T-t} g(x) + \int_t^T P_{s-t} f(s,\cdot)(x)ds - \int_t^T P_{s-t} \left[ H(\cdot,D_{x}u(s,\cdot) \right](x)ds. $$
Differentiating twice and integrating by parts leads, for all $(t,x) \in [t_0,T] \times \T^d$ to
\begin{equation}
    D^{2}_{xx}u(t,x) = D_{x}P_{T-t} D_x g (x) + \int_t^T D_{x}P_{s-t}D_{x}f(s,\cdot)(x)ds - \int_t^T D_{x}P_{s-t}D_x \left[ H(\cdot,D_{x}u(s,\cdot) \right](x)ds.
\end{equation}
Using the smoothing properties \eqref{smoothingSchauder}, we deduce that inequality 
$$ \|D^2_{xx}u(t,\cdot) \|_{\linf} \leq \frac{\|D_{x}g\|_{\linf}}{\sqrt{T-t}} + \int_t^T \frac{\|D_{x}f \|_{\linf}}{\sqrt{s-t}}ds + \int_t^T \frac{C(1 + \|D^2_{xx}u(s,\cdot) \|_{\linf})}{\sqrt{s-t}}ds $$
holds for some $C=C(\|D_{x}u\|_{\linf})>0$ and for all $t \in [t_0,T)$. As a consequence, taking \eqref{Lipschitzestimate27Mars2023} into account  we can find $C=C( \|D_{x}f \|_{\linf}, \|D_{x}g\|_{\linf}) >0$ such that inequality
$$\|D^2_{xx}u(t,\cdot)\|_{\linf} \leq C \Bigl( 1 + \frac{1}{\sqrt{T-t}} + \int_t^T \frac{\|D^2_{xx}u(s,\cdot) \|_{\linf}}{\sqrt{s-t}}ds \Bigr) $$
holds for all $t \in [t_0,T)$ and we deduce from Grönwall's Lemma \eqref{GronwallLemma} that $$\|D^2_{xx}u(t,\cdot) \|_{\linf} \leq C\Bigl(1+ \frac{1}{\sqrt{T-t}}\Bigr), $$
for some $C = C(\|D_{x}f \|_{\linf}, \|D_{x}g \|_{\linf})>0$ .

For the second part of the lemma, we easily prove by induction that
$$\sup_{t \in [t_0,T]} \|u(t,.) \|_{\cC^s} \leq C,$$
for some $C=C(\sup_{t \in [t_0,T]} \|f(t,\cdot)\|_{\cC^s}, \|g\|_{\cC^s})>0$ and we deduce the bound on $\|u(t,\cdot) \|_{\cC^{s+1}}$ as we got the bound on $\|D^2_{xx}u(t,\cdot) \|_{L^{\infty}}$.
\label{lem:linearHJBappendix3avril2023}

\end{proof}

\section{A mollification argument}
\label{subse:mollif:arg}

We provide the proof of Lemma
\ref{lem:regularisation}.

\begin{proof}

\underline{Step 1.}
We introduce the following notation. For any $n \geq 1$, we denote by 
${\mathcal P}_n({\mathbb T}^d)$ the collection of probability measures whose Fourier coefficients of order greater than $n$ are null, i.e. 
\begin{equation*} 
{\mathcal P}_n({\mathbb T}^d) := \bigl\{ m \in {\mathcal P}({\mathbb T}^d) : \widehat{m}^{\boldsymbol k} = 0, \quad \vert {\boldsymbol k} \vert_\infty \geq n \bigr\}, 
\end{equation*}
where ${\boldsymbol k}$ in the above notation is implicitly understood as a multi-index ${\boldsymbol k}=(k_1,\cdots,k_d) \in {\mathbb Z}^d$
of sup norm $\vert {\boldsymbol k} \vert_\infty := \max(\vert k_1\vert, \dots , \vert k_d\vert)$. 
Accordingly, the set ${\mathcal P}_n({\mathbb T}^d)$ 
can be identified 
with a finite-dimensional set ${\mathcal O}_n$ 
describing the Fourier coefficients of order less than $n$. 
The set ${\mathcal O}_n$ can be described in an exhaustive manner by means of Bochner's theorem.
We just refer to \cite{cecchin2022weak} for the details as the precise formulation of ${\mathcal O}_n$ does not really matter for our needs. 
The only two things that we need are 
\begin{enumerate}
\item The set ${\mathcal O}_n$ can be regarded as a subset of ${\mathbb R}^{\vert D(n) \vert} \times {\mathbb R}^{\vert D(n) \vert}$ and 
it contains an open ball ${\mathcal B}_n \times {\mathcal B}_n$ centered at the origin of 
${\mathbb R}^{\vert D(n)\vert} \times {\mathbb R}^{\vert D(n)\vert}$, where 
$D(n)$ is a certain subset of 
$\{ {\boldsymbol k} \in {\mathbb Z}^d \setminus \{0\} : \vert {\boldsymbol k} \vert_\infty < n\}$ 
and 
$\vert D(n) \vert$ denotes its cardinal.
\item
Denoting by
\begin{equation*} 
\begin{split}
{\mathscr I}_n &: 
\Bigl( (a_{\boldsymbol k})_{{\boldsymbol k} \in D(n)}, 
(b_{\boldsymbol k})_{{\boldsymbol k} \in  D(n)}\Bigr) 
\in 
{\mathbb R}^{D(n)} \times {\mathbb R}^{D(n)}
\\
&\mapsto 
\Bigl( x \in {\mathbb T}^d \mapsto 
 1+ 2 \sum_{ {\boldsymbol k} \in D(n)}
 \bigl( a_{\boldsymbol k} \cos ( 2 \pi {\boldsymbol k} \cdot x) +  b_{\boldsymbol k} \sin ( 2 \pi {\boldsymbol k} \cdot x) \bigr)
 \Bigr),
 \end{split}
\end{equation*} 
the mapping ${\mathscr I}_n$ 
is one-to-one from 
 ${\mathcal O}_n$ onto ${\mathcal P}_n({\mathbb T}^d)$ (here, the additional $1$ is regarded as the constant function on the torus, equal to $1$).  For $m \in {\mathcal P}_n({\mathbb T}^d)$, its pre-image
 by ${\mathscr I}_n$ writes
 \begin{equation*}
 {\mathscr I}_n ^{-1}(m) = \Bigl( \Re\bigl[\widehat{m}^{\boldsymbol k}\bigr],\Im\bigl[ \widehat{m}^{\boldsymbol k}\bigr] \Bigr)_{ {\boldsymbol k} \in 
 D(n)},
  \end{equation*} 
\end{enumerate}  
where $\Re(z)$ and $\Im(z)$ denote the real and imaginary parts of a complex number $z$. 

The next ingredient that is needed in the proof is the Féjer kernel of rank $n$, which we denote by
$\varphi_n : {\mathbb T}^d \rightarrow {\mathbb R}$. We recall that $\varphi_n$ is a (symmetric) probability density whose Fourier coefficients are given by 
\begin{align*}
&\widehat{\varphi}^{\boldsymbol k} = \prod_{j=1}^d \bigl( 1 - \frac{k_j}{n} \bigr), &\text{if} \ \vert {\boldsymbol k} \vert_{{\infty}} < n,
\\
&\widehat{\varphi}^{\boldsymbol k} =0, &\text{if} \ \vert {\boldsymbol k} \vert_{\infty} \geq n. 
\end{align*}
for any multi-index ${\boldsymbol k}=(k_1,\cdots,k_d) \in {\mathbb Z}^d$. In particular, 
for any $m \in {\mathcal P}({\mathbb T}^d)$, 
the convolution $m*\varphi_n$ is an element of ${\mathcal P}_n({\mathbb T}^d)$. 
Also, for $\eta \in [0,1]$ and for 
$({\boldsymbol a},{\boldsymbol b}):=( (a_{\boldsymbol k})_{{\boldsymbol k} \in D(n)}, 
(b_{\boldsymbol k})_{{\boldsymbol k} \in  D(n)}) \in {\mathcal B}_n$, the combination
\begin{equation*} 
\bigl(1- \eta \bigr) m*\varphi_n + \eta {\mathscr I}_n 
( {\boldsymbol a}, {\boldsymbol b})  
\end{equation*}  
is an element of ${\mathcal P}_n({\mathbb T}^d)$. Its pre-image by 
${\mathscr I}_n$ is 
\begin{equation*}
\begin{split}
&{\mathscr I}_n^{-1} \Bigl(
\bigl(1- \eta \bigr) m*\varphi_n + \eta {\mathscr I}_n 
( {\boldsymbol a}, {\boldsymbol b})  
\Bigr) 
\\
&\hspace{5pt} =
 \biggl(
 \Bigl[ 
 (1-\eta) \Re\bigl[\widehat{m}^{{\boldsymbol k}}\bigr]
  \widehat{\varphi}_n^{{\boldsymbol k}}
 +
 \eta a_{{\boldsymbol k}}
 \Bigr] 
 ,
 \Bigl[ 
 (1-\eta){\Im}\bigl[\widehat{m}^{{\boldsymbol k}}\bigr]  \widehat{\varphi}_n^{{\boldsymbol k}}
 +
 \eta b_{{\boldsymbol k}}
 \Bigr] 
  \Bigr)_{ {\boldsymbol k} \in 
 D(n)}.
 \end{split}
\end{equation*} 

\underline{Step 2}.
For $\eta \in (0,1)$, we now let:
\begin{equation*}
\Phi^{{\eta},n}(m) = \int_{{\mathbb R}^{\vert D(n) \vert} \times {\mathbb R}^{\vert D(n)\vert}}
\Phi \Bigl( 
\bigl(1- {\eta}\bigr) m*\varphi_n + {\eta}{\mathscr I}_n 
( {\boldsymbol a}, {\boldsymbol b}) 
\Bigr) \rho_n({\boldsymbol a}) \rho_n({\boldsymbol b}) d {\boldsymbol a} d{\boldsymbol b}, 
\end{equation*}
for $\rho_n$ a smooth kernel whose support is included in ${\mathcal B}_n$. 
By an obvious change of variable, we have 
\begin{equation*}
\begin{split}
\Phi^{{\eta},n}(m) &= \int_{{\mathbb R}^{\vert D(n) \vert} \times {\mathbb R}^{\vert D(n)\vert}}
\biggl[ \Phi \bigl( 
 {\mathscr I}_n 
( {\boldsymbol a}, {\boldsymbol b})   
\bigr) 
\\
&\hspace{5pt} \times
\rho_n\biggl(\Bigl(\frac{a_{\boldsymbol k} - (1-{\eta}) \Re\bigl[\widehat{m}^{\boldsymbol k}\bigr] {\widehat{\varphi}_n^{\boldsymbol k}}}{{\eta}}
\Bigr)_{{\boldsymbol k} \in D(n)} 
\biggr) \rho_n
\biggl(
\Bigl(\frac{b_{\boldsymbol k} - (1-{\eta}) \Im\bigl[\widehat{m}^{\boldsymbol k}
\bigr] {\widehat{\varphi}_n^{\boldsymbol k}}}{\epsilon}
\Bigr)_{{\boldsymbol k} \in D(n)}
\biggr)
\biggr]
d {\boldsymbol a} d{\boldsymbol b},
\end{split}  
\end{equation*}
from which we deduce that 
$\Phi^{{\eta},n}$ is infinitely differentiable with respect to the real and imaginary parts of the 
Fourier coefficients of the inputs.
In particlar, the latter says that $\Phi^{{\eta},n}$ can be regarded as a smooth function on 
${\mathcal O}_n$. The analysis achieved in 
\cite{cecchin2022weak} shows that 
$\Phi^{{\eta},n}$ is continuously differentiable in $m$. 
The derivative {$[\delta \Phi^{\eta,n}/\delta m](m,\cdot)$} has, at any $m \in {\mathcal P}({\mathbb T}^d)$, a finite number of 
non-zero Fourier modes and is thus smooth. This shows \textit{(3)} in the statement. 
\vskip 5pt

\underline{Step 3}. We now prove \textit{(2)} 
in the statement. Lipschitz property (w.r.t. $d_1$) is just a consequence of the following two facts: $\Phi$ is $c_1$-Lipschitz continuous w.r.t. $d_1$ and 
\begin{equation*}
d_1(m*\varphi_n,m'*\varphi_n) \leq d_1(m,m'), \quad m,m' \in {\mathcal P}({\mathbb T}^d), 
\end{equation*}
i.e. convolution by Féjer kernel is a contraction {under} $d_1$. The proof of the above inequality is quite obvious. 
For two $m,m' \in {\mathcal P}({\mathbb T}^d)$, we call $\pi$ an optimal coupling between 
$m,m'$ for the $d_1$-distance. Then, the 
probability measure defined by 
\begin{equation*}
\pi' (A \times B ) = \int_{{\mathbb T}^d \times {\mathbb T}^d} 
\pi \bigl( (A+z) \times (B+z) \bigr) d \varphi_n(z),
\end{equation*} 
for any two Borel subsets $A$ and $B$ of ${\mathbb T}^d$, is a coupling between $m*\varphi_n$ and
$m' * \varphi_n$. Obviously, 
\begin{equation*} 
d_1 \bigl( m*\varphi_n,m'*\varphi_n\bigr) 
\leq 
\int_{{\mathbb T}^d \times {\mathbb T}^d} 
\vert x-y \vert d \pi'(x,y) 
= 
\int_{{\mathbb T}^d \times {\mathbb T}^d} 
\vert x-y \vert d \pi(x,y) 
= d_1 (m,m'). 
\end{equation*} 

As for the proof of the semi-concavity property, we already know that, for any $\lambda \in [0,1]$,  
\begin{equation*} 
\Phi\bigl( \lambda m' + (1-\lambda) m \bigr)
\geq \lambda \Phi(m') + (1-\lambda) \Phi(m) 
- \frac{c_2}2 \lambda (1-\lambda) d_1^2(m,m'), \quad 
m,m' \in {\mathcal P}({\mathbb T}^d). 
\end{equation*}
Inserting this property in the definition of $\Phi^{{\eta},n}$ 
and using again the contractive property of the convolution under 
$d_1$, we deduce that $\Phi^{{\eta},n}$ is also $c_2$-semi-concave
w.r.t. $d_1$. 
\vskip 4pt

\underline{Step 4}.
We now prove \textit{(1)}. To do so, it suffices to observe that 
\begin{equation*} 
\begin{split}
d_1 \bigl( (1- {\eta}) m * \varphi_n + 
{\eta} {\mathscr I}({\boldsymbol a},{\boldsymbol b}) {,m} \bigr) 
&\leq (1- {\eta}) 
d_1 \bigl( m * \varphi_n, m \bigr) 
+ 
{\eta}
d_1 \bigl(  {\mathscr I}({\boldsymbol a},{\boldsymbol b}) * \varphi_n,m \bigr) 
\\
&\leq (1- {\eta}) 
d_1 \bigl( m * \varphi_n, m \bigr)
+ c_d {\eta}, 
\end{split} 
\end{equation*} 
for a constant $c_d$ only depending on the dimension (this following from the fact that, on the torus, 
$d_1$ is bounded by the total variation distance up to a multiplicative constant{)}. 

It then remains to see that 
\begin{equation*} 
\begin{split} 
d_1 \bigl( m * \varphi_n, m \bigr)
&\leq \sup_{f 1-\text{Lip}} \int_{{\mathbb T}^d} f {(x)}d \bigl( m*\varphi_n - m \bigr) {(x)}
\leq \sup_{f 1-\text{Lip}} \bigl\| f - f*\varphi_n \|_\infty. 
\end{split}
\end{equation*} 
The latter right-hand side tends to $0$ as $n$ tends to $\infty$, see again \cite{cecchin2022weak}. 
\end{proof}

\color{black}

\bibliographystyle{alpha}

\bibliography{/Users/sam/Documents/Bibtex/DDJ05152023Arxiv.bib}

\end{document}